\documentclass[12pt,reqno]{amsart}
\usepackage{amsfonts,amsmath,amssymb,amsthm}
\usepackage{latexsym}
\usepackage{eucal}
\usepackage[ansinew]{inputenc}
\usepackage[american]{babel}
\usepackage{amsfonts}
\usepackage{amssymb}
\usepackage[dvips]{graphicx}

\numberwithin{equation}{section}

\newtheorem{teo}{Theorem}[section]
\newtheorem{lema}{Lemma}[section]
\newtheorem{prop}{Proposition}[section]
\newtheorem{defi}{Definition}[section]

\newtheorem{coro}{Corollary}[section]

\begin{document}

\title[Stability of Periodic-Peak Travelling Waves Solutions]
{The Non-Linear  Schr\"odinger Equation with a periodic {\bf{$\delta$}}--interaction}

\author[J. Angulo]{Jaime Angulo Pava $^\dagger$}
\author[G. Ponce]{Gustavo Ponce $^\ddagger$}

 \email{angulo@ime.usp.br} \email{ponce@math.ucsb.edu}

\subjclass[2000]{76B25, 35Q51, 35Q53.}

\keywords{NLS-Dirac equation,
periodic travelling-waves, nonlinear stability}
\thanks{{\it Date}: {\bf {10/25/2010.}}}

\maketitle

{\scriptsize \centerline{$^\dagger$Department of Mathematics,
IME-USP}
 \centerline{Rua do Mat\~ao 1010, Cidade Universit\'aria, CEP 05508-090, S\~ao Paulo, SP, Brazil.}
\centerline{$^\ddagger$ Department of Mathematics, UCSB}\centerline{Santa Barbara, CA. 93106, USA.}
 }

\begin{abstract}
We study the existence and stability of the standing waves for the periodic cubic nonlinear 
Schr\"odinger equation with a point defect determined by a periodic Dirac distribution at the origin. This equation admits a smooth curve of positive periodic solutions in the form of standing waves with a profile given by the Jacobi elliptic function of dnoidal type.  Via a perturbation method and continuation argument, we obtain that in the case of an attractive defect the standing wave solutions are stable in $H^1_{per}$ with respect to perturbations which have the same period as the wave itself. In the case of a repulsive defect, the standing wave solutions are stable in the subspace of even functions of $H^1_{per}$ and unstable in $H^1_{per}$ with respect to perturbations which have the same period as the wave itself. 

\end{abstract}


\section{Introduction}

Consider the semi-linear Schr\"odinger equation (NLS)  
\begin{equation}\label{nlsp}
\partial_t u+\Delta u \pm |u|^pu=0,\;\;\;\;\;\;\;(x,t)\in \mathbb R^n\times \mathbb R,
\end{equation}
where $u=u(x,t)$ is a complex-valued function and $0<p<\infty$. 
This is a canonical dispersive equation which arises as a model in several physical situations, see for example \cite{SuSu},
\cite{CMM}, and references therein.

The mathematical study of the NLS (the local well posedness of its initial value problem (IVP) and its periodic boundary value problem (PBVP) under minimal regularity assumptions on the data, the long time behavior of their solutions, blow up and scattering results, etc)  has attracted a great deal of attention and is a very active research area  (see \cite{Ca}, \cite{Bou}, \cite{Tao}, and \cite{LP}). 

In the  one dimensional cubic case, $n=1,\;p=2$,  it was established  in \cite{ZaSh} that the NLS is a completely integrable system. Thus, using the inverse scattering theory it can be solved in the line $\mathbb R$ (IVP) and in the circle $\mathbb T$ (PBVP) (see  \cite{AS}, \cite{MA} and references therein).

Special solutions of the NLS equation \eqref{nlsp} have been widely considered in analytic, numerical and experimental works. In particular, in the focussing case ($+$ in \eqref{nlsp})  one has the \lq\lq standing waves" solutions

\begin{equation}\label{stand}
u_s(x,t)=e^{i\omega t} \phi(x),\;\;\;\;\;\;\;\;\omega>0,
\end{equation}
or their generalization \lq\lq travelling waves" solutions
\begin{equation}\label{standco}
u_{tw}(x,t)=e^{i\omega t}\,e^{i(c\cdot x- |c|^2 t)} \phi(x-2ct),\;\;\;\;\;\;\;\;\;\omega>0,\;\;c\in \mathbb R^n,
\end{equation}
with $\phi=\phi_{\omega, p}$ being the unique positive, radially symmetric solution (ground state) of the nonlinear elliptic 
problem
\begin{equation}\label{stand1}
-\Delta \phi+\omega\phi(x)-\phi^{p+1}(x)=0,\quad x\in \mathbb R^n,
\end{equation}
satisfying  the boundary condition $\phi(x)\to 0$ as $|x|\to \infty$. In the one dimensional case, $n=1$, $\phi$ is given by the explicit formula (modulo translation)
\begin{equation}\label{stand2}
\phi(x)=\phi_{\omega,p}(x)=\Big[\frac{(p+2)\omega}{2} sech^2\Big(\frac{p\sqrt{\omega}}{2}x\Big)\Big]^{\frac{1}{p}}.
\end{equation}

The stability and instability properties of the standing waves have been extensively studied.  A crucial role in the stability analysis is played by the symmetries of the NLS equation in $\mathbb R^n$. The most important ones for this purpose are :
\begin{enumerate}
\item {\it phase invariance}: $u(x,t)\to e^{i\theta}u(x,t),\;\theta\in \mathbb R$;

\item {\it translation invariance}: $u(x,t)\to u(x+y,t), \;y\in \mathbb R^n$;
 
\item  {\it Galilean invariance}: $u(x,t)\to e^{i(v \cdot x-|v|^2 t)}u(x-2vt,t),\;v\in \mathbb R^n$.
\end{enumerate}
So,  if one  considers the orbit generated by the solution $\phi=\phi_{\omega,p}$ of \eqref{stand1} and the
phase-invariance symmetries above, namely,
\begin{equation}
\label{stand3}
\Theta(\phi_{\omega,p})=\{e^{i\theta}\phi_{\omega,p}(\cdot+y): \theta\in [0, 2\pi), y\in \mathbb R^n\},
\end{equation}
is known that in the one dimensional case, $n=1$,   $\Theta(\phi_{\omega,p})$  is stable in $H^1(\mathbb R)$ by the flow of the NLS equation provided that $p<4$ and unstable for $p\geqq 4$ (for details and results in higher dimensions see Cazenave\&Lions \cite{CL}, Weinstein \cite{w3}). This means that for $p<4$, if $u_0$ is close to $\Theta(\phi_{\omega,p})$ in $H^1(\mathbb R^n)$, then the corresponding solution of \eqref{nlsp} $u(t)$ with initial data  $u_0$ remains close to the orbit $\Theta(\phi_{\omega,p})$ for each $ t\in\mathbb R$. The necessity of the rotations and space translation appearing in the stability criterium  can be seen in \cite{Ca}.

 From now on we shall restrict our attention to the one dimensional focussing  NLS
 \begin{equation}
 \label{nls1}
 i\partial_t u +\partial_x^2 u+|u|^pu=0,\;\;\;\;\;p>0.
\end{equation}

In contrast to the standing waves solutions in the line, i.e. \eqref{stand} and \eqref{standco} with  $n=1$ and $\phi$ as in  \eqref{stand2}, relatively little is known about the existence and stability of  periodic standing wave solutions, i.e., $\phi$ in \eqref{stand} being a periodic function.

A partial {\it spectral stability} analysis was carried out by Rowlands \cite{Row} for the case $p=2$ with respect to long-wave disturbances, who showed that periodic waves with real-valued profile are unstable. Similar results were also obtained for certain NLS-type equations with spatially periodic potentials by Bronski\&Rapti  \cite{BR}. The first results concerning the {\it nonlinear stability} of periodic standing waves are due to Angulo  \cite{angulo1}.  In \cite{angulo1} he  established   the existence of a smooth family of {\it dnoidal waves} for the cubic NLS equation ($p=2$ in \eqref{nls1}) of the form
\begin{equation}\label{dnoidal1}
\omega\in\Big(\frac{\pi^2}{2L^2}, +\infty\Big)\to \phi_{\omega, 0}\in H^{\infty}_{per}([-L, L]),
\end{equation}
where the profile of $\phi=\phi_{\omega}=\phi_{\omega,0}$ is given by the Jacobian elliptic function called {\it dnoidal}, $dn$ by the formula
\begin{equation}\label{dnoidal2}
\phi_{\omega}(\xi)=\eta_1dn\Big(\frac{\eta_1}{\sqrt{2}} \xi;k\Big),
\end{equation}
with $\eta_1\in (\sqrt{\omega}, \sqrt{2\omega})$ and the modulus $k\in (0,1)$ depending smoothly on $\omega$. Angulo showed  that for every $\omega>\frac{\pi^2}{2L^2}$ the $2L$-periodic wave $\phi_{\omega}$ is {\it orbitally stable} with respect to perturbations which have the same period as the wave itself, and  {\it nonlinearly unstable} with respect to perturbations which have two times the period ($4L$) as the wave itself. Indeed, the  same analysis used to obtain the instability result  provides the {\it nonlinear instability} of the dnoidal wave by perturbations which have $j$-times ($j> 2$) the period as the wave itself (for further details see also \cite{angulo1} and \cite{angulo4}). 

In \cite{GH1}-\cite{GH2} Gallay\&Haragus have shown the stability of periodic traveling waves described in \eqref{standco} for the cubic NLS equation by allowing the profile $\phi$ being complex-valued. In the case $p=4$, Angulo\&Natali \cite{AnguloNatali2} have shown the existence of a family of periodic waves of the form described in \eqref{stand} for which there is a unique (threshold)  value of the phase-velocity $\omega$ which separates the two global scenarios: stability and instability.  

In this paper we are interested in  the periodic setting  for nonlinear Schr\"odinger equation (NLS-$\delta$ henceforth) of the form
\begin{equation}\label{nls}
i\partial_t u+\partial_x^2 u+Z \delta(x) u +|u|^pu=0,
\end{equation}
where $\delta$ is the Dirac distribution at the origin, namely, $\langle \delta, v\rangle=v(0)$ for $v\in H^1$, and $Z\in \mathbb R$. The equation \eqref{nls}, $Z\neq 0$ has  been  considered in a variety of physical models with a point defect, for instance, in nonlinear optics and Bose-Einstein condensates. Indeed, the Dirac distribution is used to model an impurity, or defect, localized at the origin. Also in this case the NLS-$\delta$ equation \eqref{nls} can be viewed as a prototype model for the interaction of a wide soliton with a highly localized potential. In nonlinear optics, this models a soliton propagating in a medium with a point defect or the interaction of a wide soliton with a much narrower one in a bimodal fiber, see  \cite{ghw},  \cite{SCH},  \cite{CM}, \cite{MN}, \cite{Mn}, \cite{Ag}, \cite{BK}, \cite{DMADDKK}, \cite{ST}, and the reference therein.

Equation \eqref{nls} in the line with $p=2$ has been  considered by several authors. In a series of papers  \cite{HMZ1}, \cite{HMZ2}, \cite{HZ1}, and \cite{HZ2} the phenomenon of soliton scattering by the effect of the defect was comprehensibly studied. In particular, in \cite{HZ1} for  the equation \eqref{nls} with $p=2$ and data
\begin{equation}\label{sol}
u(x,0)=e^{icx}sech(x-x_0),\;\;\;x_0<<-1,
\end{equation}
it was shown that for the $|Z|<<1$ the corresponding solution, the traveling wave 
for $t>|x_0|/c$ remains intact. The case $Z>0$ and $|c|>>1$ was examinated in \cite{HMZ1}, \cite{HMZ2} where it was proven  how the defect separate  the soliton into two parts: one part is transmitted past the defect, the other one is reflected at the defect.  The case $Z<0$ and $|c|>>1$ was considered in \cite{DH}.

The existence of standing wave solutions   of  the equation \eqref{nls} requires that the profile $\phi=\phi_{\omega,Z,p}$ satisfy the semi-linear elliptic equation
\begin{equation}\label{ellip}
-\phi''(x)+\omega\phi(x)-Z\delta(x)\phi-|\phi(x)|^{p}\phi(x)=0,\quad x\in \mathbb R.
\end{equation}
In Fukuizumi\&Jeanjean \cite{fuje} (see also \cite{ghw}) it was deduced the formula for the unique positive even solution of \eqref{ellip}, modulo rotations :
\begin{equation}\label{ellip1}
\phi_{\omega, Z, p}(x)=\Big[\frac{(p+2)\omega}{2} sech^2\Big(\frac{p\sqrt{\omega}}{2}|x|+\tanh^{-1}\Big(\frac{Z}{2\sqrt{\omega}}\Big)\Big)\Big]^{\frac{1}{p}},\;\;\;x\in\mathbb R,
\end{equation}
if $\omega>Z^2/4$.  This solution is constructed from the known solution in the case  $Z=0$ on each side of the defect pasted together at $x=0$ to satisfy the conditions of continuity and the jump condition in the first derivative at $x=0$, $u'(0+)-u'(0-)=-Zu(0)$. So  $\phi$ belongs to the domain  of the formal expression $-\partial_x^2-Z\delta$ (see \cite{agfhr}) 
$$
\{u\in H^1(\mathbb R)\cap H^2(\mathbb R-\{0\}): u'(0+)-u'(0-)=-Zu(0)\}.
$$
Notice that there is no nontrivial solution of \eqref{ellip} in $H^1(\mathbb R)$ when $\omega\leqq Z^2/4$.

The basic symmetry associated to equation \eqref{nls} is the phase-invariance  since the translation invariance of the solutions is not hold due to the defect. Thus,  the notion of stability and instability will be based only on this symmetry and is formulated as follows:

\begin{defi}\label{dsta} For $\eta>0$, let $\phi$ be a solution of \eqref{ellip} and define
$$
U_\eta(\phi)=\Big\{v\in X: \inf_{\theta\in\mathbb R}\|v-e^{i\theta}\phi\|_X<\eta\Big\}.
$$
The standing wave $e^{i\omega t}\phi$ is (orbitally) stable in $X$ if for any $\epsilon>0$ there exists $\eta>0$ such that for any $u_0\in U_\eta(\phi)$, the solution $u(t)$ of \eqref{nls} with $u(0)=u_0$ satisfies $u(t)\in U_\epsilon(\phi)$ for all $t\in \mathbb R$. Otherwise, $e^{i\omega t}\phi$ is said to be (orbitally) unstable in $X$.
\end{defi}

Gathering  the information in \cite{fuje}, \cite{fuoht}, \cite{ghw}, and \cite{lffgs}, one can summarize the known results on the stability and instability  of standing waves associated to the solitary wave-peak  in \eqref{ellip1} as follows:
\begin{enumerate}
\item[$\bullet$] Let $Z>0$ and $\omega>Z^2/4$.
\begin{enumerate}
\item If $0<p\leqq 4$, the standing wave $e^{i\omega t}\phi_{\omega, Z, p}$ is stable in $H^1(\mathbb R)$ for any $\omega\in (Z^2/4, +\infty)$.

\item If $p\geqq 5$, there exists a unique $\omega_1>Z^2/4$ such that $e^{i\omega t}\phi_{\omega, Z, p}$ is stable in $H^1(\mathbb R)$ for any $\omega\in (Z^2/4, \omega_1)$, and unstable in $H^1(\mathbb R)$ for any $\omega\in (\omega_1, +\infty)$.
\end{enumerate}

\item[$\bullet$] Let $Z<0$ and $\omega>Z^2/4$.
\begin{enumerate}
\item If $0<p\leqq 2$, the standing wave $e^{i\omega t}\phi_{\omega, Z, p}$ is stable in $H^1_{rad}(\mathbb R)$ for any $\omega\in (Z^2/4, +\infty)$.

\item If $0<p\leqq 2$, the standing wave $e^{i\omega t}\phi_{\omega, Z, p}$ is unstable in $H^1(\mathbb R)$ for any $\omega\in (Z^2/4, +\infty)$.

\item If $2<p< 4$, there exists a  $\omega_2>Z^2/4$ such that $e^{i\omega t}\phi_{\omega, Z, p}$ is unstable in $H^1(\mathbb R)$ for any $\omega\in (Z^2/4, \omega_2)$, and stable in $H^1_{rad}(\mathbb R)$ for any $\omega\in (\omega_2, +\infty)$.

\item If $2<p< 4$,  the standing wave $e^{i\omega t}\phi_{\omega, Z, p}$ is unstable in $H^1(\mathbb R)$ for any  $\omega\in (\omega_2, +\infty)$, where $\omega_2$ is that in item (c) above.

\item if $p\geqq 4$, then the standing wave $e^{i\omega t}\phi_{\omega, Z, p}$ is unstable in $H^1(\mathbb R)$.
\end{enumerate}
\end{enumerate}

In this paper, we study the existence and nonlinear stability of periodic standing waves solutions of \eqref{nls} in the case $p=2$ and $Z\neq 0$. More precisely, we show the existence of a branch of periodic solutions, $\omega\to \varphi_{\omega, Z}$, for the semilinear elliptic equation
\begin{equation}\label{p2}
-\varphi_{\omega, Z}''+\omega \varphi_{\omega, Z}-Z \delta(x)\varphi_{\omega, Z}=\varphi_{\omega, Z}^3,
\end{equation}
where $\varphi_{\omega, Z}>0$ is a periodic real-valued function with prescribe period $2L>0$ and where $\omega$ will belong to a determined interval in $\mathbb R$ with $\omega>Z^2/4$. Our solutions $\varphi=\varphi_{\omega,Z}$ satisfy the following boundary values:
\begin{equation}\label{p3}
\begin{array}{llll}

(1) \; \varphi_{\omega, Z} (x+2L)=\varphi_{\omega, Z} (x),\quad {\text{for all}}\;x\in \mathbb R.\\

(2)\; \varphi_{\omega, Z} \in C^j(\mathbb R-\{2nL: n\in\mathbb Z\})\cap C(\mathbb R),\quad j=1,2.\\

(3)\; -\varphi_{\omega, Z}''(x)+\omega \varphi_{\omega,Z}(x)=\varphi_{\omega, Z}^3(x) \quad {\rm{for}}\;x\neq \pm 2nL,\;n\in \mathbb N.\\

(4)\; \varphi'_{\omega, Z}(0+)-\varphi'_{\omega, Z}(0-)=-Z \varphi_{\omega,Z}(0).\\
\end{array}
\end{equation}
The notation $\varphi'_{\omega, Z}(0\pm)$ in \eqref{p3} is defined as $\varphi'_{\omega, Z}(0\pm)=\lim_{\epsilon\downarrow 0}\varphi'_{\omega, Z}(\pm\epsilon)$. From the periodicity of the function $\varphi_{\omega, Z}$  one also has that $ \varphi'_{\omega,Z}(\pm 2nL+)-\varphi'_{\omega, Z}(\pm 2nL-)=-Z \varphi_{\omega,Z}(2nL)$, for $n\in \mathbb N$. We recall that if $\varphi_{\omega,Z}$ is a solution of \eqref{p2} then $\varphi_{\omega, Z}(\cdot+y)$ is not necessarily a solution of  \eqref{p2}.  Hence, our stability study for the ``{\it periodic-peaks}'' $\varphi_{\omega, Z}$ will be for the orbit generated by this solution and defined in the form
\begin{equation}\label{orb}
\Omega_{\varphi_{\omega,Z}}=\{e^{i\theta}\varphi_{\omega,Z}: \theta\in [0,2\pi]\}.
\end{equation}

From equation \eqref{p2} arises naturally the condition that our solutions $\varphi_{\omega,Z}$ need to belong to the domain of the formal expression 
\begin{equation}\label{Zdelta}
-\frac{d^2}{dx^2}-Z\delta.
\end{equation}
So, we shall develop a precise formulation for this {\it periodic} point interaction, also called $\delta$-interaction. We present a detailed study of the model of quantum mechanics \eqref{Zdelta} with a potential supported on a $\delta$ and in the framework of periodic functions.  In our  study of the \lq\lq solvability'' of this model we will describe their resolvents  explicitly in terms of the interactions strengths, $Z$, and the location of the source, $x=0$. We start by establishing the  definition of all the self-adjoint extensions of the operator $A^0=-\frac{d^2}{dx^2}$ with domain
\begin{equation}\label{specA}
D(A^0)=\{\psi\in D(A) : \delta(\psi)\equiv\psi(0)=0\},
\end{equation}
which is a densely defined symmetric operator  on $L^2_{per}([0,2L])$ with deficiency indices $(1,1)$. Here $A$ represents the self-adjoint operator $-\frac{d^2}{dx^2}$ on $L^2_{per}([0,2L])$ with the natural domain $D(A)=H^2_{per}([0,2L])$.  Using the von Neumann theory we can parametrized all the self-adjoint extensions of $A^0$ with the help of $Z$. Indeed, for $Z\in [-\infty, \infty)$ we have 
\begin{equation}\label{specA1}
\begin{cases}
\begin{aligned}
&-\Delta_{-Z}=-\frac{d^2}{dx^2}\\
&D(-\Delta_{-Z})=\{\zeta\in H^1_{\text{per}}([-L, L])\cap H^2((-L,L)-\{0\})\cap H^2((0, 2L)): \\
&\qquad\qquad\qquad\qquad \zeta'(0+)-\zeta'(0-)=-Z \zeta(0)\}.
\end{aligned}
\end{cases}
\end{equation}
These definitions are not only important to determine solutions for equation in \eqref{p2} but also for our nonlinear stability theory.

In Section 5, we will find a smooth branch of positive, even, periodic-peak solutions of \eqref{p2}, $\omega\to \phi_{\omega,Z}\in H^n_{per}([0,2L])$, such that $\phi_{\omega,Z}$ belongs to the domain of the formal expression $-\frac{d^2}{dx^2}-Z\delta$ and satisfying
\begin{equation}\label{specA2}
\lim_{Z\to 0^+} \phi_{\omega,Z}=\phi_{\omega,0}
\end{equation}
where $\phi_{\omega,0}$ is the dnoidal traveling wave defined in \eqref{dnoidal2}. The profile of $\phi_{\omega,Z}$ is based in the Jacobian elliptic function {\it dnoidal} and determined for $\omega>Z^2/4$  by the pattern
\begin{equation}\label{specA3}
\phi_{\omega,Z}(\xi)=\eta_{1,Z} dn\Big(\frac{\eta_{1,Z}}{\sqrt{2}}|\xi|+a;k\Big), \quad \xi\in [-L,L]
\end{equation}
where $\eta_{1,Z}$ and the modulus $k$ depend smoothly of $\omega$ and $Z$. The shift value $  \,a\, $ is also a smooth function of $\omega$ and $Z$ satisfies  that $\lim_{Z\to 0^+}a(\omega,Z)=0$. See Figure 3 below for a general profile of $\phi_{\omega,Z}$.

Similarly, we obtain via the theory of elliptic integrals for $Z<0$ a   smooth branch of positive, even, periodic-peak solutions of \eqref{p2}, $\omega\to \zeta_{\omega,Z}\in H^n_{per}([0,2L])$, such that $\zeta_{\omega,Z}$ belongs to the domain of the formal expression $-\frac{d^2}{dx^2}-Z\delta$ and satisfying
\begin{equation}\label{specA4}
\lim_{Z\to 0^-} \zeta_{\omega,Z}=\phi_{\omega,0}
\end{equation}
where $\phi_{\omega,0}$ is the dnoidal wave defined in \eqref{dnoidal2}. The profile of $\zeta_{\omega,Z}$ is determined for $\omega>Z^2/4$  by the pattern
\begin{equation}\label{specA5}
\zeta_{\omega,Z}(\xi)=\eta_{1,Z} dn\Big(\frac{\eta_{1,Z}}{\sqrt{2}}|\xi|-a;k\Big), \quad \xi\in [-L,L].
\end{equation}
See Figure 4 below for a general profile of $\zeta_{\omega,Z}$. We note that the periodic-peak $\zeta_{\omega,Z}$ and $\phi_{\omega,Z}$  ``converge'' to the solitary wave-peak $\phi_{\omega, Z,2}$ in \eqref{ellip1} when we consider $\eta_1\to \sqrt{2\omega}$. We refer the reader to Section 5 for the precise details on this convergence.

Our approach for the stability theory of the periodic-peak family
\begin{equation}
\label{fami}
\varphi_{\omega,Z}=
\begin{cases}
\begin{aligned}
&\phi_{\omega, Z},\quad  Z>0,\\
&\zeta_{\omega, Z},\quad  Z<0,
\end{aligned}
\end{cases}
\end{equation}
with $\phi_{\omega,Z}$ and $\zeta_{\omega,Z}$ given in \eqref{specA3}-\eqref{specA4}, it will be based in the general framework developed by Grillakis\&Shatah\&Strauss \cite{grillakis1}, \cite{grillakis2}, for a Hamiltonian system which is invariant under a one-parameter unitary group of operators. This theory  requires the following informations :
\begin{enumerate}
\item[$\bullet$] The {\it Cauchy problem}: The initial value problem associated to the NLS-$\delta$ equation is well-posedness in $H^1_{per}([0,2L])$.

\item[$\bullet$] The {\it spectral condition}:
\begin{enumerate}
\item The self-adjoint operator $\mathcal L_{2, Z}$ defined on $L^2_{per}([0,2L])$, as
\begin{equation}\label{specA6}
\mathcal L_{2, Z}\zeta =-\frac{d^2}{dx^2}\zeta+\omega \zeta-\varphi^2_{\omega, Z}\zeta
\end{equation}
with domain $\mathcal D= D(-\Delta_{-Z})$ given in \eqref{specA1},  is a nonnegative operator with the eigenvalue zero being simple and with eigenfunction $\varphi_{\omega, Z}$.

\item The self-adjoint operator $\mathcal L_{1, Z}$ defined on $L^2_{per}([0,2L])$, by
\begin{equation}\label{specA7}
\mathcal L_{1, Z}\zeta =-\frac{d^2}{dx^2}\zeta+\omega \zeta-3\varphi^2_{\omega, Z}\zeta
\end{equation}
with domain $\mathcal D= D(-\Delta_{-Z})$ given in \eqref{specA1},  has a trivial kernel for all $Z\in \mathbb R-\{0\}$.

\item The number of negative eigenvalues of the operator $\mathcal L_{1, Z}$.
\end{enumerate}

\item[$\bullet$] The {\it slope condition}: The sign of $\partial_\omega\int_{-L}^L \varphi^2_{\omega, Z}(\xi)d\xi$.

\end{enumerate}

In general, to count the number of negative eigenvalues of linear operator is a delicate issue. In the case of the self-adjoint operator $\mathcal L_{1, Z}$ our strategy is based in two basic facts. The first one is that in the case $Z=0$, the spectrum of the self-adjoint operator $\mathcal L_0\equiv \mathcal L_{1, 0}$ defined on  $L^2_{per}([0,2L])$  by
\begin{equation}\label{specA8}
\mathcal L_{0}\zeta =-\frac{d^2}{dx^2}\zeta+\omega \zeta-3\phi_{\omega,0}^2\zeta
\end{equation}
with domain $H^2_{per}([0,2L])$ and $\omega>\pi^2/{2L^2}$, has already been described  in \cite{angulo1} and in \cite{AnguloNatali1}: there is only one negative eigenvalue which is simple, zero is a simple eigenvalue  with eigenfunction $\frac{d}{dx}\phi_{\omega,0}$. The rest of the spectrum is positive and discrete. The second  is that for $Z$ small, $\mathcal L_{1, Z}$ can be considered as a {\it real-holomorphic perturbation} of $\mathcal L_{0}$. So, we have that the spectrum of $\mathcal L_{1, Z}$ depends holomorphically on the spectrum of  $\mathcal L_{0}$. Then we obtain that for $Z<0$ there are exactly two negative eigenvalues of $\mathcal L_{1, Z}$ and exactly one for $Z>0$. We refer the reader to Subsection 6.1 for the precise details  on these statements.

Our main result  is the following:

\begin{teo}\label{main} Let $\omega>\frac{Z^2}{4}$ and $\omega>\frac{\pi^2}{2L^2}$. We have for $\omega$ large:
\begin{enumerate}
\item For $Z>0$, the  dnoidal-peak standing wave $e^{i\omega t}\varphi_{\omega,Z}$ is stable in $H^1_{per}([-L,L])$.

\item For $Z<0$, the  dnoidal-peak standing wave  $e^{i\omega t}\varphi_{\omega,Z}$ is unstable in $H^1_{per}([-L,L])$.

\item For $Z<0$, the  dnoidal-peak standing wave  $e^{i\omega t}\varphi_{\omega,Z}$ is stable in $H^1_{per, even}([-L,L])$.
\end{enumerate}
\end{teo}

The restriction about $\omega$ being large in Theorem \ref{main} is due to technical reasons in proving  the strictly increasing property of the mapping $\omega\to\|\varphi_{\omega,Z}\|^2$ (see Theorem \ref{con1} in Section 6.2).

The local well-posedness of the Cauchy problem for \eqref{nls} with $p=2$ in $H^1_{per}([0, 2L])$ is an consequence from Theorem 3.7.1 in \cite{Ca} and the theory spectral established in Section 3 for the  operator $-\partial_x^2-Z\delta$ for $Z\neq 0$.  The global existence of solutions is an immediate consequence of the following conserved quantities for \eqref{nls}: the  energy and the charge, respectively,
\begin{equation}\label{conse}
\begin{aligned}
E(v)&=\frac12 \int |v'(x)|^2\;dx-\frac{Z}{2}\int\delta(x)|v(x)|^2dx-\frac{1}{4}
\int |v(x)|^{4}\;dx,\\
Q(v)&=\frac12\int |v(x)|^{2}\;dx.
\end{aligned}
\end{equation}

This paper is organized as follows. Section 3 is devoted to establish a spectral theory for the operator $-\partial_x^2-Z\delta$ for $Z\neq 0$. Our analysis is based in the theory of von Neumann for self-adjoint extensions. Section 4 is concerned with the periodic well-posedness theory for \eqref{nls}, $p=2$, in $H^1_{per}([0, 2L])$. Section 5 describe the construction, via the implicit function theorem, of a smooth curve of periodic-peak for equation \eqref{p2}. Finally, in Section 6, the stability and instability theory of the dnoidal-peak is established.

\section{Notation} For any complex number $z\in\mathbb C$, we denote by $\Re\, z$ and $\Im\, z$ the real part and imaginary part of $z$, respectively. For $s\in\mathbb R$, the Sobolev space $H^s_{\rm per}([0,2L])$ consists of all periodic distributions $f$ such that $\|f\|_{H^s}^2= 2L\underset{k=-\infty}{\overset{\infty}{\sum}}(1+k^2)^s|\widehat{f}(k)|^2<\infty$. For simplicity, we will use the notation $H^s_{\rm per}$ in several places and  $H^0_{\rm per}=L^2_{\rm per}$. We remark that $L^2_{\rm per}$ and $H^1_{\rm per}$ are regarded as real Hilbert space with inner products
\begin{equation}\label{inner}
\langle f,g\rangle_{L^2}=\Re \int_{-L}^{L}f(x)\overline{g(x)}dx,\;\;\langle f,g\rangle_{H^1}=\langle f,g\rangle_{L^2}+\langle \partial_x f, \partial_x g\rangle_{L^2}.
\end{equation}
We denote $\|f\|_{L^2}=\|f\|$ and $\langle f,g\rangle_{L^2}=\langle f,g\rangle$. For $\Omega$ being an open set of  $ \mathbb R$, $H^n(\Omega)$, $n\geqq 1$, represents the classical local Sobolev space. $[H^s_{per}]'$, 
the topological dual of $H^s_{per}$, is isometrically isomorphic to $H^{-s}_{per}$ for all $s \in \mathbb{R}$. The duality is implemented concretely by the pairing
$$
(f,g) = 2L\sum_{k=-\infty}^{\infty}\widehat{f}(k)\overline{\widehat{g}(k)}, \ \ \ for \ \ \ f \in H^{-s}_{per}, \ \  g \in H^s_{per}. 
$$
Thus, if $f \in L^2_{per}$ and $g \in H^s_{per} $, with $s\geq 0,$ it follows that $(f,g) = \langle f,g\rangle$. The convolution for $f, g\in  L^2_{per}$ is defined by
$$
f\star g(x)=\frac{1}{2L}\int_{-L}^L f(x-y)g(y)dy.
$$
The normal elliptic integral of first type (see \cite{byrd}) is defined by
\begin{equation*}
\int\limits_0^y\frac{dt}{\sqrt{(1-t^2)(1-k^2t^2)}}=
\int\limits_0^{\varphi}\frac{d\theta}{\sqrt{1-k^2\sin^2\theta}}=F(\varphi,k)
\end{equation*}
where $y=\sin \varphi$ and  $k\in (0,1)$. $k$ is called the modulus and $\varphi$ the argument. When $y=1$, we denote $F(\pi/2,k)$ by $K=K(k)$. The three basic Jacobian elliptic functions are denoted by $sn(u;k)$, $cn(u;k)$
and $dn(u;k)$ (called, snoidal, cnoidal and dnoidal, respectively), and are defined via the previous
elliptic integral. More precisely, let
\begin{equation}\label{jef3}
u(y;k):=u=F(\varphi,k)
\end{equation}
then $y=sin \varphi:=sn(u;k)=sn(u)$ and
\begin{equation}\label{jef4}
\begin{aligned}
cn(u;k)&:=\sqrt{1-y^2}=\sqrt{1-sn^2(u;k)}\\
dn(u;k)&:=\sqrt{1-k^2y^2}=\sqrt{1-k^2sn^2(u;k)}.
\end{aligned}
\end{equation}
In particular,  we have that 
$1\geqq dn(u;k)\geqq k'\equiv \sqrt{1-k^2}$ and the following asymptotic formulas: 
$sn(x;1)=tanh (x)$, $cn(x;1)=sech(x)$ and $dn(x;1)=sech(x)$. 

Finally,  $\varphi(0\pm)=\lim_{\epsilon\downarrow 0}\varphi(\pm\epsilon)$.


\section{The one-center periodic $\delta$-interaction in one dimension}

In this section we develop a precise formulation for the {\it periodic} point interaction determined by the formal linear differential operator
\begin{equation}\label{deltaop}
-\frac{d^2}{dx^2}+\gamma \delta\equiv -\frac{d^2}{dx^2}+\gamma (\delta,\cdot)\delta ,
\end{equation}
defined on functions on the torus $\mathbb T=\mathbb R/{2\pi}$.  $\gamma$ is denominated the coupling constant or strength attached to the point source located at $x=0$. 

Our main purpose here is to study the \lq\lq solvability" of  this model. So, we will show that their resolvents can be given explicitly in terms of the interactions strengths, $\gamma$, and the location of the specific source, $x=0$. As a consequence the spectrum and  the eigenfunctions can be determined explicitly. Our method is  based on the concept of self-adjoint operator extensions of densely defined symmetric operators and so the von Neumann extension theory will be our main tool. 

The basic idea  behind the study of  models as in \eqref{deltaop} is that, once  their hamiltonian have been well defined, they can serve as corner stones for more complicated and more realistic interactions, obtained by various perturbations/approximations, such that as the point interaction models \eqref{nls}. In our case, such theory is essential for finding the  right profile of the solutions for equation in \eqref{ellip} and for the domain of the self-adjoint operators $\mathcal L_{1,Z}, \mathcal L_{2,Z}$ in \eqref{specA6}-\eqref{specA7}, which are the core of our stability theory.

For  $A^0$ being a densely defined symmetric operator on a Hilbert space  and ${A^0}^*$ denoting its adjoint, we consider the subspaces 
\begin{equation}\label{sa1}
\mathcal D_+=\text{Ker}({A^0}^*-i),\quad{\rm{and}}\quad\mathcal D_-=\text{Ker}({A^0}^*+i),
\end{equation}
$\mathcal D_+$ and $\mathcal D_-$ are called the {\it deficiency subspaces} of $A^0$. The pair of numbers $n_+$,  $n_-$, given by  
$$
n_+(A^0)=\text{dim}[\mathcal D_+], \quad{\rm{and}}\quad n_-(A^0)=\text{dim}[\mathcal D_-]
$$ 
are called the {{\it deficiency indices}  of $A^0$. 

Let $A=-\frac{d^2}{dx^2}$ and  we consider the periodic Sobolev spaces on $[0,2\pi]$, $H^s_{\text{per}}\equiv H^s_{\text{per}}([0,2\pi])$. 

\begin{lema} $A$ is a self-adjoint operator on $L^2_{\text{per}}([0,2\pi])$ with the domain $D(A)=H^2_{\text{per}}$.
 \end{lema}

Next, since $\delta\in H^{-2}_{\text{per}}- L^{2}_{\text{per}}$ we have the following.

\begin{lema} \label{sp0} The restriction $A^0\equiv A|_{D(A^0)}$, where
\begin{equation}\label{sp1}
D(A^0)=\{\psi\in D(A) : (\delta, \psi)\equiv\psi(0)=0\},
\end{equation}
is a densely defined symmetric operator with deficiency indices $(1,1)$. Namely,
\begin{enumerate}
\item symmetric: $\langle A^0\psi,\varphi\rangle =\langle \psi, A^0\varphi\rangle$ for $\psi, \varphi\in  D(A^0)$;

\item dense: $\overline{D(A^0)}=L^{2}_{\text{per}}$;

\item deficiency elements: 
\begin{equation}\label{sp2}
\begin{cases}
\begin{aligned}
&{\rm{for}}\;\;\lambda=i, \;\;\;\quad g_i\equiv (A-i)^{-1}\delta,\\
&{\rm{for}}\;\;\lambda=-i, \quad g_{-i}\equiv (A+i)^{-1}\delta,
\end{aligned}
\end{cases}
\end{equation}
$g_{\pm i}\in D({A^0}^*)$ and  ${A^0}^*g_{\pm i}=\pm ig_{\pm i}$. Moreover, $n_+(A^0)_{-}(A^0)=1$.
\end{enumerate}
\end{lema}

\begin{proof} 
\begin{enumerate}

\item The symmetric property of $A^0$ follows immediately from that of the operator  $A$.

\item The operator $A$ is densely defined and thus for every $f\in L^2_{per}$  there exists $\{f_n\}\subset H^2_{per}$ such that $\lim_{n\to +\infty}\|f-f_n\|=0$. The functional $\delta$ is not a bounded functional on the space $L^2_{per}$. Then  there exists a sequence $\{\psi_n\}\subset H^2_{per}$  with $\|\psi_n\|=1$ such that  $
\delta(\psi_n)=(\delta,\psi_n)=\psi_n(0)\to \infty$, as $n\to\infty$.
Since $\delta$ is a bounded linear on $H^2_{per}$, we can choose this sequence such that
$$
\lim_{n\to +\infty}\frac{(\delta,f_n)}{(\delta,\psi_n)}=0.
$$
Define the sequence $\zeta_n=f_n-(\delta,f_n)\psi_n/(\delta,\psi_n)$.
Then $\{\zeta_n\}\subset D(A^0)$ and 
$$
\|\zeta_n-f\|\leqq \|f_n-f\|+\Big |\frac{(\delta,f_n)}{(\delta,\psi_n)}\Big |\to 0 \quad{\rm{as}}\;\; n\to\infty.
$$
Thus the operator $A^0$ is densely defined.

\item  Since $(A-i)^{-1}\in B(H^{-2}_{\text{per}}; L^{2}_{\text{per}})$ we have $g_i\equiv (A-i)^{-1}\delta\in  L^{2}_{\text{per}}$.  Since $\widehat{\delta}(k)=1/{2\pi}$,   for $\psi\in D(A_0)\subset D(A)$ we obtain
\begin{equation}\label{sp3}
\begin{aligned}
\langle A^0\psi, g_i\rangle&=\langle A\psi, (A-i)^{-1}\delta\rangle=2\pi\sum_{k\in \mathbb Z}k^2\widehat{\psi}(k)\overline{\frac{1}{k^2-i}\widehat{\delta}(k)}\\
&=\psi(0)+2\pi\sum_{k\in \mathbb Z}\widehat{\psi}(k)\overline{i\widehat{g_i}(k)}=\psi(0)+\langle\psi, ig_i\rangle=\langle\psi, ig_i\rangle.
\end{aligned}
\end{equation}
 So, $g_i\in D({A^0}^*)$ and  ${A^0}^*g_i=ig_i$. A similar analysis show that $g_{-i}\in D({A^0}^*)$ and ${A^0}^*g_{-i}=-ig_{-i}$. 

\item The deficiency element $g_i$ is unique (up to multiplication by complex numbers).  We introduce the following norm $\|\cdot\|_{2,*}$ in the space $H^2_{per}([0,2\pi])$, which is equivalent to the standard norm in this space,
\begin{equation}\label{norm}
\|f\|_{2,*}^2 \equiv \|(-\partial_x^2-i) f\|^2=2\pi\sum_{k\in \mathbb Z}|(k^2-i)\widehat{f}(k)|^2=
\langle (\partial_x^4+1)^{1/2}f, (\partial_x^4+1)^{1/2}f\rangle.
\end{equation}
Since $\delta$ is a bounded linear on $(H^2_{per}, \|\cdot\|_{2,*}) $, the kernel $\mathcal K(\delta)=\{f\in H^2_{per}: \delta(f)=f(0)=0\}=D(A_0)$,  it is a hyperplane of codimension 1. Next for $h_0\equiv (A+i)^{-1}g_i\in H^2_{per}$ we have $h_0\perp \mathcal K(\delta)$. In fact, for $f\in \mathcal K(\delta)$
\begin{equation}\label{ker}
\langle (\partial_x^4+1)^{1/2}h_0, (\partial_x^4+1)^{1/2}f\rangle=\sum_{k\in \mathbb Z}\overline{\widehat{f}(k)}=\overline{f(0)}=0.
\end{equation}
Next, suppose $f_0\in D({A^0}^*)$ such that ${A^0}^*f_0=if_0$. Let $\psi\in D(A_0)\subset D(A)$, then $\langle A \psi, f_0\rangle=\langle A^0 \psi, f_0\rangle=\langle \psi, {A^0}^*f_0\rangle=\langle \psi, i f_0\rangle$. Therefore, $\langle (A+i) \psi, f_0\rangle=0$. Now, we show that $h_1\equiv (A+i)^{-1}f_0\in H^2_{per}$ satisfies that $h_1\perp \mathcal K(\delta)$. Let $\psi\in \mathcal K(\delta)$, then from the above analysis we obtain
$$
\langle (\partial_x^4+1)^{1/2}\psi, (\partial_x^4+1)^{1/2}h_1\rangle=2\pi\sum_{k\in \mathbb Z}
(k^2-i)\widehat{\psi}(k)\overline{\widehat{f_0}(k)}=\langle (A-i) \psi, f_0\rangle=0.
$$
So, there exists $\lambda\in \mathbb C$ such that $f_0=\lambda g_i$. This completes the proof of the Lemma.

\end{enumerate}
\end{proof}

\subsection{Deficiency elements $g_{\pm i}$}

Next we are interested in the profile of $g_{\pm i}$ which  will be crucial in our stability theory. We consider $\|g_{\pm i}\|=1$.  From \eqref{sp2} it follows that $g_{\pm i}$ represents the fundamental solution associated to $A\mp i$, respectively.  Next, we shall determine a formula for $g_{-i}\in L^{2}_{\text{per}}([0,2\pi])$. From \eqref{sp2} will be sufficient to find $\mathcal K_i \in L^{2}_{\text{per}}([0,2\pi])$ such that
\begin{equation}\label{sp4}
\widehat{\mathcal K_i}(k)=\frac{1}{k^2+i},
\end{equation}
since
$$
\widehat{g_{-i}}(k)=\frac{1}{2\pi}\frac{1}{k^2+i}
$$
implies $g_{-i}=\frac{1}{2\pi}\mathcal K_i(x)$  (we can also to obtain this formula via the following equality in the distributional sense, $g_{-i}=(A+i)^{-1}\delta=\delta\star \mathcal K_i=\frac{1}{2\pi}\mathcal K_i$). The formula for $g_{i}$ is obtained from relation $g_{i}=\overline{g_{-i}}$. Next, we find explicitly $g_{-i}$.  So,  for $\psi\in C^\infty_{\text{per}}$ we solve
\begin{equation}\label{sp6}
\Big(-\frac{d^2}{dx^2}+i\Big)h=\psi.
\end{equation}
We start by finding a specific base for the homogeneous equation
\begin{equation}\label{sp7}
y''-iy=0. 
\end{equation}
For $\beta=\frac{1+i}{\sqrt{2}}$ and $\mu=\frac{-1-i}{\sqrt{2}}$ the general solution for the second-order equation in \eqref{sp7} is given by $y(\xi)=ae^{\beta \xi} + be^{\mu \xi}$. We consider the following base $\mathcal B=\{y_1, y_2\}$ for the set of solutions of \eqref{sp7},
\begin{equation}\label{sp9}
y_1(\xi)=\cosh(\beta \xi),\;\;\;\;\;\;\;\;
y_2(\xi)=\cosh\Big(\beta (\xi+\pi)\Big).
\end{equation}
So, we have that  the Wronskian is given by $W(y_1,y_2)=\beta\sinh(\beta\pi)$.  Next for $f=-\psi$ we find a particular {\it periodic solution} $y_p$ of the equation
\begin{equation}\label{sp10}
y''-iy=f,
\end{equation}
which via the variational parameters method is given by
\begin{equation}\label{sp11a}
y_p=u_1y_1+u_2y_2,
\end{equation}
where
\begin{equation}\label{sp11b}
\begin{cases}
\begin{aligned}
&u'_1=-\frac{y_2f}{W}=-\frac{ \cosh\Big(\beta(\xi+\pi)\Big) f}{W},\\
&u'_2=\frac{y_1f}{W}=\frac{\cosh(\beta \xi) f}{W}.
\end{aligned}
\end{cases}
\end{equation}
Then
\begin{equation}\label{sp12a}
\begin{cases}
\begin{aligned}
&u_1(\xi)=-\frac{1}{W}\int_0^\xi \cosh\Big(\beta(x+\pi)\Big)f(x)dx+\alpha_0,\\
&u_2(\xi)=-\frac{1}{W}\int_\xi^{2\pi} \cosh(\beta x)f(x)dx+\beta_0,
\end{aligned}
\end{cases}
\end{equation}
for $\alpha_0, \beta_0$ integration constants to be chosen later. So, after some calculations we obtain for $\xi\in \mathbb R$ the formula
\begin{equation}\label{sp11c}
\begin{aligned}
&y_p(\xi)=-\frac{1}{2W}\int_0^{2\pi} \cosh\Big(2\pi\beta\Big(\frac{\xi-x}{2\pi}-\Big [\frac{\xi-x}{2\pi}\Big]-\frac12\Big)\Big)f(x)dx\\
&\quad+\alpha_0\cosh(\beta \xi)+ \beta_0\cosh(\beta(\xi+\pi))-\frac{1}{2W}\int_0^{2\pi} \cosh\Big(\beta\Big(\xi+x+\pi\Big)\Big)f(x)dx.
\end{aligned}
\end{equation}
Here $[\cdot]$ stands for the integer part. Next we can choose $\alpha_0, \beta_0$ such that the second line in \eqref{sp11c} will be zero, more exactly we have the following choices 
\begin{equation*}
\begin{aligned}
\alpha_0&=-\frac{1}{2\beta \sinh^2(\beta \pi)}\int_0^{2\pi}\sinh(\beta x)\psi(x)dx,\\
\beta_0&=\frac{1}{2\beta \sinh^2(\beta\pi)}\int_0^{2\pi}\sinh\Big(\beta(x+\pi)\Big)\psi(x)dx.
\end{aligned}
\end{equation*}
Therefore we have that $y_p$ is a periodic function with  a minimal period $2\pi$ and has the convolution expression
\begin{equation}\label{sp12}
y_p(\xi)=\mathcal K_{i}\star\psi(\xi)\qquad \xi\in\mathbb R,
\end{equation}
where $\mathcal K_{i}\in L^2_{per}([0,2\pi])$ is defined by
\begin{equation}\label{sp13}
\mathcal K_{i}(x)=\frac{2\pi}{2\beta \sinh(\beta \pi)}\cosh\Big(2\pi\beta\Big(\frac{x}{2\pi}-\Big [\frac{x}{2\pi}\Big ]-\frac12\Big)\Big),\qquad x\in \mathbb R.
\end{equation}
Therefore $\mathcal K_{i}$ satisfies \eqref{sp4}. So, we get the profile
\begin{equation}\label{sp14}
g_{-i}(x)=\frac{1}{2\beta \sinh(\beta \pi)}\cosh\Big(\beta\Big(|x|-\pi\Big)\Big), \quad{\rm{for}}\;\; x\in [-\pi,\pi].
\end{equation}

Lastly, we obtain the expression for  the deficiency element $g_{-i}$. For $\sigma=1/(2\beta \sinh(\beta\pi))$  and $x\in [-\pi,\pi ]$  
\begin{equation}\label{sp17}
g_{-i}(x)=\sigma\Big[\cosh\Big(\frac{|x|-\pi}{\sqrt{2}}\Big)\cos\Big(\frac{|x|-\pi}{\sqrt{2}} \Big)+i\sinh\Big(\frac{|x|-\pi}{\sqrt{2}} \Big)\sin\Big(\frac{|x|-\pi}{\sqrt{2}} \Big)\Big]. 
\end{equation}
See Figure 1 and Figure 2 below for the profile of the real and imaginary parts of $g_{-i}$, $\Re(g_{-i})$ and $\Im(g_{-i})$, respectively. 

\begin{figure}[h]
\centering
\begin{minipage}[b]{0.5\linewidth}
\includegraphics[angle=0,scale=0.3]{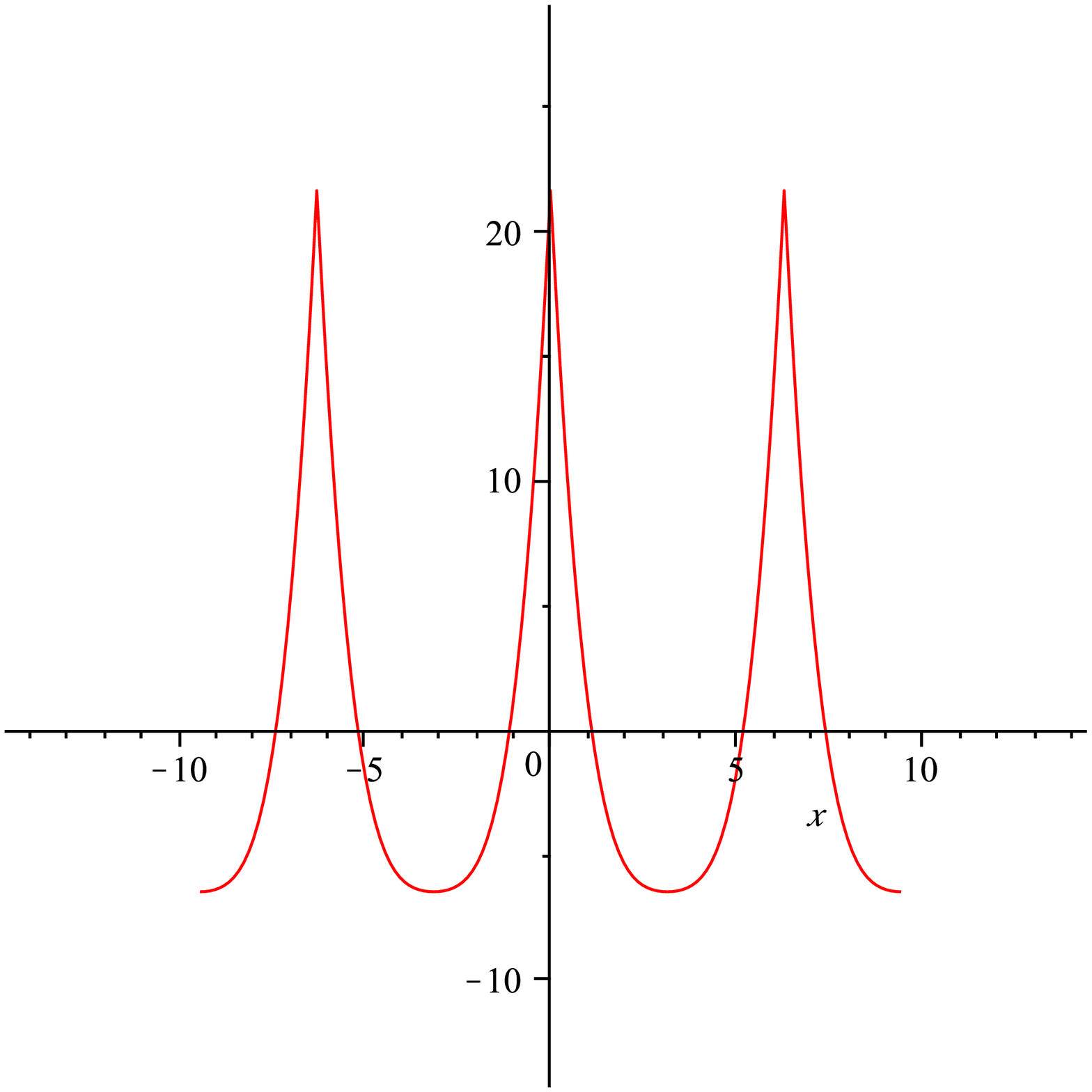}
\caption{Graphic of the function $\Re(g_{-i})$ given by
\eqref{sp17}}
\end{minipage}\hfill
\begin{minipage}[b]{0.5\linewidth}
\includegraphics[angle=0,scale=0.25]{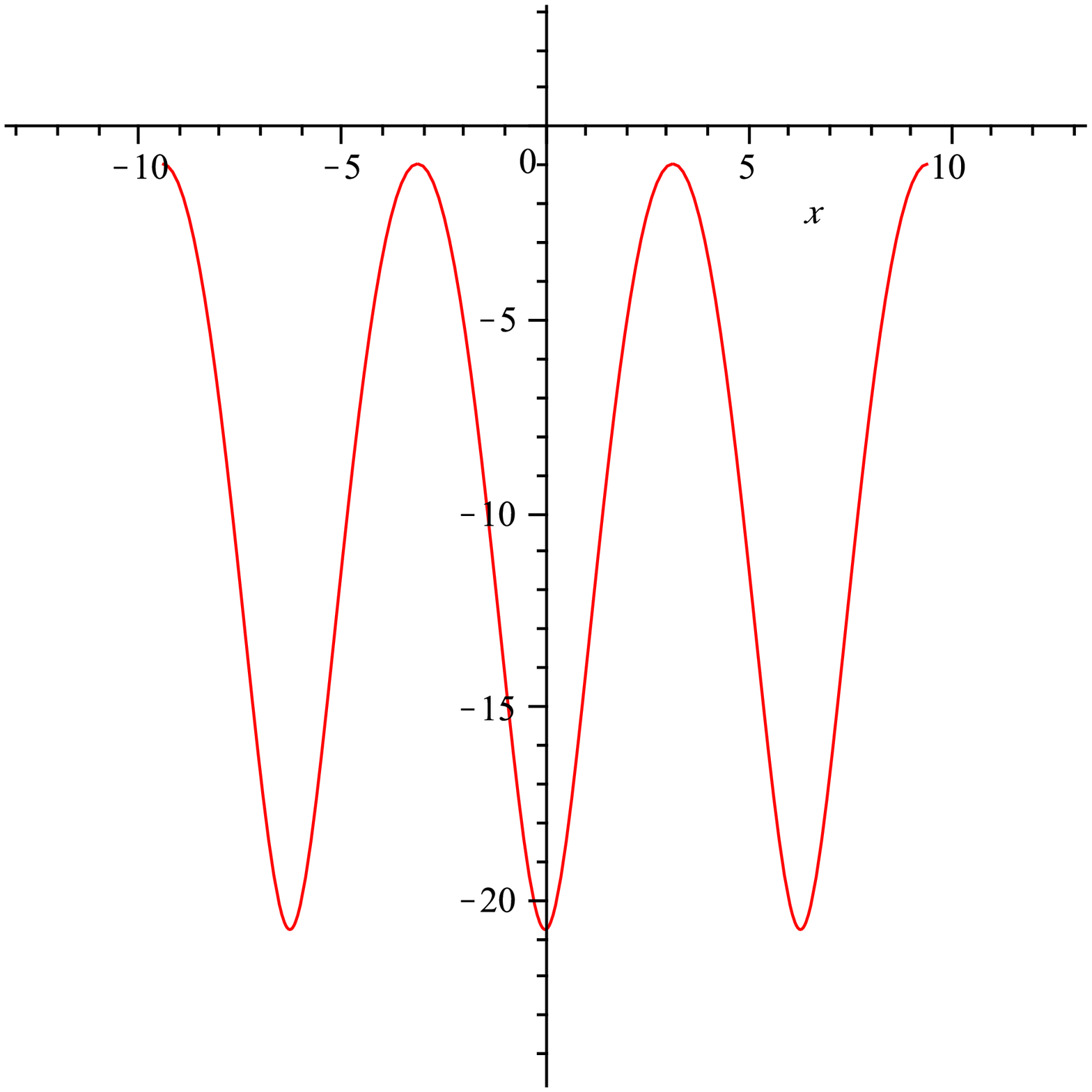}
\caption{Graphic of the function $\Im(g_{-i})$ given by
\eqref{sp17}}
\end{minipage}
\end{figure}

In the next subsection we need to use that the deficiency elements $g_{\pm i}$ have $L^2_{\text{per}}([0,2\pi])$-norm equal to $1$. So, for $\|g_i\|^2=\|g_{-i}\|^2=\theta$ with
$$
\theta=\frac{\sqrt{2}}{4}\frac {\sinh(\sqrt{2}\pi)+\sin(\sqrt{2}\pi)}{\cosh(\sqrt{2}\pi)-\cos(\sqrt{2}\pi)} 
$$
we obtain the normalized deficiency elements $\widetilde{g}_{\pm i}=\frac{g_{\pm i}}{\|g_{\pm i\|}}$. But for convenience of notation we will continue to use $g_{\pm i}$.

{\bf{Remark:}} We note that   $\Re(g_{-i})$ has the peaks in $\pm 2n\pi$, $n\in \mathbb Z$, and $\Im(g_{-i})$ is a smooth periodic function.

\subsection{Self-adjoint extensions of $A^0$}

In this subsection we present explicitly  all the self-adjoint extensions of the symmetric operator $A^0$ defined in Lemma \ref{sp0}, which will be parametrized by the strength $\gamma$. 

From Lemma \ref{sp0}  we have that the {{\it deficiency indices} and   {\it deficiency subspaces} of $A^0$ are given by
\begin{equation}\label{sa2}
n_{+-}=1\quad{\rm{and}}\quad \mathcal D_+=[g_i],\;\; \mathcal D_-=[g_{-i}],
\end{equation}
where $g_{-i}$ is given by \eqref{sp14} and $g_i=\overline{g_{-i}}$. Next, let $B$ be a closed symmetric extension of $A^0$. Then for $\varphi\in D(B^*)$, we have 
$$
\langle \psi,B^*\varphi\rangle=\langle B\psi,\varphi\rangle=\langle A\psi,\varphi\rangle\quad{\rm{for\; all}}\;\;\psi\in D(A^0).
$$
Thus $\varphi\in D({A^0}^*)$ and $B^*\varphi={A^0}^*\varphi$, therefore we obtain the basic relation
\begin{equation}\label{sa3}
A^0\subseteq B\subseteq B^*\subseteq{A^0}^*.
\end{equation}
So, from \eqref{sa2}-\eqref{sa3}  and from the von Neumann extension theory for symmetric operators \cite{RS} we have that all the closed symmetric extensions of $A^0$ are self-adjoint and coincides with the restriction of the operator ${A^0}^*$. Moreover, there is a one-one correspondence between self-adjoint extensions of $A^0$ and unitary maps from $\mathcal D_+$ onto
$\mathcal D_-$. Hence,  if $U$ is  such an isometry with initial space $\mathcal D_+$ then there exists  $\theta\in[0,2\pi)$ such that
$$
U(\lambda g_i)=\lambda e^{i\theta}g_{-i},\quad {\rm{for\; all}}\;\lambda\in \mathbb C.
$$
Then via this identification for $\theta\in[0,2\pi)$  the corresponding self-adjoint extension $A^0(\theta)$ of $A^0$ is defined as follows;
\begin{equation}\label{sa4}
\left\{
\begin{aligned}
&D(A^0(\theta))=\{\psi+\lambda g_i+\lambda e^{i\theta}g_{-i} :\psi\in D(A^0), \lambda\in \mathbb C\},\\
&A^0(\theta)(\psi+\lambda g_i+\lambda e^{i\theta}g_{-i})={A^0}^*(\psi+\lambda g_i+\lambda e^{i\theta}g_{-i})=A^0\psi+i\lambda g_i-i\lambda e^{i\theta}g_{-i}.
\end{aligned}\right.
\end{equation}

For our purposes we will parametrize the self-adjoint extensions  $A^0(\theta)$ with the strength parameter $\gamma\in \mathbb R\cup\{+\infty\}$ instead of the parameter $\theta$ appeared in the von Neumann formulas \eqref{sa4}. So, we obtain from \eqref{sp14} that for $\zeta\in D(A^0(\theta))$,  in the form $\zeta=\psi+\lambda g_i+\lambda e^{i\theta}g_{-i}$, we have the basic expression
\begin{equation}\label{sa5}
\zeta'(0+)-\zeta'(0-)=-\lambda(1+e^{i\theta}).
\end{equation}
Next we find $\gamma$ such that $\gamma\zeta(0)=-\lambda(1+e^{i\theta})$. Indeed, after some calculations we find the formula
\begin{equation}\label{sa6}
\gamma(\theta)=\frac{-2\cos(\theta/2)}{\Re[\coth(\beta\pi)e^{i(\frac{\theta}{2}-\frac{\pi}{4})}]},
\end{equation}
which can be  write as
\begin{equation}\label{sa7}
\gamma(\theta)=\frac{-4|\sinh(\beta\pi)|^2\cos(\theta/2)}{\sinh(\sqrt{2}\pi)\cos(\frac{\theta}{2}-\frac{\pi}{4})+\sin(\sqrt{2}\pi)\sin(\frac{\theta}{2}-\frac{\pi}{4})}.
\end{equation}
Therefore, if $\theta$ varies  in $[0,2\pi)$, $\gamma=\gamma(\theta)$ varies  in $\mathbb R\cup \{+\infty\}$. For the unique $\theta_0 \in [0,2\pi)$ such that $\Re[\coth(\beta \pi)e^{i(\frac{\theta_0}{2}-\frac{\pi}{4})}]=0$ we have $
\lim_{\theta\uparrow \theta_0}\gamma(\theta)=+\infty$.


So, from now on we parametrize all self-adjoint extensions of $A^0$ with the help of $\gamma$. Thus we get,

\begin{teo}\label{self} All self-adjoint extensions of $A^0$ are given for $-\infty<\gamma\leqq +\infty$ by 
\begin{equation}\label{sa8}
\begin{aligned}
-\Delta_{\gamma}&=-\frac{d^2}{dx^2}\\
D(-\Delta_{\gamma})&=\{\zeta\in H^1_{\text{per}}([-\pi, \pi])\cap H^2((-\pi,\pi)-\{0\})\cap H^2((2n \pi, 2(n+1)\pi): \\
&\qquad\qquad\qquad \zeta'(0+)-\zeta'(0-)=\gamma \zeta(0)\}.
\end{aligned}
\end{equation}
The special case $\gamma=0$ just leads to the kinetic energy hamiltonian $-\Delta$ in $L^2_{\text{per}}([-\pi,\pi])$, 
\begin{equation}\label{sa9}
-\Delta=-\frac{d^2}{dx^2},\qquad D(-\Delta)= H^2_{\text{per}}([-\pi, \pi]),
\end{equation}
whereas the case $\gamma=+\infty$ yields a {\it Dirichlet-periodic} boundary condition at zero, 
\begin{equation}\label{sa9a}
D(-\Delta_{+\infty})=\{\zeta\in H^1_{\text{per}}([-\pi, \pi])\cap H^2((-\pi,\pi)-\{0\})\cap H^2((2\pi n, 2(n+1)\pi):  \zeta(0)=0\}.
\end{equation}

\end{teo}

\begin{proof}  By the arguments sketched above we obtain
\begin{equation}\label{sa10}
A^0(\theta)\subset -\Delta_{\gamma}
\end{equation}
with $\gamma=\gamma(\theta)$ given in \eqref{sa7}. But $-\Delta_{\gamma}$ is easily seen to be symmetric in the corresponding domain $D(-\Delta_{\gamma})$ for all $-\infty<\gamma\leqq +\infty$, which implies the relation
$$
A^0(\theta)\subset -\Delta_{\gamma}\subset (-\Delta_{\gamma})^* \subset A^0(\theta).
$$
It completes the proof of the Theorem.
\end{proof}

{\bf{Remarks:}} 

\begin{enumerate}

\item  Since $-\frac{d^2}{dx^2}g_{\pm i}(x)=\pm i g_{\pm i}(x)$, for $x\neq 2n\pi$ ($n\in \mathbb Z$), we have
for $\zeta\in D(A^0(\theta))$  the relation
\begin{equation*}
\begin{aligned}
A^0(\theta)(\zeta)(x)&=A^0(\theta)(\psi+\lambda g_i+\lambda e^{i\theta}g_{-i})(x)\\
&=-\frac{d^2}{dx^2}\psi(x)-\lambda \frac{d^2}{dx^2}g_i(x)-\lambda e^{i\theta}\frac{d^2}{dx^2}g_{-i}(x)=-\Delta_\gamma \zeta(x),
\end{aligned}
\end{equation*}
which implies relation \eqref{sa10}.
\item For $\zeta\in D(A^0(\theta))$ it follows $\zeta\in H^2((2n+1)\pi, (2n+3)\pi)-\{2(n+1)\pi\})$ for 
$n\in \mathbb Z$. Obviously we have $\zeta'(2(n+1)\pi+)-\zeta'(2(n+1)\pi-)=\gamma \zeta(2(n+1)\pi)$.

\item For a characterization of the domains $D(-\Delta_{\gamma})$ for all $-\infty<\gamma\leqq +\infty$ we refer to the reader to Theorem \ref{resol4} below and remarks associated to  it.
In particular,  from \eqref{expan} we obtain for the extreme case  $\gamma=+\infty$ that  all element $\zeta\in D(-\Delta_{+\infty})$ has the decomposition
$$
\zeta(x)=\psi_{-i}(x)-\frac{2\beta}{\coth(\beta\pi)}\psi_{-i}(0)g_{-i}(x), 
$$
where $\psi_{-i}\in H^2_{per}([-\pi,\pi])$.

\item The expression  {\it Dirichlet-periodic boundary condition at zero} emerges  in a natural form since every element $\zeta\in D(-\Delta_{+\infty})$ belongs to the space $H^1_{\text{per}}([-\pi, \pi])\cap H^2((-\pi,\pi)-\{0\})\cap H_0^2((0,2\pi))$.

\item By definition $-\Delta_\gamma$ describes a {\it periodic $\delta$-interaction of strength $\gamma$ centered at zero}. In other words, equation \eqref{sa8} is the precise formulation of the {\it formal linear differential operator} \eqref{deltaop}, namely, for $\zeta\in D(-\Delta_\gamma)$, $-\frac{d^2}{dx^2} \zeta=(-\frac{d^2}{dx^2}+\gamma\delta ) \zeta$ in a distributional sense. 
 
 \item An informal calculation shows that  the jump condition in \eqref{sa8}  is \lq\lq quite natural''. Indeed, consider the Schr\"odinger equation $-\zeta''+\gamma \delta\zeta=\lambda \zeta$ and ``integrate'' from $-\epsilon$ to $\epsilon$, then
$$
-\zeta'(\epsilon)+\zeta'(-\epsilon)+\gamma\zeta(0)=\lambda\int_{-\epsilon}^\epsilon\zeta (x)dx.
$$
If $\epsilon\to 0^+$  we obtain $\zeta'(0+)-\zeta'(0-)=\gamma\zeta(0)$.
\end{enumerate}

\subsection{Resolvents and spectrum for $-\frac{d^2}{dx^2} +\gamma \delta$}

In this subsection we  study the solvability of  the model in \eqref{deltaop} in a periodic context. So, we will show that  their resolvents can be given explicitly in terms of the interactions strengths $\gamma$. It will be shown that the spectrum and the eigenfunctions can be given explicitly. Here we use the Krein's formula  for  the resolvents of two self-adjoint extensions of one symmetric operator (see \cite{ak}). The results to be established here can be use for showing that the family of self-adjoint operator  $\mathcal L_{1,Z}$ and  $\mathcal L_{2,Z}$ in \eqref{specA6}-\eqref{specA7}, are real-analytic in the sense of Kato (see section 6).

We star with the following basic result.

\begin{teo}\label{resol1}  The resolvent of $-\Delta$ in $L^2_{per}([-\pi, \pi])$ is given for $k^2\in \rho(-\Delta)$ by
\begin{equation}
(-\Delta-k^2)^{-1}f=J_k\star f,\quad k\neq n,\;n\in \mathbb Z,
\end{equation}
where the integral kernel $J_k\in  L^2_{per}([-\pi, \pi])$ is defined by
\begin{equation}\label{inker}
J_k(\xi)=\frac{2\pi}{2ik\sinh(ik\pi)}\cosh\Big(ik\Big(|\xi|-\pi\Big)\Big),\quad{\rm{for}}\;\;\xi\in [-\pi,\pi].
\end{equation}
with $k$ being a square root of $k^2$.
\end{teo}

\begin{proof} The proof follows the same ideas explained  in subsection 3.1,  so it will be omitted.
\end{proof}

{\bf Remark:}} From \eqref{inker}) it follows  that the set of singularities of $J_k$, $\{n\,:\,n\in\mathbb Z\}$, produces the well-known set of eigenvalues associated to the operator $-\Delta$, namely, $\{n^2\,:\,n\in\mathbb N\}$.

\vskip0.1in

Next,  we shall describe  the resolvent of the self-adjoint operators $-\Delta_\gamma$.

\begin{teo}\label{resol2}  The resolvent of $-\Delta_\gamma$ in $L^2_{per}([-\pi, \pi])$ is given for $k\neq n$, $n\in \mathbb Z$, by 
\begin{equation}\label{resol2a}
\begin{aligned}
& (-\Delta_\gamma-k^2)^{-1}= (-\Delta-k^2)^{-1}-\frac{1}{4\pi^2}\frac{2i\gamma k}{\gamma \coth(ik\pi)+2ik}\langle\cdot, \overline{J_k}\rangle J_k,\\
&\qquad k^2\in \rho(-\Delta_\gamma),\quad k\;\;\text{being a square root of}\;\;k^2,\quad -\infty<\gamma\leqq +\infty.
\end{aligned}
\end{equation}
Therefore, $-\Delta_\gamma$ has a compact resolvent for $-\infty<\gamma\leqq +\infty$.
 \end{teo}
 
\begin{proof}  Let $\gamma \neq +\infty$ and  $k$  such that $\gamma \coth(ik\pi)\neq -2ik$. For $k\neq n$, $n\in \mathbb Z$,  and $f\in L^2_{per}([-\pi, \pi])$ define
\begin{equation}\label{resol3a}
h_\gamma(x)=[(-\Delta-k^2)^{-1}f](x)-\frac{1}{4\pi^2}\frac{2i\gamma k}{\gamma \coth(ik\pi)+2ik}\langle f, \overline{J_k}\rangle J_k(x).
\end{equation}
It is easy  to see that $h_\gamma\in H^1_{per}([-\pi,\pi])\cap H^2((-\pi,\pi)-\{0\}))$. Since $J'_k(0+)-J'_k(0-)=-2\pi$ we obtain
\begin{equation}\label{resol3b}
h'_\gamma(0+)-h'_\gamma(0-)=\frac{1}{2\pi}\frac{2i\gamma k}{\gamma \coth(ik/2)+2ik}\int J_k(y)f(y)dy=\gamma h_\gamma(0).
\end{equation}
Therefore equation \eqref{resol3b} implies that $h_\gamma\in D(-\Delta_\gamma)$. Next, since for $x\in \mathbb R-2\pi\mathbb Z$
\begin{equation}\label{resol3c}
-J''_k(x)-k^2J_k(x)=0,
\end{equation}
from Theorem \ref{self}  follows that
\begin{equation}\label{resol3d}
[(-\Delta_\gamma-k^2)h_\gamma](x)=-h''_\gamma(x)-k^2h_\gamma(x)=f(x),\quad{\rm{for}}\;\;x\in \mathbb R-2\pi\mathbb Z.
\end{equation}
Hence we obtain \eqref{resol2a}.

Let $\gamma =+\infty$ and $k\neq n$, $n\in \mathbb Z$. Since $\coth(ir\pi)=0$ if and only if $r\in \mathbb R$ and $r\in \mathbb Z$,  the following formula for the resolvent of $-\Delta_{+\infty}$ is well defined  
\begin{equation}\label{resolinf}
(-\Delta_{+\infty}-k^2)^{-1}= (-\Delta-k^2)^{-1}-\frac{1}{4\pi^2}\frac{2ik}{\coth(ik\pi)}\langle\cdot, \overline{J_k}\rangle J_k.
\end{equation}

Finally, combining  \eqref{resol2a} and \eqref{resolinf} we obtain that $-\Delta_\gamma$ has a compact resolvent for $-\infty<\gamma\leqq +\infty$,  so the spectrum of $-\Delta_\gamma$, $\sigma(-\Delta_\gamma)$, is a infinity enumerable set of eigenvalues $\{\mu_n\}_{n\geqq 0}$ such that
$$
\mu_0< \mu_1\leqq\mu_2\leqq \cdot\cdot\cdot
$$
and $\mu_n\to +\infty$ as $n\to\infty$.
\end{proof}

{\bf{Remarks:}} \begin{enumerate} \item $J_k\notin D(-\Delta_\gamma)$ for $k$ such that $\gamma \coth(ik\pi)\neq -2ik$. Indeed, 
\begin{equation}\label{rema3}
J'_k(0+)-J'_k(0-)=-2\pi\neq \gamma J_k(0).
\end{equation}

\item $J_k\in H^1_{\text{per}}([-\pi, \pi])\cap H^2((-\pi,\pi)-\{0\})\cap H^2((2\pi n,2(n+1)\pi))$, and satisfies \eqref{resol3c} in $(-\pi,\pi)-\{0\}$ with $J'_k(\pm \pi)=0$.
\end{enumerate}

\vskip0.1in

Next we have additional domain properties of $-\Delta_\gamma$ and point out the locality of the periodic $\delta$-interactions.

\begin{teo}\label{resol4r} The domain $D(-\Delta_\gamma)$, $-\infty< \gamma\leqq +\infty$, consiste of all elements $\zeta$ of the type
\begin{equation}\label{resol3}
\zeta(x)=\psi_k(x)-\frac{1}{2\pi}\frac{2i\gamma k}{\gamma \coth(ik\pi)+2ik}\psi_k(0)J_k(x),\qquad x\in \mathbb R-2\pi\mathbb Z,
\end{equation}
where $\psi_k\in D(-\Delta)=H^2_{per}([-\pi,\pi])$, $k^2\in \rho(-\Delta_\gamma)$, $k$ being a square root of $k^2$, and $k\neq n$. The decomposition \eqref{resol3} is unique and with $\zeta\in D(-\Delta_\gamma)$ of this form it follows that
\begin{equation}\label{resol4}
(-\Delta_\gamma-k^2)\zeta=(-\Delta-k^2)\psi_k.
\end{equation}
Also if  $\zeta\in D(-\Delta_\gamma)$ such  that $\zeta=0$ in an open set $\mathcal O\subset \mathbb R$, then $-\Delta_\gamma\zeta=0$ in $\mathcal O$.
\end{teo}

\begin{proof} Since $(-\Delta-k^2)^{-1}(L^2_{per})=D(-\Delta)$, one has that
$$
(-\Delta_\gamma-k^2)^{-1}(-\Delta-k^2)D(-\Delta)=(-\Delta_\gamma-k^2)^{-1}(L^2_{per})=D(-\Delta_\gamma).
$$
Therefore for every $\zeta\in D(-\Delta_\gamma)$ there exists $\psi_k\in D(-\Delta)$ such that from \eqref{resol2a} we obtain
\begin{equation}\label{resol4a}
\begin{aligned}
\zeta&=(-\Delta_\gamma-k^2)^{-1}(-\Delta-k^2)\psi_k\\
&\quad=\psi_k - \frac{1}{4\pi^2}\frac{2i\gamma k}{\gamma \coth(ik\pi)+2ik}\langle (-\Delta-k^2)\psi_k, \overline{J_k}\rangle J_k.
\end{aligned}
\end{equation}
Next we prove $\langle (-\Delta-k^2)\psi_k, \overline{J_k}\rangle=2\pi\psi_k(0)$. Indeed, combining \eqref{resol3c}, the Remark-(1) after the proof of Theorem \ref{resol2}, and $\psi_k\in H^2_{per}$ it follows that
\begin{equation}\label{resol4b}
\begin{aligned}
\langle (&-\Delta-k^2)\psi_k, \overline{J_k}\rangle=\lim_{\epsilon \downarrow 0}\int_{-\pi}^{-\epsilon}(-\psi''_k(x)-k^2\psi_k(x))J_k(x)dx +\\
&\lim_{\epsilon \downarrow 0}\int_{\epsilon}^{\pi}(-\psi''_k(x)-k^2\psi_k(x))J_k(x)dx=\psi_k(0)[J'_k(0-)-J'_k(0+)]=2\pi\psi_k(0).
\end{aligned}
\end{equation}

Next, we  prove uniqueness of the decomposition in \eqref{resol3}. Let $\zeta=0$, so
$$
\psi_k(x)=\frac{1}{2\pi}\frac{2i\gamma k}{\gamma \coth(ik\pi)+2ik}\psi_k(0)J_k(x), \quad x\in (-\pi,\pi)-\{0\}.
$$
Since $\psi_k\in H^2_{per}$ it follows immediately that $\psi_k\equiv 0$. Now, relation \eqref{resol4} simply follows from the equality
$$
(-\Delta_\gamma-k^2)^{-1}(-\Delta-k^2)\psi_k=\psi_k-\frac{1}{4\pi^2}\frac{2i\gamma k}{\gamma \coth(ik\pi)+2ik}\langle (-\Delta-k^2)\psi_k, \overline{J_k}\rangle J_k=\zeta.
$$

To prove locality we assume first $2\pi\mathbb Z \cap \mathcal O=\emptyset $. Then is immediate that the relation 
$$
\Big(-\frac{d^2}{dx^2}-k^2\Big)J_k (x)=0,\qquad{\rm{for}}\;\;x\in \mathcal O
$$
and \eqref{resol4} imply that for $x\in \mathcal O$
\begin{equation}\label{local}
\begin{aligned}
(-\Delta_\gamma \zeta)(x)&=k^2\zeta(x)+ \Big(-\frac{d^2}{dx^2}-k^2\Big)\psi_k (x)\\
& =\frac{1}{2\pi} \frac{2i\gamma k}{\gamma \coth(ik\pi)+2ik}\psi_k(0)\Big(-\frac{d^2}{dx^2}-k^2\Big)J_k(x)=0.
\end{aligned}
\end{equation} 
On the other hand, if there exists $n\in \mathbb Z$ such that $2\pi n\in \mathcal O$ then $\zeta(0)=\zeta(2\pi n)=0$. Therefore from the definition of $D(-\Delta_\gamma)$ we have $\zeta'(0+)=\zeta'(0-)$ and so $\zeta'(0)$ exists. But relation \eqref{resol3} then implies that $J_k'(0)$ exists if $\psi_k(0)\neq 0$. Hence we need to have that $\psi_k(0)= 0$ and then $\zeta=\psi_k\in H^2_{per}([-\pi,\pi])$. So it follows $\zeta\in D(A^0)$ and it implies that
$$
-\Delta_\gamma \zeta(x)=-\frac{d^2}{dx^2}\zeta(x)=0,\quad {\rm{for}}\;\;x\in \mathcal O.
$$
This completes the proof of the Theorem.
\end{proof}

\begin{coro}\label{resol5a} The domain $D(-\Delta_\gamma)$, $-\infty< \gamma\leqq +\infty$, consists of all elements $\zeta$ of the type
\begin{equation}\label{expan}
\zeta(x)=\psi(x)-\frac{2\gamma \beta}{\gamma \coth(\beta \pi)+2\beta}\;\psi(0)g_{-i}(x),\qquad x\in \mathbb R-2\pi\mathbb Z
\end{equation}
for $\psi\in H^2_{per}([-\pi,\pi])$.
\end{coro}
\begin{proof}  For $k$ being a root of $k^2=-i$ in Theorem \ref{resol4r} we have that  $\beta=ik$ is a root of $i$, and so the corresponding function $J_k$ in \eqref{resol3} is given by $g_{-i}=\frac{1}{2\pi}\mathcal K_{i}$.
\end{proof}

{\bf{Remark:}}  From  Theorem \ref{resol4r} we obtain that if $\zeta\in D(-\Delta_\gamma)$ and $\zeta(0)=0$ then $\zeta\in H^2_{per}([-\pi, \pi])$.

\vskip0.1in

Next, we deduce some spectral properties of $-\Delta_\gamma$.  This information will be relevant for  our well-posedness results.

\begin{teo}\label{resol5b} Let $-\infty<\gamma\leqq +\infty$. Then the spectrum of $-\Delta_\gamma$ is discrete $\{\theta_{j,\gamma}\}_{j\geqq 1}$ and such that $\theta_{1,\gamma} < \theta_{2,\gamma}\leqq \theta_{3,\gamma}\leqq \cdot\cdot\cdot$.

If $-\infty<\gamma<0$, $-\Delta_\gamma$ has precisely one negative, simple eigenvalue, i.e.,
\begin{equation}\label{resol7}
\sigma_{\text{p}}(-\Delta_\gamma)\cap(-\infty,0)=\{-\mu^2_\gamma\}
\end{equation}
where $\mu_\gamma$ is positive and satisfies $\gamma=-2\mu_\gamma tanh(\mu_\gamma\pi)$. The function
\begin{equation}\label{resol8}
\psi_\gamma(\xi)=\frac{J_{i\mu_\gamma}(\xi)}{\|J_{i\mu_\gamma}\|}=\frac{2\pi}{2\|J_{i\mu_\gamma}\|\mu_\gamma\sinh(\mu_\gamma\pi)}\cosh\Big(\mu_\gamma\Big(|\xi|-\pi\Big)\Big),\quad{\rm{for}}\; \xi\in[-\pi,\pi]
\end{equation}
is the strictly positive (normalized) eigenfunction associated to the eigenvalue $-\mu^2_\gamma$.  The nonnegative eigenvalues  (are nondegenerated) are ordered in the increasing form
$$
0<\kappa^2_1<1<\kappa_2^2<2^2<\cdot\cdot\cdot<\kappa_j^2<j^2<\cdot\cdot\cdot
$$
where for $j\geqq 1$, $\kappa_{j}$ is the only solution of the equation
\begin{equation}\label{resol8a}
cot(\kappa \pi)=\frac{2\kappa}{\gamma}
\end{equation}
in the interval $(j-\frac12, j)$. The eigenfunction associated with $\kappa_{j}$ is $J_{\kappa_{j}}\in D(-\Delta_\gamma)$. The sequence $\{j^2\}_{j\geqq 1}$ is the classical set of eigenvalues associated to the operator $-\Delta$ with associated eigenfunctions $\{\sin(j x): j\geqq 1\}\subset D(-\Delta_\gamma)$.

If $\gamma> 0$, $-\Delta_\gamma$ has no negative eigenvalues and the positive  eigenvalues (are nondegenerated) are ordered in the increasing form
\begin{equation}\label{resol9}
0<k^2_1<1<k_2^2<2^2<\cdot\cdot\cdot<k_j^2<j^2<\cdot\cdot\cdot
\end{equation}
where for $j\geqq 0$, the eigenvalue $k_{j+1}$ is the only solution of the equation
\begin{equation}\label{resol8b}
cot(k\pi)=\frac{2k}{\gamma}
\end{equation}
in the interval $(j, j+\frac12)$. The eigenfunction associated with $k^2_{j+1}$ is $J_{k_{j+1}}\in D(-\Delta_\gamma)$. The sequence $\{j^2\}_{j\geqq 1}$ is the classical set of eigenvalues associated to the operator $-\Delta$ with associated eigenfunctions $\{\sin(j x): j\geqq 1\}\subset D(-\Delta_\gamma)$.  

Zero is not eigenvalue of $-\Delta_\gamma$ for all $\gamma\neq 0$.

For $\gamma=+\infty$, $\sigma(-\Delta_{+\infty})=\{j^2\}_{j\geqq 1}$ and with associated eigenfunctions $\{\sin(j x): j\geqq 1\}\subset -\Delta_{+\infty}$. The eigenvalues are nondegenerated.

\end{teo}

\begin{proof} We divide the proof into several steps.
 \begin{enumerate}

\item For $\gamma\neq 0$ the eigenvalues $j^2$, $j\in \mathbb N$, $j\geqq 1$, are simple. In fact, it is immediate that $\psi_j(x)=\sin(j x)\in D(-\Delta_\gamma)$ and $-\Delta_\gamma\psi_j=-\psi''_j(x)=j^2\psi_j(x)$ for $x\in \mathbb R$. We known that the next equation 
$$
-\psi''=j^2\psi\qquad {\rm{on}}\;\; (0,2\pi) \;\;({\rm{similarly \;\;on}}\;\;(-2\pi,0))
$$
for $\psi\neq 0$, has exactly two linearly independent solutions. So $j^2$ is simple on $D(-\Delta_\gamma)$. Moreover,  for every $\psi\in D(-\Delta_\gamma)$ satisfying $-\Delta_\gamma\psi=j^2\psi$, we have $\psi(x)=\alpha \sin (j x)+\beta \cos(jx)$ for $x\in (0,2\pi)$ and $x\in (-2\pi,0)$. Then $\gamma \psi(0)=\psi'(0+)-\psi'(0-)=0$, implies $\beta=0$.

\item For $\gamma<0$, $-\mu^2_\gamma$ is the unique negative eigenvalue for $-\Delta_\gamma$. Suppose $\lambda>0$ such that $\lambda\neq \mu_\gamma$ and for $-\lambda^2$  there exists $\psi_0\in D(-\Delta_\gamma)-\{0\}$ satisfying $-\Delta_\gamma\psi_0=-\lambda^2\psi_0$. Define
$$
p_\gamma(x)=[(-\Delta+\lambda^2)^{-1}\psi_0](x)-\frac{1}{4\pi^2}\frac{2\gamma \lambda}{\gamma \coth(\lambda\pi)+2\lambda}\langle \psi_0, \overline{J_{i\lambda}}\rangle J_{i\lambda}(x).
$$
Then as in the proof of  Theorem \ref{resol2},  $p_\gamma\in D(-\Delta_\gamma)$ and $[(-\Delta_\gamma+\lambda^2)p_\gamma](x)=\psi_0(x)$ for $x\in (-\pi,\pi)-\{0\}$. Hence,
$$
\|\psi_0\|^2=\langle (-\Delta_\gamma+\lambda^2)p_\gamma,\psi_0\rangle=
\langle p_\gamma, (-\Delta_\gamma+\lambda^2)\psi_0\rangle=0,
$$
which is a contradiction.

\item For $\gamma> 0$, $-\Delta_\gamma$ has no negative eigenvalues. Suppose $\lambda<0$, $\zeta\neq 0$ and $-\Delta_\gamma\zeta=\lambda\zeta$. Then, from integration by parts 
\begin{equation}\label{resol10}
\begin{aligned}
\lambda\|\zeta\|^2=\int_{-\pi}^{\pi} \zeta(x) & (-\Delta_\gamma\zeta(x))dx=\lim_{\epsilon\downarrow 0}\int_{-\pi}^{-\epsilon} \zeta(x)(-\zeta''(x))dx\\
&+\lim_{\epsilon\downarrow 0}\int_{\epsilon}^{\pi} \zeta(x)(-\zeta''(x))dx=\gamma \zeta^2(0)+\|\zeta'\|^2.
\end{aligned}
\end{equation}
From \eqref{resol10} we obtain $\zeta(x)=0$ for all $x\in\mathbb R$, which is a contradiction.

\item Next we show that zero is not eigenvalue of $-\Delta_\gamma$ for $\gamma >0$. Indeed, let $\zeta\in D(-\Delta_\gamma)$, $\zeta\neq 0$, and $-\Delta_\gamma\zeta=0$.  Then zero will be the first eigenvalue of the self-adjoint operator $-\Delta_\gamma$ and so it is simple. Moreover,  $\zeta$ can be choose  as being an even positive function. Then $\zeta$ is symmetric with regard to the line $\xi=\pi$ and therefore $\zeta'(\pm \pi)=0$ (recall $\zeta\in H^2((2n\pi,2(n+1)\pi))$, $n\in\mathbb Z$).  Hence, from integration by parts 
\begin{equation}\label{resol10a}
0=\int_{-\pi}^{\pi} \zeta(x)(-\Delta_\gamma\zeta(x))dx=\gamma \zeta^2(0)+\|\zeta'\|^2,
\end{equation}
which implies  $\zeta(x)=0$ for all $x\in\mathbb R$.

\item The eigenvalues  $\kappa^2_j$ and  $k^2_j$ satisfying the relations \eqref{resol8a} and \eqref{resol8b} respectively, for $j\geqq 1$, are simple.  The proof follows the ideas of how to obtain the function $J_k$ in \eqref{inker}. More exactly, if $f$ satisfies $-\Delta_\gamma f=k^2_j f$ then there is $\beta\in \mathbb R$ such that $f=\beta J_{k_j}.$

\item From  formula in \eqref{resolinf} it follows $\coth(ik\pi)=0$ if and only if $k\in \mathbb R$ and $k$, $n\in \mathbb Z$. Therefore, since for $j\geqq 1$, $j\in \mathbb N$, $\sin(j x)\in D(A^0)\subset D(-\Delta_{+\infty})$ and for all $x\in\mathbb R$
\begin{equation*}
-\Delta_{+\infty} \sin(j x)=A^0(\sin(j x))=-\frac{d^2}{dx^2}\sin(j x)=j^2\sin(j x),
\end{equation*}
the nondegeneracy  of the eigenvalues is immediate.
\end{enumerate}
The proof of the Theorem is completed.
\end{proof}

{\bf{Remarks:}} 

\begin{enumerate}

\item From the formula for the resolvent in \eqref{resol2a} we obtain for $\gamma<0$ the explicit structure of the residuum at $k$ satisfying $2ik=-\gamma \coth(ik\pi)$.

\item From the definition of the domain $D(-\Delta_\gamma)$, $\gamma \neq 0$,  the only periodic constant function in this set is the zero function.

\item We can  give a general proof of  that zero is not eigenvalue of $-\Delta_\gamma$ : Suppose $f\in  D(-\Delta_\gamma)-\{0\}$ such that $-\Delta_\gamma f=0$. Then $f''(x)=0$ for all $x\in (0,2\pi)$, hence since $f$ is periodic we need to have $f\equiv r$, $r$ a real constant. So, from the jump condition $r=0$ which is a contradiction.

\item From the min-max principle we obtain that for $\gamma <0$
\begin{equation}\label{resol11}
\lambda=\inf \Big\{\|v_x\|^2+\gamma \int \delta(x)|v(x)|^2dx: \|v\|=1, v\in H^1_{per}\Big\}
\end{equation}
is given by $\lambda=-\mu_\gamma^2$ and the corresponding positive eigenfunction is $\psi_\gamma$ in \eqref{resol8}.

\end{enumerate}

\section{Global Well-Posedness in $H^1_{per}$}

Our notion of well-posedness for the equation NLS-$\delta$ in an arbitrary functional space $Y$ is the existence, uniqueness, persistence property (i.e. the solution describes a continuous curve in $Y$ whenever $u_0\in Y$) and the continuous dependence of the solution upon the data. The following proposition is concerned with the well-posedness of equation \eqref{nls} in $H^1_{per}([0, 2L])$.

\begin{prop}\label{lwp1}  For any $u_0\in H^1_{per}([0, 2L])$, there exists $T>0$ and a unique solution $u\in C([-T,T]; H^1_{per}([0, 2L]))\cap C^1([-T,T]; H^{-1}_{per}([0, 2L]))$ of \eqref{nls}, such that $u(0)=u_0$.  For each $T_0\in (0,T)$ the mapping
$$
u_0\in H^1_{per}([0, 2L]) \to u\in C([-T_0,T_0]; H^1_{per}([0, 2L]))
$$
is continuous. Moreover, since $u$ satisfies the conservation of the energy and the charge defined in \eqref{conse}, namely,
$$
E(u(t))=E(u_0),\quad Q(u(t))=Q(u_0),
$$
for all $t\in [0,T)$, we can choose $T=+\infty$.

If an initial data $u_0$ is even the solution $u(t)$ is also even.
\end{prop}

\begin{proof} We apply Theorem 3.7.1 of \cite{Ca} to our problem. Indeed,  from  Theorem \ref{resol5b} we have $-\Delta_{-Z}\geqq -\beta$, where $\beta=\mu^2_{-Z}$, if $Z>0$ and $\beta=0$ if $Z<0$. So,  for the self-adjoint operator $\mathcal A\equiv \Delta_{-Z}-\beta$  on $X=L^2_{per}([0, 2L])$ with domain
$\mathcal D(\mathcal A)=\mathcal D(-\Delta_{-Z})$ we have $\mathcal A\leqq 0$. Moreover, in our situation, we may take the space  $X_{\mathcal A}=H^1_{per}([0, 2L])$ with norm 
$$
\|u\|^2_{X_{\mathcal A}}=\|u_x\|^2+ (\beta+1)\|u\|^2-Z|u(0)|^2,
$$
which is equivalent to $H^1_{per}([0, 2L])$ norm (see \eqref{resol11}). So, it is very easy to see that the uniqueness of solutions and the conditions (3.7.1), (3.7.3)-(3.7.6) in \cite{Ca} hold choosing $r=\rho=2$. Finally, the condition (3.7.2) in \cite{Ca} with $p=2$ is satisfied because of $\mathcal A$ is a self-adjoint operator on $L^2_{per}([0, 2L])$.
\end{proof}

\section{Periodic travelling-wave for NLS-$\delta$}

In this section we construct   positive periodic solutions for the elliptic equation  \eqref{p2} such that the conditions in \eqref{p3} are satisfied. These solutions belong to the domain of the operator $-\frac{d^2}{dx^2}-Z \delta$, $Z\neq 0$. Our analysis is based in the theory of elliptic integral,  the theory of Jacobi elliptic functions and the implicit function theorem. 

\subsection{The quadrature method}

We start by writing  \eqref{p3}-(3) in quadratic form. Indeed, for $\varphi=\varphi_{\omega, Z}$ and $x\neq \pm 2n L$ we obtain
 \begin{equation}\label{s1}
 [\varphi'(x)]^2=\frac12[-\varphi^4(x)+2\omega \varphi^2(x)+4B_{\varphi}] \equiv \frac12 F(\varphi(x)),
\end{equation}
where  $F(t)=-t^4+2\omega t^2+4B_{\varphi}$ and $B_{\varphi}$ is a integration constant. We factor $F(\cdot)$ as 
 \begin{equation}\label{s1a}
F(\varphi)=(\eta_1^2-\varphi^2)(\varphi^2-\eta_2^2)=2[\varphi']^2,
\end{equation}
where $\eta_1, \eta_2$ are the positive zeros of the polynomial
$F$. We assume without loss of
generality that $\eta_1>\eta_2>0$.  So, $\eta_2\leqq
\varphi(\xi)\leqq\eta_1$  and 
\begin{equation}\label{s2}
2\omega=\eta_1^2+\eta^2_2,\;\;\;\;\;\;
4B_{\varphi}=-\eta_1^2\eta_2^2.
\end{equation}

Next, since $\varphi$ is continuous one has 
\begin{equation}\label{s3}
[\varphi'(0+)]^2=\frac12 F(\varphi(0))\quad{\rm{and}}\quad [\varphi'(0-)]^2=\frac12 F(\varphi (0)).
\end{equation}
Then $|\varphi'(0+)|=|\varphi'(0-)|$, which as we will show below implies that $\varphi'(0+)=-\varphi'(0-)$, and so from \eqref{p3}-(4) 
\begin{equation}\label{s4}
\varphi'(0+)=-\frac{Z}{2}\varphi(0).
\end{equation}
The case  $\varphi'(0+)=\varphi'(0-)$ can not happen. Indeed, from \eqref{p3}-(4) it follows $\varphi(0)=0$ and so $\varphi'(0)$ exists. Therefore from \eqref{s1a} $[\varphi'(0)]^2=- \eta_1^2\eta_2^2/2$
which is a contradiction.

Next, we obtain restrictions on the value of $\varphi(0)$. From \eqref{s1} and \eqref{s4} we need to have 
\begin{equation}\label{s5}
\frac{Z^2}{4} \varphi^2(0)=\frac12 F(\varphi(0))>0,
\end{equation}
and so $\eta_1>\varphi(0)>\eta_2$.  Next, since  $\max_{t\in\mathbb R} F(t)=\omega^2+4B_{\varphi}$ (which is attained for $t>0$ in $t_0=\sqrt{\omega}$), we obtain the condition
\begin{equation}\label{s6}
\frac{Z^2}{4} \varphi^2(0)\leqq \frac{\omega^2+4B_{\varphi}}{2}=\frac{(\omega-\eta_2^2)^2}{2},
\end{equation}
and from \eqref{s5}
\begin{equation}\label{s7}
\varphi^2(0)= \frac{-(2\omega-\frac{Z^2}{2})\pm\sqrt{(2\omega-\frac{Z^2}{2})^2+
16B_{\varphi}}}{-2}.
\end{equation}
Since $\varphi(0)\in \mathbb R$ we need to have $(2\omega-\frac{Z^2}{2})^2+
16B_{\varphi}>0$. We start by considering the case of sign ``$\;-\;$'' in the square root in \eqref{s7}, then:
\begin{enumerate}
\item For $2\omega-\frac{Z^2}{2}>0$, it follows from \eqref{s2} that $
(2\omega-\frac{Z^2}{2})^2>4\eta_1^2\eta_2^2$ and so $\eta_1^2+\eta_2^2-2\eta_1\eta_2>\frac{Z^2}{2}$. Hence,
\begin{equation}\label{s8}
\eta_1-\eta_2>\frac{|Z|}{\sqrt{2}}.
\end{equation}

\item From \eqref{s7} we have  as $Z\to 0$ the asymptotic behavior
\begin{equation}\label{s9}
\varphi^2(0)\to \frac{-2\omega-\sqrt{4\omega^2+16B_{\varphi}}}{-2}=\frac{-\eta_1^2-\eta_2^2-(\eta_1^2-\eta_2^2)}{-2}=\eta_1^2.
\end{equation}

\item For $2\omega-\frac{Z^2}{2}<0$  we obtain from \eqref{s7} that $16B_{\varphi}>0$, which is not possible from \eqref{s2}.
\end{enumerate}
Now, we consider the case of sign ``$\;+\;$'' in the square root in \eqref{s7}, then:
\begin{enumerate}
\item For $2\omega-\frac{Z^2}{2}<0$ we have  $\varphi^2(0)<0$, which is a contradiction.

\item For $2\omega-\frac{Z^2}{2}>0$  we still have relation \eqref{s8}, but as $Z\to 0$ we obtain $\varphi^2(0)\to \eta_2^2$.
\end{enumerate}

\vskip0.1in

We are interested only in  the  sign ``$\;-\;$'' in \eqref{s7} for  our stability theory.

\subsection{Profile of positive periodic peaks for $Z>0$}

Next we are interested in finding a even periodic profile  solution, $\phi_{\omega,Z}$ for \eqref{p2}  such that the peaks will be happen  in $\pm 2nL$, $n\in \mathbb Z$, $\eta_1>\phi_{\omega,Z}(0)\geqq \phi_{\omega,Z}(\xi)\geqq \eta_2$ for all $\xi$, and
$$
\lim_{Z\to 0^+} \phi_{\omega,Z}=\phi_{\omega,0},
$$
where $\phi_{\omega,0}$ is the dnoidal traveling wave defined in \eqref{dnoidal2}.  Without loss of generality we can assume $2L=1$.

We start our analysis by considering an additional variable $\psi$ via the relation
\begin{equation}\label{s10}
\phi^2(\xi)=\theta-\alpha \sin^2\psi(\xi),
\end{equation}
with $\phi=\phi_{\omega,Z}$ and $\theta$, $\alpha$, constants to be choosen later. So for $\xi\neq \pm n$, $n\in \mathbb Z$, we obtain the equality
\begin{equation}\label{s11}
2\phi(\xi)\phi'(\xi)=-2\alpha \psi'(\xi)\sin\psi(\xi) \cos\psi(\xi).
\end{equation}
Therefore, from \eqref{p3}-(4) we get the identity
\begin{equation}\label{s12}
-2Z \phi^2(0)=\alpha  [\psi'(0-)-\psi'(0+)]\sin 2\psi(0),
\end{equation}
and so the phase-function $\psi$ satisfies the conditions
\begin{equation}\label{s13}
\psi'(0+)-\psi'(0-)\neq 0\quad{\rm{and}}\quad \psi(0)\neq\frac{k\pi}{2}, \;k\in \mathbb Z.
\end{equation}

From \eqref{s11} and \eqref{s1a} it follows that 
\begin{equation}\label{s14}
(\psi')^2\alpha^2\sin^2 \psi\cos^2 \psi=\frac12(\theta-\alpha \sin^2 \psi)(\eta_1^2-\theta+\alpha \sin^2 \psi)(\theta-\eta_2^2-\alpha \sin^2 \psi).
\end{equation}
By choosing $\theta=\alpha=\eta_1^2$ we have $\phi^2(\xi)=\eta^2_1 \cos^2\psi(\xi)$, and obtain the  ordinary differential equation
\begin{equation}\label{s17}
[\psi'(\xi)]^2=\frac12(\eta^2_1-\eta^2_2)(1-\eta^2\sin^2\psi(\xi))
\end{equation}
for  $\xi\neq \pm n$, $n\in \mathbb Z$, and
\begin{equation}\label{s16}
\eta^2\equiv \frac{\eta_1^2}{\eta_1^2-\eta_2^2}>1.
\end{equation}

From a basic analysis (see Appendix) one has  that  $\psi'(\xi)=0$ if and only if $\xi=s$, where $s$ is the unique point  in $(0,1)$ s.t. $\phi(s)=\eta_2$. Moreover, $\psi$ has minimal period 1 and  if $\phi$ is even then $\psi$ is also even.  Then we obtain that $\psi'(\xi)>0$ for $\xi\in (0,s)$ which implies that we have a peak in  $\pm n\mathbb Z$ for $\psi$ in the form ``$\vee$''.

Now for $\ell^2=(\eta_1^2-\eta_2^2)/2$, it follows from \eqref{s17} and from the behavior of $\psi$ that
\begin{equation}\label{s20}
\psi'(\xi)=-\ell \sqrt{1-\eta^2\sin^2\psi(\xi)}\qquad \text{for}\;\;\xi\in (-s, 0),
\end{equation}
and so for
\begin{equation}\label{s21}
F(\xi)\equiv -\int_{\psi(0)}^{\psi(\xi)}\frac{dt}{\sqrt{1-\eta^2\sin^2 t}}
\end{equation}
we have  $F(\xi)=\ell \xi$ for $\xi\in (-s, 0)$. Therefore from the equality
\begin{equation}\label{s22}
\ell \xi=-\int_{0}^{\psi(\xi)}\frac{dt}{\sqrt{1-\eta^2\sin^2 t}}+\int_{0}^{\psi(0)}\frac{dt}{\sqrt{1-\eta^2\sin^2 t}}
\end{equation}
and from the relations $k=1/\eta$ and $\sin \beta(\xi)\equiv\eta \sin \psi(\xi)$,  we obtain from the  theory of Jacobian elliptic functions (see \cite{byrd}) that
\begin{equation}\label{s23}
\begin{aligned}
\ell \xi&= -k \text{sn}^{-1}(\sin \beta(\xi);k) +k \text{sn}^{-1}(\sin \beta(0);k)\\
&= -k \text{sn}^{-1}(\eta \sin \psi(\xi);k) +k \text{sn}^{-1}(\eta \sin \psi(0);k).
\end{aligned}
\end{equation}
Now from \eqref{s23} one has
\begin{equation}\label{s24}
a\equiv \text{sn}^{-1}(\eta \sin \psi(0);k),\qquad [\text{sn}(a;k)=\eta \sin \psi(0)]
\end{equation}
and consequently  the formula
\begin{equation}\label{s25}
\text{sn}\Big(-\frac{\ell\xi}{k}+a;k\Big)=\eta \sin\psi(\xi),\quad{\rm{for}}\;\xi\in (-s,0).
\end{equation}

Next we obtain the exactly value of $a$. Indeed, from   the identities
$$
\text{sn}^{-1}(y;k)=\text{cn}^{-1}(\sqrt{1-y^2};k)=\text{dn}^{-1}(\sqrt{1-k^2y^2};k)
$$
follows,
\begin{equation}\label{s25a}
a=\text{dn}^{-1}\Big(\sqrt{1-k^2\eta^2\sin^2\psi(0)};k\Big)=\text{dn}^{-1}(\cos \psi(0);k)=\text{dn}^{-1}\Big(\frac{\phi(0)}{\eta_1};k\Big).
\end{equation}
Here $\phi(0)$ is given by the formula in the right hand side of \eqref{s7}. The shift-value $a$ depends of the values of $Z$ and $\omega$. Moreover, since $1>\phi(0)/\eta_1>k'\equiv \sqrt{1-k^2}$, $1\geqq\text{dn}(x)\geqq k'$, for all $x\in\mathbb R$ and  $k'=\text{dn}(K)$, it follows that $a$ is well-defined and $a\in [0, K]$. Hence, from  \eqref{s25} and  \eqref{s10} we obtain the profile
\begin{equation}\label{s26}
\phi(\xi)=\eta_1\text{dn}(-\eta \ell \xi +a;k)=\eta_1\text{dn}(-\frac{\eta_1}{\sqrt{2}} \xi +a;k), \quad{\rm{for}}\;\xi\in (-s,0).
\end{equation}

Similarly, from \eqref{s17}  
\begin{equation}\label{s27}
\phi(\xi)=\eta_1\text{dn}(-\frac{\eta_1}{\sqrt{2}} \xi +a;k), \quad{\rm{for}}\;\xi\in (0,s).
\end{equation}
Therefore, one obtains the  peak-function
\begin{equation}\label{s28}
\phi(\xi)=\phi(\xi; \eta_1, \eta_2, Z)=\eta_1\text{dn}\Big(\frac{\eta_1}{\sqrt{2}} |\xi| +a;k\Big), \quad{\rm{for}}\;\xi\in (-s,s),
\end{equation}
where the shift $a$ is given by \eqref{s25a} and $\eta_1, \eta_2$ satisfy

\begin{equation}\label{s29}
\begin{cases}
\begin{aligned}
&k^2=\frac{\eta_1^2-\eta_2^2}{\eta_1^2},\quad \eta_1^2+\eta_2^2=2\omega,\\
&\quad  0<\eta_2<\eta_1,\quad {\rm{and}}\quad \eta_1-\eta_2>\frac{|Z|}{\sqrt{2}}.
\end{aligned}
\end{cases}
\end{equation}

Next we shall  determine the exactly value of $s>0$ such that $\phi(s)=\eta_2$. From \eqref{s28} it follows
$$
\text{dn}^2\Big(\frac{\eta_1}{\sqrt{2}} s +a;k\Big)=\frac{\eta_2^2}{\eta_1^2}=1-k^2\equiv k'^2.
$$
Then,  since $dn(K)=k'$ and $\text{dn}$ has a minimal period $2K$ one has $\frac{\eta_1}{\sqrt{2}} s +a=(2n+1)K(k)$ for $n\in \mathbb Z$ and so one can choose
\begin{equation}\label{s30}
s\equiv \frac{\sqrt{2}}{\eta_1}(K-a).
\end{equation}
We note that if $Z\to 0^+$ then $\phi(0)\to \eta_1$, and so $a\to 0$. Then we conclude that $s(Z)\to \frac{\sqrt{2}}{\eta_1}K$ as $Z\to 0^+$. Lastly, relation
\begin{equation}\label{s31}
\phi(2s)=\eta_1\text{dn}\Big(\frac{\eta_1}{\sqrt{2}} s +K\Big)=\eta_1\text{dn}(2K-a)=\eta_1\text{dn}(-a)=\eta_1\text{dn}(a)=\phi(0)
\end{equation}
implies that the profile $\phi$ in \eqref{s28} can be extend to  all the line as a continuous periodic function satisfying the conditions in \eqref{p3} with a  minimal period $2s$  (see Figure 3). In the next subsection 5.4 we will show that it is possible to choose $s=1/2$.

\begin{figure}[!htb]
\includegraphics[angle=0,scale=0.2]{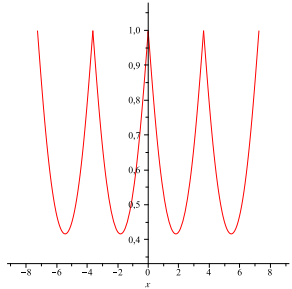}
\caption{Profile of the periodic dnoidal-peak $\phi$ in \eqref{s28}.}
\end{figure}

From the above analysis we have, at least formally,  that
$$
\lim_{Z\to 0^+} \phi_{\omega,Z}(x)=\phi_{\omega,0}(x),
$$
where $\phi_{\omega,0}$ is the dnoidal traveling wave defined in \eqref{dnoidal2}. The last equality must be understood in the following sense:  for $\omega>\frac{\pi^2}{2L^2}$ fixed and $x\in (0, 2L)$, there is a $\delta>0$ such that for $Z\in(0,\delta)$ and $Z^2/4<\omega$ we have that the family of periodic-peaks  $\phi_{\omega,Z}$, with minimal period $2s$, are all defined in $x$. We note that this type of convergence is not convenient for our purposes, because the period of  $\phi_{\omega,Z}$ is changing .

\subsection{Positive periodic peaks for $Z<0$}

We shall find a even  periodic-peak, $\zeta_{\omega,Z}$, with peaks   in $\pm 2nL$, $n\in \mathbb Z$,  $\eta_1> \zeta(0)>\eta_2$, $
\eta_1\geqq \zeta(\xi)\geqq \eta_2$ for all $\xi \neq \pm 2nL$, and
$$
\lim_{Z\to 0^-} \zeta_{\omega,Z}=\phi_{\omega,0},
$$
where $\phi_{\omega,0}$ is the dnoidal traveling wave defined in \eqref{dnoidal2}. Suppose $2L=1$. Next we consider $f$ via the relation
\begin{equation}\label{sne1}
\zeta^2(\xi)=\eta_1^2 \sin^2f(\xi)\geqq\eta_2^2.
\end{equation}
So for $\xi\neq \pm n$, $n\in \mathbb Z$, 
\begin{equation}\label{sne2}
2\zeta(\xi)\zeta'(\xi)=\eta_1^2 f'(\xi)\sin2f(\xi).
\end{equation}
Therefore, from \eqref{p3}-(4) we get 
\begin{equation}\label{sne3}
-2Z \zeta^2(0)=\eta_1^2  [f'(0+)-f'(0-)]\sin 2f(0),
\end{equation}
and so the phase-function $f$ satisfies the conditions
\begin{equation}\label{sne4}
f'(0+)-f'(0-)\neq 0\quad{\rm{and}}\quad f(0)\neq\frac{k\pi}{2}, \;k\in \mathbb Z.
\end{equation}
From \eqref{sne2} and \eqref{s1a} one obtains
\begin{equation}\label{sne5}
[f']^2=\frac12\eta_2^2\Big(\frac{\eta_1^2}{\eta_2^2} \sin^2f - 1\Big).
\end{equation}
For $Z<0$ we know that $\zeta'(0+)>0$, so let  $p\in (0,1)$ be the first value such that $\zeta(p)=\eta_1$ ($p$ is a maximum point for $\zeta$). Therefore, $\sin^2f(p)=1$ implies $f(p)=\frac{\pi}{2}$. Now, we have the following assumptions and behavior for  $f$:

 (a) From \eqref{sne1}, $f(\xi)\neq 0$ for all $\xi$. So, if we suppose $f$ being strictly positive and $f(\xi)\in (0, \frac{\pi}{2}]$ (it suffices  to have $f(0)\in (0, \frac{\pi}{2})$, we obtain that $f$ is periodic with  period 1. 
 
 (b) From \eqref{sne3}, \eqref{sne4}, and \eqref{sne5}  we have $f'(0+)=-f'(0-)$. So, since $Z<0$ it follows $f'(0+)>0$ and for $\xi\in (0, p)$, $f'(\xi)>0$. By evenness $f'(\xi)<0$ for $\xi\in (-p, 0)$ and therefore $f$ has a peak in zero in the form ``$\vee$''.

\vskip0.1in

Next we build a periodic peak in the form  ``$\vee$''. From  \eqref{sne5} it follows for $\xi\in(-p, 0)$ that
\begin{equation}\label{sne6}
f'(\xi)=-\frac{\eta_2}{\sqrt{2}}\sqrt{a^2\sin^2f(\xi)-1},
\end{equation}
with $a^2=\eta_1^2/{\eta_2^2}$ and $\frac{\pi}{2}>f(\xi)>f(0)\geqq \sin^{-1}(\eta_2/{\eta_1})$. Next, define for $\xi\in (-p,0)$
\begin{equation}\label{sne7}
G(\xi)\equiv -\int_{f(0)}^{f(\xi)}\frac{d\nu}{\sqrt{a^2\sin^2\nu -1}},
\end{equation}
then from \eqref{sne6}  it follows that  $G(\xi)=\frac{\eta_2}{\sqrt{2}}\xi$ for $\xi\in(-p, 0)$.  Therefore,  from the equality
\begin{equation}\label{sne8}
 \frac{\eta_2}{\sqrt{2}}\xi=-\int_{f(0)}^{\frac{\pi}{2}}\frac{d\nu}{\sqrt{a^2\sin^2 \nu -1}} +\int_{f(\xi)}^{\frac{\pi}{2}}\frac{d\nu}{\sqrt{a^2\sin^2 \nu -1}},
\end{equation}
and from Byrd\&Friedman (pg. 167), we obtain  that the relations
\begin{equation}\label{sne9}
k^2=\frac{\eta_1^2-\eta_2^2}{\eta_1^2}, \quad \sin(\tau(0))=\frac{1}{k}\cos f(0), \quad\text{and}\quad \sin(\tau(\xi))=\frac{1}{k}\cos f(\xi),
\end{equation}
imply
\begin{equation}\label{sne10}
 \frac{\eta_2}{\sqrt{2}}\xi= -\frac{\eta_2}{\eta_1} sn^{-1}(\sin \tau(0);k)+\frac{\eta_2}{\eta_1}sn^{-1}(\sin \tau(\xi);k).
\end{equation}
Moreover, since $Z^2\zeta^2(0)/4=F(\zeta(0))/2$, we have from \eqref{s5} and \eqref{s7} that $\zeta(0)=\phi(0)$ and so  \eqref{s25a} implies
\begin{equation}\label{sne10a}
b\equiv sn^{-1}(\sin \tau(0);k)=dn^{-1}(\sin f(0);k)=dn^{-1}\Big(\frac{\phi(0)}{{\eta_1}};k\Big)=a.
\end{equation}
Then from \eqref{sne10} we obtain
\begin{equation}\label{sne11}
 k^2sn^2\Big(\frac{\eta_1}{\sqrt{2}}\xi +a;k\Big)= k^2\sin^2 \tau(\xi)=\cos^2 f(\xi)=1-\sin^2 f(\xi),
\end{equation}
and so from \eqref{sne1} we find the profile
\begin{equation}\label{sne12}
\zeta(\xi)=\eta_1dn\Big(\frac{\eta_1}{\sqrt{2}}\xi +a;k\Big),\quad\text{for}\;\xi\in(-p,0).
\end{equation}
Similarly, we obtain 
\begin{equation}\label{sne13}
\zeta(\xi)=\eta_1dn\Big(-\frac{\eta_1}{\sqrt{2}}\xi +a;k\Big),\quad\text{for}\;\xi\in(0,p).
\end{equation}
Then, we obtain the following  peak-function defined initially for $\xi\in(-p,p)$ and with the ``$\vee$''-profile in zero,
\begin{equation}\label{sne14}
\zeta(\xi)=\eta_1dn\Big(\frac{\eta_1}{\sqrt{2}}|\xi| -a;k\Big).
\end{equation}

Since $\zeta(p)=\eta_1$ follows that $dn\Big(\frac{\eta_1}{\sqrt{2}}p -a;k\Big)=1$ and so $p=p(Z) =\sqrt{2}a/{\eta_1}$ (the first one such that $\zeta'(p)=0$). We note that if $Z\to 0^-$ then $\zeta(0)\to \eta_1$ and so $p(Z)\to 0$ as $Z\to 0^-$. Now, if we see the profile $\zeta$ in \eqref{sne14} defined in all the line, we have that   for 
\begin{equation}\label{sne15}
p_0\equiv \frac{\sqrt{2}}{\eta_1}(K+a),
\end{equation}
follows
\begin{equation}\label{sne16}
\zeta(p_0)=\eta_1dn\Big(\frac{\eta_1}{\sqrt{2}}p_0 -a;k\Big)=\eta_1dn(K;k)=\eta_1k'=\eta_2.
\end{equation}
Then $\zeta(\pm p_0)=\eta_2$ ($\zeta'(\pm p_0)=0$). 

Therefore we can build a even periodic peak for the NLS-$\delta$ satisfying the conditions in \eqref{p3} with a minimal period $2p_0$, and 
it being the periodicity of the even-profile $\zeta(\xi)$ in \eqref{sne14} with $\xi\in [-p_0,p_0]$ (see Figure 4). In the next subsection 5.5 we will show that it is possible to choose $p_0=1/2$.

\begin{figure}[!htb]
\centering
\includegraphics[angle=0,scale=0.2]{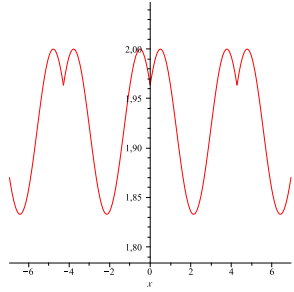}
\caption{Profile of the periodic dnoidal-peak $\zeta$ in \eqref{sne14}.}
\end{figure}

\vskip0.2in

{\bf{Remarks:}}\begin{enumerate}
\item For the \lq\lq convergence'' of the periodic-peak $\phi_{\omega,Z}$ and   $\zeta_{\omega,Z}$ to the solitary wave peak \eqref{ellip1}, with $p=2$  we  consider for a determined parameter ($\eta_2$ is our case) the minimal period $2s$ in \eqref{s30} or $2p_0$ in \eqref{sne15} sufficiently large. Indeed, from \eqref{s29} and \eqref{sne9} we obtain for all $Z$ that $k^2(\eta_2)\to 1$ for $\eta_2\to 0$  and so $K(k(\eta_2))\to +\infty$. From \eqref{s7} (with ``-'') and $2\sqrt{\omega}>|Z|$ it follows that 
\begin{equation}\label{s32}
\phi^2_{\omega,Z}(0)=\zeta^2_{\omega,Z}(0)\to 2\omega-\frac{Z^2}{2}\qquad\text{as}\;\eta_2\to 0.
\end{equation}
So, in the case of $\phi_{\omega, Z}$ (a similar result is obtained for $\zeta_{\omega, Z}$) that  \eqref{s32} implies 
\begin{equation}\label{s35}
\sin^2\psi(0)=1-\frac{\phi_{\omega, Z}^2(0)}{\eta_1^2}\to \frac{Z^2}{4\omega}\quad{\rm{as}}\,\, \eta_2\to 0,
\end{equation}
since $\eta_1^2\to 2\omega$ as $\eta_2\to 0$. Then, combining \eqref{s24}, \eqref{s35}, and $sn(\cdot;1)=\tanh(\cdot)$ one gets that 
$$
a(\eta_2; Z)\to \tanh^{-1}\Big(\frac{Z}{2\sqrt{\omega}}\Big)\quad{\rm{as}}\,\, \eta_2\to 0.
$$
Lastly, since $\text{dn}(\cdot;1)=\text{sech}(\cdot)$ we obtain the convergence (uniformly on compact-set) 
\begin{equation}\label{s36}
\phi_{\omega, Z}(\xi)\to \phi_{\omega, Z,2}(\xi), \quad{\rm{as}}\,\, \eta_2\to 0.
\end{equation}


\item Since $\phi_{\omega, Z, 2}^2(0)=2\omega-\frac{Z^2}{2}$ and the value $\phi_{\omega, Z}^2(0)=\zeta_{\omega, Z}^2(0)$ is a increasing function of $\eta_2$, it follows  from \eqref{s32} that the peaks associated to $\phi_{\omega, Z}, \zeta_{\omega, Z}$ resemble that of $\phi_{\omega, Z, 2}$ in a neighborhood of zero and approximating to it  for below.

\item As  $Z\to 0^-$ one has  $\zeta_{\omega, Z}(0)\to \eta_1$ and so $a\to 0$. Thus, from \eqref{sne15}, $p_0(Z)\to \frac{\sqrt{2}}{\eta_1}K^+$ as $Z\to 0^-$. Hence 
$$
\lim_{Z\to 0^-} \zeta_{\omega,Z}(x)=\phi_{\omega,0}(x),
$$
where $\phi_{\omega,0}$ is the dnoidal traveling wave defined in \eqref{dnoidal2}.

\end{enumerate}

\subsection{Dnoidal-peak solutions to the NLS-$\delta$ with an arbitrary minimal period}

In subsections 5.2 and 5.3 we found  dnoidal-peak profiles \eqref{s28} and \eqref{sne14} with a minimal period $2s$ and $2p_0$. 
Next we shall see that the equality $s=L$ can be obtained by any {\it a priori} $L$.  In the analysis below we consider the case $Z>0$, but a similar result can be established for $Z<0$.

We start by  defining the general notations to be used in the next subsections. For $4\omega>Z^2$ it follows from \eqref{s29} that for all  $Z$, 
\begin{equation}\label{s37}
0<\eta_2<\theta(\omega,Z)<\sqrt{\omega}<\lambda(\omega,Z)<\eta_1<\sqrt{2\omega}
\end{equation}
with
\begin{equation}\label{s38}
\theta(\omega,Z)=-\frac{\sqrt{2}}{4}|Z|+\sqrt{\omega-\frac18 Z^2}\quad{\rm{and}}\quad \lambda(\omega,Z)=\frac{\sqrt{2}}{4}|Z|+\sqrt{\omega-\frac18 Z^2}.
\end{equation}
For $\eta\in (0,\theta(\omega,Z))$ we define the functions:
\begin{equation}\label{s39}
k^2(\eta, \omega)=\frac{2\omega-2\eta^2}{2\omega-\eta^2},
\end{equation}
and 
\begin{equation}\label{s40}
T_{-}(\eta, \omega, Z)=\frac{2\sqrt{2}}{\sqrt{2\omega-\eta^2}}[K(k(\eta, \omega))-a(\eta, \omega, Z)]
\end{equation}
where 
\begin{equation}\label{s41}
a(\eta, \omega, Z)=\text{dn}^{-1}\Big(\frac{\Phi(\eta, \omega, Z)}{\sqrt{2\omega-\eta^2}};k(\eta, \omega)\Big),
\end{equation}
with $\Phi(\eta, \omega, Z)$  defined by (see \eqref{s7}) 
\begin{equation}\label{s42}
\Phi^2(\eta, \omega,  Z)= \frac{(2\omega-\frac{Z^2}{2})+\sqrt{(2\omega-\frac{Z^2}{2})^2-4\eta^2(2\omega-\eta^2)}}{2}.
\end{equation}
We note that the functions $a$ and $\Phi$ defined above are independent of the sign of $Z$. We will denote  them by  $a(\eta)$, $\Phi(\eta)$ or   $a(\eta,\omega)$, $\Phi(\eta, \omega)$ depending of the context. 

\vskip0.2in

{\bf{Remark:}} For $\eta\in (0, \theta(\omega,Z))$ we obtain the  condition (\ref{s6}), namely, $\frac{Z^2}{4}\Phi^2(\eta, \omega,  Z)\leqq \frac{(\omega-\eta^2)^2}{2}$.

\begin{teo}\label{funda1}  For $Z\neq 0$ and $\omega>Z^2/4$ fixed, the mappings for $\eta_2\in (0,\theta(\omega,Z))$
\begin{equation}\label{s42a}
\eta_2 \to a(\eta_2),\quad \eta_2 \to \Phi(\eta_2), \quad{\rm{and}}\quad \eta_2\to \frac{\Phi(\eta_2)}{\sqrt{2\omega-\eta_2^2}}
\end{equation}
are  well defined. Moreover, they are strictly increasing,  strictly decreasing and strictly decreasing functions respectively. Also, one has that
\begin{equation}\label{s43}
\lim_{\eta_2\to 0}T_-(\eta_2)=+\infty,
\end{equation}
and 
\begin{equation}\label{s47}
\lim_{\eta_2\to \theta}T_-(\eta_2)=\frac{2\sqrt{2}}{\lambda(\omega,Z)}[K(k_0)-a_0]\equiv T_0(\omega,Z),
\end{equation}
where
\begin{equation}\label{s44}
k_0^2\equiv k_0^2(\omega,Z)=\frac{\sqrt{2}|Z|\sqrt{\omega-\frac{Z^2}{8}}}{\lambda^2(\omega,Z)},
\end{equation}
and  $a_0\equiv a_0(\omega,Z)\in (0, K(k_0))$ is defined by
\begin{equation}\label{s46}
\text{dn}(a_0;k_0)=\frac{\sqrt{\omega-\frac{Z^2}{4}}}{\lambda(\omega,Z)}.
\end{equation}

Finally,  the mapping  $\eta_2\in (0, \theta(\omega, Z))\to T_-(\eta_2)$ is a strictly decreasing function and so $T_-(\eta_2)\in (T_0(\omega, Z),+\infty)$. Moreover, for $\eta_2\in (0, \theta(\omega, Z))$ it follows that 
$$
\sqrt{2\omega-\eta_2^2}-\eta_2>\frac{|Z|}{\sqrt{2}}.
$$
\end{teo}

\begin{proof} The proof of the Theorem is immediate. In fact, the  inequality 
$$
0>\frac{Z^2}{2}-\sqrt{2}|Z|\sqrt{\omega-\frac{Z^2}{8}}>2\eta_2^2-2\omega+\frac{Z^2}{2}
$$
implies that
$$
1>\frac{\Phi(\eta_2)}{\sqrt{2\omega-\eta_2^2}}>\sqrt{1-k^2(\eta_2)}\equiv k'(\eta_2),
$$
and so $a$ is well defined.  Now, from \eqref{s39} and \eqref{s42} we have for $\eta_2\to 0$ that $k(\eta_2)\to 1$
and $\Phi^2(\eta_2)\to 2\omega-\frac{Z^2}{2}$. Therefore, 
$$
\alpha=\lim_{\eta_2\to 0}a(\eta_2)<\infty
$$ 
with $\alpha$ satisfying $\text{sech}(\alpha)=\sqrt{1-\frac{Z^2}{4\omega}}$. So, combining  \eqref{s40} and the fact that $K(k(\eta_2))\to +\infty$ as $k(\eta_2)\to 1$ we obtain \eqref{s43}.

Now, for $\eta_2\to \theta(\omega,Z)$ one has  $k^2(\eta_2)\to k_0^2$. Since the mappping $\eta_2\to k^2(\eta_2)$ is strictly decreasing it follows that
\begin{equation}\label{s45}
k(\eta_2)\in (k_0,1), \qquad\text{for}\;\;\eta_2\in (0, \theta(\omega,Z)).
\end{equation}
We note that the condition $\omega>Z^2/4$ implies that the right hand side of \eqref{s46} is bigger than $k_0'\equiv \sqrt{1-k_0^2}$ and so $a_0$ is well-defined.  The above considerations yield the limit in \eqref{s47}. 

From Figure 5 ($Z=1/2$), $\eta_2\in (0, \theta(\omega, Z))\to a(\eta_2)$ is a strictly increasing function. The  decreasing property of the other functions in (\ref{s42a}) follows immediately.

\begin{figure}[!htb] 
\includegraphics[angle=0,scale=0.3]{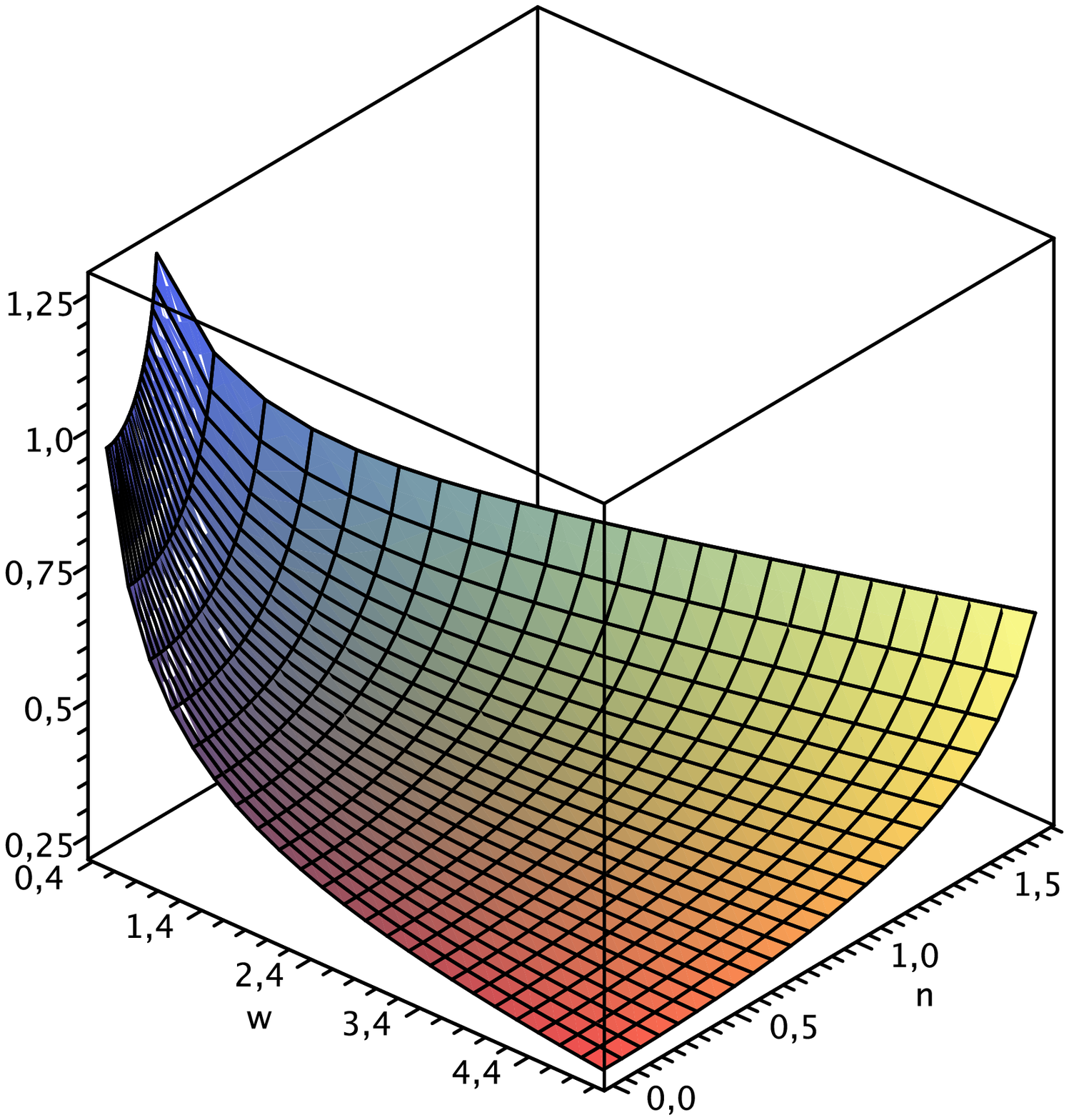}
\caption{graph of  $(\eta,\omega)\to a (\eta,\omega)$.}
\end{figure}

The fact that the mapping  $\eta_2\in (0, \theta(\omega, Z))\to T_-(\eta_2)$ is a strictly decreasing function follows from the analysis in Theorem \ref{s51b} below.
\end{proof}

\vskip0.2in

{\bf{Remark:}} Figure 5 shows that the mapping $\omega \to a(\eta, \omega)$ is a strictly decreasing function. Moreover, $a(\eta, \omega)\to 0$ as $\omega\to +\infty$. This latter can be seen easy from formula (\ref{s41}).

Next we study, for $Z$ fixed, the behavior of the mapping $\omega\in (\frac{Z^2}{4}, +\infty)\to T_0(\omega, Z)$ given in \eqref{s47} (Figure 6 shows a general  profile of $(\omega,Z)\to T_0(\omega,Z)$). From \eqref{s44} one has  for $\omega\to +\infty$ that $k_0^2\to 0$, then $K(k_0)\to \frac{\pi}{2}$ and $a_0\to 0$. So, 
\begin{equation}\label{s48}
\lim_{\omega\to +\infty}T_0(\omega, Z)=0.
\end{equation}
Now since $\beta(\omega, Z)\equiv K(k_0)-a_0$, we have that $\beta(\frac{Z^2}{4},Z)$ is well defined and so
\begin{equation}\label{s49}
\lim_{\omega\to \frac{Z^2}{4}}T_0(\omega,Z)=\frac{4}{|Z|}\beta\Big(\frac{Z^2}{4},Z\Big).
\end{equation}

From Figure 6, for $Z$ fixed,  $\omega\to T_0(\omega, Z)$  is a strictly decreasing function. \eqref{s48} is a key result for our future analysis. In fact, for $Z>0$ fixed and $L>0$ there  exists $\omega>\frac{Z^2}{4}$ such that $2L>T_0(\omega, Z)$. Consequently, from Theorem \ref{funda1} there is a unique $\eta_2=\eta_2(\omega)\in (0,\theta(\omega, Z))$ such that
\begin{equation}\label{s50}
2s=T_-(\eta_2)=2L.
\end{equation}
In particular, for $Z =1$ and $L=1/2$ there is a unique $\omega_0>\frac14$, $\omega_0\approx 5.2$, such  $T_0(\omega_0,1)=1$ and for all
$\omega>\omega_0$ we have $1>T_0(\omega, Z)$. For $L$ large $\omega_0\to \frac14^+$, and for $L$ small $\omega_0$ is large. Also   for $\omega$ fixed,
$$
T_0(\omega, Z)\to \frac{\sqrt{2}}{\sqrt{\omega}}\pi\qquad \text{as}\; Z\to 0^+,
$$
and $Z\to T_0(\omega, Z)$ is a strictly increasing function. Then $\omega$ must satisfy $\omega>\frac{\pi^2}{2L^2}$ (see the theory in Angulo \cite{angulo1} for the case $Z=0$ in \eqref{p2}). 
\begin{figure}[!htb]
\centering
\includegraphics[angle=0,scale=0.4]{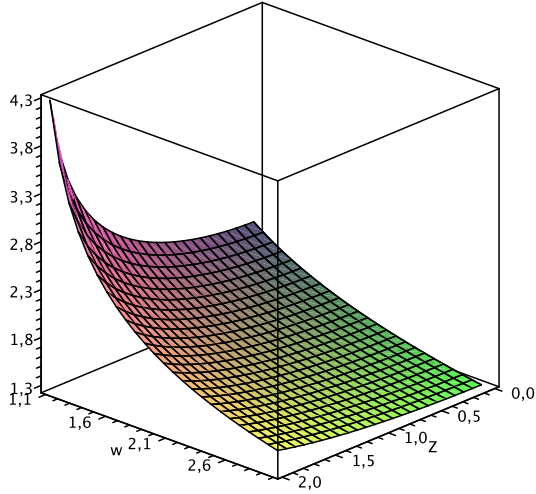}
\caption{Profile of the function $(\omega, Z)\to T_0(\omega, Z)$  in \eqref{s47}}
\end{figure}

From the above analysis, if we define  $\eta_1\equiv\sqrt{2\omega-\eta_2^2}$ for $\eta_2$ satisfying \eqref{s50}, $k^2$ and $a$ via the relations in \eqref{s39} and \eqref{s41}, respectively, we obtain the  peak-function
\begin{equation}\label{s51}
\phi(\xi)=\eta_1\text{dn}\Big(\frac{\eta_1}{\sqrt{2}} |\xi| +a;k\Big),
\end{equation}
for $\xi\in (-L,L)$. It satisfies item  (4) and item  (3) (for $\xi\in(-L,L)-\{0\}$) in \eqref{p3}. Moreover, since $\phi(\pm L)=\eta_2$ we can extend   $\phi$ to all the line as a  even periodic function with a minimal period $2L$ and in the interval $[0, 2L]$ it is symmetric with regard to the line $\xi=L$.  Hence,  we have obtained a periodic dnoidal-peak solution for equation \eqref{p2} which satisfies all  the properties in \eqref{p3} and it belongs to the domain of $-\frac{d^2}{dx^2}-Z\delta$.

\subsection{Smooth curve of periodic peaks to the NLS-$\delta$ with $Z\neq 0$}

In this section we construct a smooth curve of positive periodic peak solutions  of \eqref{p2} with $Z$ fixed. These solutions  $\varphi=\varphi_{\omega, Z}$  have   {\it a priori} fundamental period $2L$,  satisfy  the conditions in \eqref{p3}, and  $ \varphi_{\omega, Z}\in D(-\frac{d^2}{dx}-Z \delta)$.  Moreover, for $\omega>\frac{Z^2}{4}$ and $\omega$ fixed one has that
\begin{equation}\label{Z0}
\lim_{Z\to 0} \varphi_{\omega,Z}(x)=\phi_{\omega,0}(x) \qquad\text{for}\;\;x\in [-L,L],
\end{equation}
where $\phi_{\omega,0}$ is the dnoidal traveling wave defined in \eqref{dnoidal2} with a minimal period $2L$ and $\omega>\frac{\pi^2}{2L^2}$.  Our analysis  will show also that the mapping $Z\to \varphi_{\omega,Z}$ is analytic. This will be an essential  in our stability theory. Also we shall need to show  that  the map $\omega\to \eta_2(\omega)\in (0,\theta(\omega, Z))$ is smooth. 

First we consider the case $Z\neq 0$ and small. First, we shall establish  a result obtained by  Angulo in \cite{angulo1} for the cubic NLS. For $\eta\in (0,\sqrt{\omega})$, we define
\begin{equation}\label{s60}
F(\eta,\omega)=\frac{2\sqrt{2}}{\sqrt{2\omega-\eta^2}}K(k(\eta,\omega)).
\end{equation}
From the properties of  $K$ it follows that for $\omega>0$ fixed and $\eta\in (0,\sqrt{\omega})$ that the mapping $\eta\to F(\eta, \omega)$ is a strictly decreasing function and satisfies $F(\eta, \omega)>\sqrt{2}\pi/\sqrt{\omega}$. Hence, for $L>0$ fixed and $\omega_0>\frac{\pi^2}{2L^2}$ there exists a unique $\eta_0$ such that $F(\eta_0, \omega_0)=2L$. The following Theorem has been obtained in \cite{angulo1}.

\begin{teo}\label{conver0}
Let $L>0$  fixed. Consider $\omega_0>\frac{\pi^2}{2L^2}$ and $\eta_{0}=\eta({\omega_0})\in (0,\sqrt{\omega_0})$ such that $F(\eta_{0}, \omega_0)=2L$. Then there are  intervals $J_0(\omega_0)$ around $\omega_0$ and  $N_0(\eta_{0})$ around $\eta_{0}$, and a unique smooth function $\Lambda_0:J_0(\omega_0)\to N_0(\eta_{0})$ such that $\Lambda_0(\omega_0)=\eta_{0}$ and for $\eta\equiv \Lambda_0(\omega)$ one has $F(\eta, \omega)=2L$. Moreover,
$$
N_0(\eta_{0})\times J_0(\omega_0)\subseteq \{(\eta, \omega): \omega>\frac{\pi^2}{2L^2},\; \eta\in (0,\sqrt{\omega})\}.
$$
Furthermore, $J_0(\omega_0)$ can be taken equal to $(\frac{\pi^2}{2L^2},+\infty)$. For  $\eta_1=\eta_1(\omega)=\sqrt{2\omega-\eta^2}$, the dnoidal wave solution $\phi_{\omega, 0}$ defined in \eqref{dnoidal2}  has fundamental period $2L$ and satisfies the equation
$$
-\phi''_{\omega, 0}(x)+\omega \phi_{\omega, 0}(x)-\phi^3_{\omega, 0}(x)=0\quad\text{for all}\;\;x\in\mathbb R.
$$ 
Also,  $\omega\in (\frac{\pi^2}{2L^2},+\infty)\to \phi_{\omega, 0}\in H^n_{per}([0,2L])$ is a smooth function for all $n\in \mathbb N$.
\end{teo}


\subsubsection{Smooth curve of periodic peaks to the NLS-$\delta$ with $Z>0$}
 
We shall show that for $Z>0$ fixed,there exists a smooth curve $\omega \to \phi_{\omega, Z}\in H^1_{per}([-L,L])$ . Moreover, the convergence in \eqref{Z0} can be justified at least for $Z\to 0^+$.  The proof will be a consequence of the implicit function theorem and Theorem \ref{conver0}. We recall that $\omega>Z^2/4$.

\begin{teo}\label{conver1}
Let $L>0$  fixed,  $\delta$ small, $\delta<\frac{\pi^2}{2L^2}$, and  $Z\in (-\delta, \delta)$. Let  $\omega_0> \frac{\pi^2}{2L^2}$ and $\eta_0$ be the unique $\eta_0\in (0, \sqrt{\omega_0})$ such that $F(\eta_0, \omega_0)=2L$. Then,
\begin{enumerate}
\item there are an rectangle $R(\omega_0, 0)=J(\omega_0)\times (-\delta_0, \delta_0)$ around $(\omega_0, 0)$,  an interval $N_1(\eta_0)$ around $\eta_0$,  and a unique smooth function $\Lambda_1:R(\omega_0, 0)\to N_1(\eta_{0})$ such that $\Lambda_1(\omega_0, 0)=\eta_{0}$ and
 \begin{equation}\label{funda}
\frac{2\sqrt{2}}{\eta_{1,Z}}[K(k)-a(\omega, Z))]=2L,
\end{equation}
where $\eta^2_{1,Z}\equiv 2\omega-\eta_{2, Z}^2$ for $(\omega, Z)\in R(\omega_0, 0)$ and  $\eta_{2,Z}=\Lambda_1(\omega, Z)$.

\item $ J(\omega_0)=(\frac{\pi^2}{2L^2}, +\infty)$ and $k\in (k_0,1)$, $k_0$ defined in \eqref{s44}.
\item $N_1(\eta_0)\times R(\omega_0, 0)\subset \mathbb G=\{(\eta, \omega, Z): \omega>\frac{\pi^2}{2L^2}, 2L>T_0(\omega,Z), \eta\in (0,\theta(\omega, Z))\}$.

\item For $Z=0$ we have $a(\omega, 0)=0$ and so from Theorem \ref{conver0}  it follows that $\Lambda_1(\omega, 0)=\Lambda_0(\omega)$. Therefore, $\lim_{Z\to 0^+}\eta_{2,Z}(\omega)=\eta(\omega)$.

\item For $Z\in (0, \delta_0)$ we denote $\eta_{2,Z}$ by $\eta_{2,+}$. Then the dnoidal-peak solution $\phi_{\omega, Z}$ in \eqref{s28} with $\eta_1$ being  $\eta_{1, +}=(2\omega^2-\eta_{2, +}^2)^{1/2}$, has minimal period $2L$ and satisfies for $\omega>\frac{\pi^2}{2L^2}$
\begin{equation}\label{conver3}
\lim_{Z\to 0^+} \phi_{\omega, Z}(x)=\phi_{\omega,0}(x),\quad \text{for}\; x\in [-L,L]. 
\end{equation}

\item $Z\to \phi_{\omega, Z}\in H^1_{\text{per}}([-L,L]) $ is real-analytic.

\end{enumerate}
\end{teo}

\begin{proof} The proof is a consequence of the implicit function theorem applied to the mapping
$$
G(\eta, \omega, Z)= \frac{2\sqrt{2}}{\sqrt{2\omega-\eta^2}}[K(k(\eta, \omega))-a(\eta, \omega, Z))]
$$
with domain $\mathbb G$. From \eqref{s48}  follows  $\mathbb G \neq\emptyset$. Moreover, if $(\eta_0, \omega_0, Z)\in \mathbb G$ then for all $ \omega>\omega_0$ we obtain $(\eta_0, \omega, Z)\in \mathbb G$ ( $G(\eta_0, \omega_0, 0)=2L$). Next, we claim  that $\partial_\eta G(\eta_0, \omega_0, 0)<0$. Indeed, from  Theorem 2.1 in Angulo \cite{angulo1} we have $\partial_\eta F(\eta,\omega)<0$ since $\partial_\eta k(\eta, \omega)$ is a strictly decreasing function of $\eta$, since $\partial_\eta a(\eta,\omega, Z)>0$ (see Theorem \ref{funda1} or Figure  5 above with a $Z$ fixed) we prove the claim. Theorem \ref{conver0} implies item (2) above.

Finally, since the functions $a$ in \eqref{s41}  and $\Lambda_1$  are analytic,  the mapping $(\omega, Z)\to \phi_{\omega, Z}$ is analytic for $(\omega, Z)\in \{(\omega, Z): \omega>Z^2/4\}$, $Z$ small. This finishes the proof of the Theorem. 
\end{proof}

\begin{coro}\label{s61} Consider the mapping $\Lambda_1:R(\omega_0, 0)\to N_(\eta_{0})$ obtained in Theorem \ref{conver1}. Then for $Z$ fixed, the mapping $\omega\to \eta_{2,+}(\omega)=\Lambda_1(\omega, Z)$  is a strictly decreasing function. Moreover,  for $k$ and $a$ defined in \eqref{s39} and \eqref{s41} respectively, one has  $\frac{d}{d\omega} k(\omega)>0$ and $\frac{d}{d\omega} a(\omega)<0$.
\end{coro}

\begin{proof} Let $Z$ fixed. Since $G(\Lambda_1(\omega, Z),\omega, Z)=2L$ one has that $\eta_{2,+}'(\omega)=-\frac{\partial G/{\partial \omega}}{\partial G/{\partial \eta}}$. 
Using that $\partial_\omega G(\eta,\omega, Z)<0$ (see Figure 11) we obtain $\eta_{2,+}'(\omega)<0$. Next, for $a(\omega)\equiv a(\Lambda_1(\omega, Z),\omega, Z)$ therefore
$$
\frac{d}{d\omega}a(\omega)=\frac{\partial a}{\partial \Lambda_1}\frac{d\Lambda_1}{d \omega}+ \frac{\partial a}{\partial \omega}<0 
$$
since $\frac{\partial a}{\partial \Lambda_1}>0$ and $ \frac{\partial a}{\partial \omega}<0$ (see Figure 5). Finally, from the formula of $k$ in \eqref{s39} it follows  immediately that $k(\omega)\equiv k(\Lambda_1(\omega, Z),\omega, Z)$ is a strictly increasing function. This completes the proof of the Corollary. 
\end{proof}

By using Maple's software we can give a general profile of $G(\eta,\omega, Z)$, $Z$ fixed. For instance, for $L=1/2$ and $Z =1$  the analysis in subsection 5.4 tell us that  $1>T_0(\omega,1)$ for all $\omega>\omega_0\approx 5.2$ with $T_0(\omega_0,1)=1$. So, we obtain the profile of $G(\eta,\omega, 1)$ for $\omega >\omega_0$ and $\eta\in (0, \theta(\omega,1))$ given by   Figure 7.  We observe  that $\partial_\omega G(\eta,\omega)<0$ and $\partial_\eta G(\eta,\omega)<0$.

\begin{figure}[!htb]
\centering
\includegraphics[angle=0,scale=0.3]{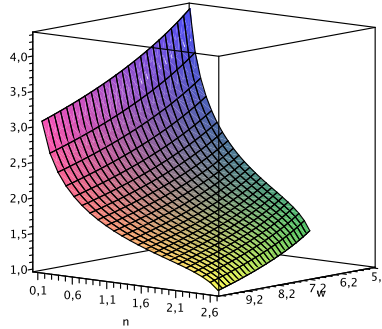}
\caption{Graph of  $G (\eta,\omega, 1)$}
\end{figure}

In the next section,  we will need to use that the mapping $Z\to \phi_{\omega, Z}$ is analytic for $Z> 0$ (we recall that this property  is   local type). So, by using an argument similar  to that provided  in the proof of Theorem \ref{conver1}  and the analysis in subsection 5.4 we obtain :

\begin{teo}\label{s51b}
Let $L>0$ fixed and $Z_0>0$. Consider $\omega_0>\frac{Z_0^2}{4}$ such that
$2L>T_0(\omega_0, Z_0)$ and $\omega_0> \frac{\pi^2}{2L^2}$. Let $\eta_{2,0}=\eta_{2,0}(\omega_0, Z_0)\in (0,\theta(\omega_0, Z_0))$ the unique value such that $T_-(\eta_{2,0}, \omega_0, Z_0)=2L$. Then,
\begin{enumerate}
\item there are an rectangle $S(\omega_0, Z_0)=H(\omega_0)\times I(Z_0)$ around $(\omega_0, Z_0)$,  an interval $N_2(\eta_{2,0})$ around $\eta_{2,0}$,  and a unique smooth function $\Lambda_2:S(\omega_0, Z_0)\to N_1(\eta_{2,0})$ such that $\Lambda_2(\omega_0, Z_0)=\eta_{2, 0}$ and $T_-(\eta_{2,+}, \omega, Z)=2L$ for $\eta_{2,+}=\Lambda_2(\omega, Z)$.

\item $H(\omega_0)$ can be choosen as $(\mu(Z, L),+\infty)$, where $\mu(Z, L)>\frac{Z_0^2}{4}$ and $\mu(Z, L) >\frac{\pi^2}{2L^2}$. For $Z=0$ we have $\mu(0, L)=\frac{\pi^2}{2L^2}$.

\item the dnoidal-peak solution in \eqref{s28}, $\phi_{\omega, Z}(\xi)\equiv \phi(\xi;\eta_{1, +})$, determined by $\eta_{1,+}\equiv (2\omega-\eta_{2,+}^2)^{1/2}$  satisfies the properties in \eqref{p3}. Moreover, the mapping
\begin{equation}\label{s56a}
Z\to \phi_{\omega, Z}\in H^1_{\text{per}}([-L,L])
\end{equation}
is real-analytic.
\end{enumerate}
\end{teo}

\begin{coro}\label{s61a} For $Z$ fixed, the mapping $\omega\to \eta_{2, +}(\omega)=\Lambda_2(\omega, Z)$ is a strictly decreasing function. Moreover,  for $k$ and $a$ defined in \eqref{s39} and \eqref{s41} respectively, one has that $\frac{d}{d\omega} k(\omega)>0$ and $\frac{d}{d\omega} a(\omega)<0$
\end{coro}

\begin{coro}\label{s62a} For $Z\geqq 0$ fixed, consider the mapping $a:(\mu(Z,L), +\infty)\to \mathbb R$ obtained in  Theorem \ref{conver1} and Theorem \ref{s51b}. Then $a(\omega)\to 0$ as $\omega\to +\infty$.
\end{coro}

\begin{proof} We consider $a$ given by Theorem \ref{s51b}. From Corollary \ref{s61a} it follows  that for $\omega>\omega_1$ and $\eta_{2, +}(\omega)=\Lambda_2(\omega, Z)$
$$
0\leqq \frac{\eta_{2, +}^2(\omega)}{\omega}\leqq \frac{\eta_{2, +}^2(\omega_1)}{\omega}.
$$
Thus,  $k^2(\omega)\to 1^+$ and $\frac{\Phi(\omega)}{\sqrt{2\omega-\eta_{2,+}^2}}\to 1^+$ as $\omega\to +\infty$. Therefore,  \eqref{s41} yields the identity
\begin{equation}\label{s63}
\lim_{\omega\to+\infty} a(\omega)=sech^{-1}(1)=0.
\end{equation}
\end{proof}

\subsubsection{Smooth curve of periodic peaks to the NLS-$\delta$ with $Z<0$}

The following Theorem shows that for $Z<0$, fixed there exists a smooth curve $\omega \to \zeta_{\omega, Z}\in H^1_{per}([-L,L])$ and  that the convergence in \eqref{Z0} for $Z\to 0^-$ is possible.  The proof is  similar to that of Theorem \ref{conver1} and Theorem \ref{s51b}, so we  shall only describe the main points in the argument. We start  by defining
\begin{equation}\label{negaZ0}
T_+(\eta_2, \omega)=\frac{2\sqrt{2}}{\sqrt{2\omega-\eta_2^2}}[K(k(\eta_2, \omega))+a(\eta_2, \omega)]
\end{equation}
and
\begin{equation}\label{negaZ1}
T_1(\omega, Z)=\frac{2\sqrt{2}}{\lambda(\omega,Z)}[K(k_0)+a_0]
\end{equation}
where $T_1(\omega, Z)=\lim_{\eta_{2}\to \theta} T_+(\eta_2, \omega)$. From  \eqref{s48} and $\lim_{\omega\to +\infty} a_0(\omega)=0$ it follows that  
$
\lim_{\omega\to +\infty}T_1(\omega,Z)=0$. So, since the mapping $\omega\to T_1(\omega,Z)$ is a strictly decreasing function we obtain a unique $\omega_1>\frac{Z^2}{4}$ such that $2L>T_1(\omega_1,Z)$ and  for every $\omega>\omega_1$, $2L>T_1(\omega,Z)$. Now, for $\omega$ chosen in this form one finds a unique $\eta_{2,1}=\eta_{2,1}(\omega)\in (0, \theta(\omega, Z))$ such $T_+(\eta_{2,1}, \omega)=2L$.  Moreover,  since $T_1(\omega, Z)=T_0(\omega, Z)+\frac{4\sqrt{2}}{\lambda(\omega,Z)}a_0\to \frac{\sqrt{2}}{\sqrt{\omega}}\pi$ as $Z\to 0^-$, we obtain {\it a priori} the condition $\omega>\frac{\pi^2}{2L^2}$. 

We have the following theorem of existence.

\begin{teo}\label{negaZ3}
Let $L>0$  fixed and $Z_0<0$. Consider $\omega_1>\frac{Z_0^2}{4}$
such that $2L>T_1(\omega_1, Z_0)$ and $\omega_1>\frac{\pi^2}{2L^2}$. Let $\eta_{2,1}=\eta_{2,1}(\omega_1, Z_0)\in (0,\theta(\omega_1, Z_0))$ the unique value such that $T_+(\eta_{2,1})=2L$. Then,
\begin{enumerate}
\item there are an rectangle $W(\omega_1, Z_0)=Q(\omega_1)\times V(Z_0)$ around $(\omega_1, Z_0)$,  an interval $N_2(\eta_{2,1})$ around $\eta_{2,1}$,  and a unique smooth function $\Lambda_3:W(\omega_1, Z_0)\to N_2(\eta_{2,1})$ such that $\Lambda_3(\omega_1, Z_0)=\eta_{2, 1}$ and $T_+(\eta_{2,-}, \omega, Z)=2L$ for $\eta_{2,-}=\Lambda_3(\omega, Z)$,

\item $Q(\omega_1)$ can be choosen as $(\nu(Z, L),+\infty)$, where $\nu(Z, L)>\frac{Z_0^2}{4}$ and $\nu(Z, L) >\frac{\pi^2}{2L^2}$. For $Z=0$ we have $\nu(0, L)=\frac{\pi^2}{2L^2}$,

\item for $Z=0$ we have  $a(\omega, 0)=0$ and so from Theorem \ref{conver0} we have $\Lambda_3(\omega, 0)=\Lambda_0(\omega)$. Therefore, $\lim _{Z\to 0^-}\eta_{2,-}(\omega)=\eta(\omega)$,

\item the dnoidal-peak solution in \eqref{sne14}, $\zeta_{\omega, Z}(\xi)\equiv \zeta(\xi;\eta_{1, -})$, determined by $\eta_{1,-}\equiv (2\omega-\eta_{2,-}^2)^{1/2}$  satisfies the properties in \eqref{p3}. Moreover, the mapping
\begin{equation}\label{s56}
Z\to \zeta_{\omega, Z}\in H^1_{\text{per}}([-L,L])
\end{equation}
is real-analytic,

\item from the condition $\omega>\frac{\pi^2}{2L^2}$ we obtain 
\begin{equation}\label{negaZ9}
\lim_{Z\to 0^-} \zeta_{\omega, Z}(\xi)=\phi_{\omega,0}(\xi), \quad\text{for}\; \xi\in [-L,L].
\end{equation}
\end{enumerate}
\end{teo}

\begin{coro}\label{negaZ10} For $Z<0$ fixed, the mapping $\omega\to \eta_{2, -}(\omega)=\Lambda_2(\omega, Z)$ is a strictly decreasing function. Moreover,  for $k$ and $a$ defined in \eqref{s39} and \eqref{s41} respectively, one has  that  $\frac{d}{d\omega} k(\omega)>0$ and $\frac{d}{d\omega} a(\omega)<0$
\end{coro}

\begin{proof} For $T_+$ defined in \eqref{negaZ0}, it follows  that $\partial_\eta T_+<0$  and $\partial_\omega T_+<0$. Then for $\Lambda_3$ satisfying $T_+(\Lambda_3(\omega, Z), \omega)=2L$  and $a(\omega)=a(\Lambda_3(\omega, Z), \omega, Z)$  we obtain that $\Lambda_3'(\omega)<0$ and $a'(\omega)<0$.
\end{proof}

\begin{coro}\label{s62b} For $Z\leqq 0$ fixed, consider the mapping $a:(\nu(Z,L), +\infty)\to \mathbb R$ determined by Theorem \ref{negaZ3}. Then $a(\omega)\to 0$ as $\omega\to +\infty$.
\end{coro}

\section{Stability of Dnoidal-Peak  for NLS-$\delta$}

In this section we study the stability  of the orbit 
\begin{equation}\label{orb1}
\Omega_{\varphi_{\omega,Z}}=\{e^{i\theta}\varphi_{\omega,Z}: \theta\in [0,2\pi)\},
\end{equation}
generated by the smooth curve of dnoidal-peak $\omega\to \varphi_{\omega,Z}$, where
\begin{equation}\label{pp}
\varphi_{\omega,Z}=
\begin{cases}
\begin{aligned}
&\phi_{\omega, Z},\quad  Z>0,\\
&\zeta_{\omega, Z},\quad  Z<0
\end{aligned}
\end{cases}
\end{equation}
with $\phi_{\omega, Z}$  and $\zeta_{\omega, Z}$  are given by Theorems \ref{conver1},  \ref{s51b} and \ref{negaZ3}. Moreover, 
\begin{equation}\label{pp1}
\lim_{Z\to 0} \varphi_{\omega, Z}(\xi)=\phi_{\omega,0}(\xi), \quad\text{for}\; \xi\in [-L,L],
\end{equation}
where $\phi_{\omega,0}$ being the dnoidal-wave solution to the cubic Schr\"odinger equation determined by Theorem \ref{conver0}.

We start obtaining  the spectral information associated to the operators in \eqref{specA6} and \eqref{specA7} necessary  to establish our  stability theorem.

\subsection{The basic linear operators $\mathcal L_{1, Z}$ and $\mathcal L_{2, Z}$}  For $u\in H^1_{per}([0,1])$ we write  $u=u_1+iu_2$. Let $H_{\omega, Z}$ be defined by
\begin{equation}\label{spec0}
H_{\omega, Z} u= \mathcal L_{1,Z}u_1+i \mathcal L_{2, Z}u_2
\end{equation}
where the linear operators $\mathcal L_{i, Z}$, $i=1,2$, are defined as: 
\begin{equation}\label{spec1}
\begin{aligned}
\mathcal D\equiv D(\mathcal L_{i, Z})=\{\zeta\in H^1_{\text{per}}([-L, L])&\cap H^2((-L,L)-\{0\})\cap H^2((2nL,2(n+1)L)):\\
& \zeta'(0+)-\zeta'(0-)=-Z \zeta(0)\},
\end{aligned}
\end{equation}
and for $\zeta\in \mathcal D$
\begin{equation}\label{spec2}
\begin{aligned}
&\mathcal L_{1, Z}\zeta =-\frac{d^2}{dx^2}\zeta+\omega \zeta-3\varphi^2_{\omega, Z}\zeta,\\
&\mathcal L_{2, Z}\zeta =-\frac{d^2}{dx^2}\zeta+\omega \zeta-\varphi^2_{\omega, Z}\zeta.
\end{aligned}
\end{equation}

We claim that $\mathcal L_{i, Z}$ are self-adjoint operators on $L^2_{per}([-L,L])$ with domain $\mathcal D$. Since the multiplication operator $
\mathcal M\zeta =(\omega-3\varphi^2_{\omega,Z})\zeta$
is obviously symmetric and bounded on $L^2_{per}([-L,L])$ and  $D(-\Delta_{-Z})\subset D(\mathcal M)=L^2_{per}([-L,L])$, it follows from the Stability Self-Adjoint Theorem (see Kato \cite{kato}) that $\mathcal L_{1,Z}=-\frac{d^2}{dx^2}\zeta+ \mathcal M\zeta$ is a self-adjoint operator on $L^2_{per}([-L,L])$ with domain $
D(\mathcal L_{1, Z})=D(-\Delta_{-Z})=\mathcal D$. A similar result holds for $\mathcal L_{2, Z}$.  

We note that the linear operators $\mathcal L_{1, Z}$ and $\mathcal L_{2, Z}$ are related with the the second variation of $G_{\omega,Z}=E+\omega Q$ at $\varphi_{\omega,Z}$. More exactly, let $u=\zeta+i\psi$ with $\zeta, \psi\in \mathcal D$ and $v=v_1+iv_2\in H^1_{per}$ then
\begin{equation}\label{hessian}
\langle G''_{\omega, Z}(\varphi_{\omega,Z})u,v\rangle= \langle H_{\omega, Z}u,v\rangle=
 \langle \mathcal L_{1, Z}\zeta+i\mathcal L_{2, Z}\psi,v\rangle=\langle\mathcal L_{1, Z} \zeta, v_1\rangle + \langle\mathcal L_{2, Z} \psi, v_2\rangle.
 \end{equation}
Next we give a idea of the proof of the equality in \eqref{hessian}. For $\zeta$  and $v_1$ we define 
$$
\mathcal Q(\zeta, v_1)=\omega\int \zeta v_1dx-3\int\varphi_{\omega, Z}^2\zeta v_1dx.
$$
Thus,
\begin{equation}\label{spec5}
\begin{aligned}
\langle\mathcal L_{1, Z} \zeta, v_1\rangle  
&=-\lim_{\epsilon \downarrow 0}\int_{-L}^{-\epsilon}\big(\frac{d^2}{dx^2}\zeta\big) v_1dx-\lim_{\epsilon \downarrow 0}\int_{\epsilon}^{L}\big(\frac{d^2}{dx^2}\zeta\big) v_1dx+\mathcal Q(\zeta, v_1)\\
&=\lim_{\epsilon \downarrow 0}[\zeta'(\epsilon) v_1(\epsilon)-\zeta'(-\epsilon) v_1(-\epsilon)]+ \int \zeta'\;v'_1dx+\mathcal Q(\zeta, v_1)\\
&=[\zeta'(0+)-\zeta'(0-)] v_1(0)+ \int \zeta'\; v'_1dx+\mathcal Q(\zeta, v_1)\\
&=-Z \zeta(0) v_1(0)+\int \zeta'\; v_1'dx+\mathcal Q(\zeta, v_1).
\end{aligned}
\end{equation}
Similarly, we obtain
$$
\langle\mathcal L_{2, Z} \psi, v_2\rangle=-Z \psi(0) v_2(0)+\int \psi'\; v_2'dx+\omega\int \psi v_2dx-\int\varphi_{\omega, Z}^2\psi v_2dx.
$$
A simple calculation shows that
$$
\langle G''_{\omega,Z}(\varphi_{\omega, Z})(\zeta, \psi),(v_1,v_2)\rangle=\langle\mathcal L_{1, Z} \zeta, v_1\rangle + \langle\mathcal L_{2, Z} \psi, v_2\rangle.
$$

\subsection{Some spectral structure of $\mathcal L_{1, Z}$ and $\mathcal L_{2, Z}$}

This subsection is concerned with some specific spectral structure of the linear operators $\mathcal L_{i, Z}$. By convenience we will denote $\mathcal L_{i, Z}$ only by $\mathcal L_i$.

\begin{lema}\label{spec6} Let $Z\in \mathbb R$ and $\omega>Z^2/4$. Then,

\begin{enumerate}

\item $\mathcal L_2$ is  a nonnegative operator with a discrete spectrum, 
$\sigma(\mathcal L_2)=\{\lambda_n\,:\,n\geqq 0\}$, ordered in the increasing form
\begin{equation}\label{specL2}
0=\lambda_0<\lambda_1\leqq \lambda_2\leqq \lambda_3\leqq \lambda_4\cdot\cdot\cdot.
\end{equation}
The eigenvalue zero is simple  with eigenfunction $\varphi_{\omega, Z}$.

\item $\mathcal L_1$ is a operator with a discrete spectrum, 
$\sigma(\mathcal L_1)=\{\alpha_n\,:\,n\geqq 0\}$, ordered in the increasing form
\begin{equation}\label{specL1}
\alpha_0<\alpha_1\leqq \alpha_2\leqq \alpha_3\leqq \alpha_4\cdot\cdot\cdot.
\end{equation}

\end{enumerate}
\end{lema}

\begin{proof} From Section 3, Theorem \ref{resol5b}, it follows  that the operators $\mathcal L_i$ have a compact resolvent and so its spectrum is discrete satisfying \eqref{specL2} and \eqref{specL1}. Since $\varphi_{\omega, Z}\in \mathcal D$ and satisfies equation \eqref{p2} we obtain that $\mathcal L_2\varphi_{\omega, Z} =0$ for all $Z$. Moreover, $\varphi_{\omega, Z}$ being positive it corresponds to the first eigenvalue of $\mathcal L_2$ which is simple. \end{proof}

Next we have  the following kernel-structure of  $\mathcal L_1$. 

\begin{lema}\label{spec6a1} Let $Z\in \mathbb R-\{0\}$ and $\omega>Z^2/4$. Then $\mathcal L_1$ has a trivial kernel.
\end{lema}

\begin{proof}   Let $v\in \mathcal D$ such that $\mathcal L_1 v=0$ and $Z>0$. Therefore $\varphi_{\omega, Z}=\phi_{\omega, Z}$. We claim  that the subspace $v$ of $L_{per}^2([0, 2L])$-solutions  of the problem
\begin{equation}\label{spec7}
\begin{cases}
\begin{aligned}
&v\in H^2(0,2L)\\
&\mathcal L_1 v(x)=0\quad{\rm{for}}\;\;x\in (0,2L),
\end{aligned}
\end{cases}
\end{equation}
is one dimensional. From  item (3) in \eqref{p3} it follows that that $\Lambda_1(x)\equiv \phi'_{\omega,Z}(x)$ for $x\in (0,2L)$ satisfies  problem \eqref{spec7}.  We consider the transformation
$$
\Lambda(x)=v(\beta x),\quad{\rm{for}}\;\;\beta=\frac{\sqrt{2}}{\eta_1},\;\;x\in (0,2(K-a))
$$
with $a$ defined in \eqref{s41}. Then from Theorems \ref{conver1} and  \ref{s51b} we have $\beta x\in (0,2L)$ and so  \eqref{spec7} implies
that \begin{equation}\label{spec8}
\Lambda''(x)+[\sigma-6k^2 sn^2(x+a;k)]\Lambda (x)=0\quad{\rm{for}}\;\;x\in (0, 2(K-a)),
\end{equation}
where $ \sigma=(6\eta_{1,+}^2-2\omega)/\eta_{1,+}^2=4+k^2$. Now, for $\Upsilon(x)=\Lambda(x-a)$ with $x\in (a,2K-a)$ we have that $\Upsilon$ satisfies the following Lam\'e's equation
\begin{equation}\label{spec8a}
\Phi''(x)+[\sigma-6k^2 sn^2(x;k)]\Phi(x)=0, \quad{\rm{for}}\;\;x\in (a, 2K-a).
\end{equation}

Next, from Angulo \cite{angulo1} (see Lemma \ref{specL0} below) the periodic eigenvalue problem in $L^2_{per}([0,2K])$
\begin{equation}\label{spec9}
\begin{cases}
\begin{aligned}
&\Phi''(x)+[\lambda-6k^2 sn^2(x;k)]\Phi (x)=0,\quad    x\in (0, 2K)\\
&\Phi(0)=\Phi(2K(k)), \;\; \Phi'(0)=\Phi'(2K(k)),\;\;k\in (0,1)
\end{aligned}
\end{cases}
\end{equation}
has the first three eigenvalues $\lambda_0, \lambda_1, \lambda_2$  simple and the rest of the eigenvalues are distributed in the form $\lambda_3\leqq \lambda_4<\lambda_5\leqq \lambda_6<\cdot\cdot\cdot $ and  satisfying $\lambda_3=\lambda_4, \lambda_5=\lambda_6,...$, i.e., they are double eigenvalues and so for these values of $\lambda$ all solutions of \eqref{spec9} have period $2K(k)$. In particular,  $\lambda_1=4+k^2$ and  $\Phi_1(x)=sn(x;k)cn(x;k)=C_0\frac{d}{dx} dn(x;k)$,  
for $x\in [0, 2K(k)]$, $k\in (0,1)$.  Now, from Floquet theory (see pg. 7 in \cite{Ea}) the other solution for \eqref{spec8a} with $\rho= \lambda_1$ is of the form $\Psi(x)=x\Phi_1(x)+p_2(x),$ where $\Psi$ is even and  $p_2(x)$ has period $2K(k)$. We recall that  $\{\Phi_1, \Psi\}$ is a linearly independent (LI) set of solutions for  \eqref{spec8a} on $\mathbb R$ and so it is a base of solutions for \eqref{spec8a} on any interval $(c,d)$. Then the following functions  on $(0,2L)$, for $\eta_1=\eta_{1, +}$,
\begin{equation}\label{spec10a}
\begin{cases}
\begin{aligned}
&\Lambda_1(x)= \frac{d}{dx}\phi_{\omega,Z}(x),\quad \text{and}\\
&\Lambda_2(x)= \Big(\frac{\eta_{1}}{\sqrt{2}} x+a\Big)\Lambda_1(x)+p_2\Big(\frac{\eta_{1}}{\sqrt{2}} x+a\Big)
\end{aligned}
\end{cases}
\end{equation}
are a LI set of solutions for \eqref{spec7} on $(0,2L)$. Therefore, there are $\alpha, \beta\in \mathbb R$ such that
\begin{equation}\label{spec10b}
v(x)=\alpha \Lambda_1(x)+\beta \Lambda_2(x),  \quad x\in (0,2L).
\end{equation}
Suppose $\beta\neq 0$. Then,  since $v$ and $\Lambda_1$ are periodic of period $2L$, we have that $\Lambda_2$ is a periodic function with period $2L$, which is not possible.  Therefore $v(x)=\alpha \Lambda_1(x)$ on $(0,2L)$. By establishing a similar problem to \eqref{spec7} on $(-2L,0)$, we can show that for $x\in (-2L,0)$, $v(x)=\alpha \Lambda_1(-x)$,which proves the claim. Moreover, $v$ is even and so   $v'(0+)=-v'(0-)$. Then,
\begin{equation}\label{spec10c}
v'(0+)=-\frac{Z}{2}v(0).
\end{equation}

Next  we will prove that  $\Lambda_1$ does not satisfy  condition \eqref{spec10c}. Indeed, we know that $\phi'_{\omega, Z}(0+)=-\frac{Z}{2}\phi_{\omega, Z}(0)$ and from \eqref{p3} 
$$
\phi''_{\omega, Z}(0+)=\phi''_{\omega, Z}(0)=\lim_{x\to 0^+} \phi''_{\omega, Z}(x)=\omega \phi_{\omega, Z}(0)-\phi^3_{\omega, Z}(0).
$$
Suppose now that $\phi''_{\omega, Z}(0+)=-\frac{Z}{2}\phi'_{\omega, Z}(0+)$. Then it follows that 
$$
\frac{Z^2}{4}\phi_{\omega, Z}(0)=\omega \phi_{\omega,Z}(0)-\phi^3_{\omega,Z}(0),
$$
which together with  \eqref{s5} and \eqref{s7} implies that for $\eta_2=\eta_{2,+}$, $
2\omega-\frac{Z^2}{2}=2\eta_1\eta_2$.  Hence, we can conclude from $2\omega=\eta^2_1+\eta^2_2$, that
$$
 \frac{Z^2}{2}=2\omega-2\eta_1\eta_2=(\eta_1-\eta_2)^2,
$$
which is a contradiction with the inequality in \eqref{s8}. Therefore $\alpha =0$ and so for $Z>0$ one has that  $Ker(\mathcal L_1)=\{0\}$. The case $Z<0$ follows similarly.  This finishes the proof.
\end{proof}

The next result will be used more later, but it is also interesting by itself. 

\begin{lema} \label{pone} Let $Z\in \mathbb R-\{0\}$ and $\omega>Z^2/4$. If $\lambda$ is an simple eigenvalue for $\mathcal L_1$  then the eigenfunction associated is either even or odd.
\end{lema}

\begin{proof} let $v\in D(\mathcal L_1)-\{0\}$ such that $ \mathcal L_1v=\lambda v$. Then, since $\varphi_{\omega, Z}$ is even, we also have for $\zeta(x)\equiv v(-x)$ the relation  $\mathcal L_1\zeta(x)=\lambda \zeta(x)$. Then there exists $\beta\in \mathbb R$ such that  $v(x)=\beta v(-x)$ for $x\in \mathbb R$. If $v(0)\neq 0$  then $\beta =1$ and thus $v$ is even. If $v(0)=0$ from \eqref{spec1} we have that $v\in H^2_{per}$ (see Remark after Theorem \ref{resol4r}) and so $v'(x)$ exists for $x\in \mathbb R$. Then we get that $v'(0)=-\beta v'(0)$ and from the Cauchy uniqueness principle $v'(0)\neq 0$ (in other way, $v\equiv 0$). Therefore $\beta=-1$ and so $v$ is a odd function.
\end{proof}

{\bf{Remark:}}  For the case $Z>0$, the even function
$$
\Upsilon_1(x)=2x sn(x) cn(x) -\frac{1+k'^2}{k'^2}E(x) sn(x) cn(x)+  \frac{1}{k'^2}dn(x)((sn^2(x)-k'^2cn^2(x))
$$
with $E(x)=\int_0^x dn^2(z) dz$, satisfies equation \eqref{spec8a} for all $x\in \mathbb R$ (see Figure 8 below). We note that it profile has the property that for $x\in (a, 2K-a)$, $\Upsilon_1$ it is not symmetric with respect to $x=K$. It  property can be used for an alternative  proof of Theorem \ref{spec6a1}.  
\begin{figure}[!htb]
\centering
\includegraphics[angle=0,scale=0.2]{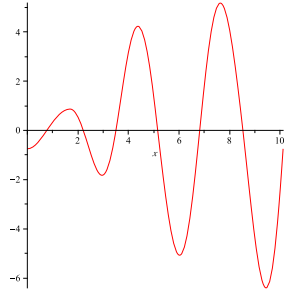}
\caption{Profile of $\Upsilon_1$ on $[0,6K(k)]$ with $k=0.5$.}
\end{figure}

\subsection{Counting the negative eigenvalues for $\mathcal L_{1, Z}$}

In this subsection we use  the theory of perturbation for linear operators to determinate   the number of negative eigenvalues of  $\mathcal L_{1, Z}$ for $Z\neq 0$.  Since the domain $\mathcal D$ of these operators is changing with $Z$ we will use the theory of analytic perturbation for linear operators (see \cite{kato} and \cite{RS}) and  some arguments found in \cite{lffgs}. Our study will be divided into four  steps:

\begin{enumerate}

\item[{(I)}] From our analysis in Section 5  it follows that by fixing $\omega>\frac{\pi^2}{2L^2}$ one has that 
\begin{equation}\label{spec13}
\lim_{Z\to 0}\varphi_{\omega, Z}=\phi_{\omega,0} \qquad{\rm{in}}\;\; H^1_{per}([-L,L])
\end{equation}
where $\phi_{\omega,0}$ represents  the dnoidal periodic solution in \eqref{dnoidal2}.

\item[{(II)}] The linear operators $\mathcal L_i$ in \eqref{spec2} are  the self-adjoint operators on $L^2_{per}([-L,L])$ associated with  the following bilinear forms defined for $v, w\in H^1_{per}([-L,L])$,
\begin{equation}\label{spec14}
\begin{aligned}
\mathcal Q_{\omega, Z}^1(v,w)&=\int_{-L}^Lv_x w_xdx +\omega \int_{-L}^Lvwdx  -Zv(0)w(0)-\int_{-L}^L3\varphi_{\omega, Z}^2vwdx\\
\mathcal Q_{\omega, Z}^2(v,w)&=\int_{-L}^Lv_x w_xdx +\omega \int_{-L}^Lvwdx  -Zv(0)w(0)-\int_{-L}^L\varphi_{\omega, Z}^2vwdx.
\end{aligned}
\end{equation}
Since these forms have the same domain $D(\mathcal  Q_{\omega, Z}^i)=H^1_{per}([-L,L])$ and  they are symmetric, bounded from below and closed, from the theory of representation of forms by operators (The First Representation Theorem in  \cite{kato}, VI. Section 2.1), one has that there are two self-adjoint  operators $\widetilde{\mathcal L_1}: D(\widetilde{\mathcal L_1})\subset L^2_{per}([-L,L])\to L^2_{per}([-L,L])$ and $\widetilde{\mathcal L_2}: D(\widetilde{\mathcal L_2})\subset L^2_{per}([-L,L])\to L^2_{per}([-L,L])$ such that
\begin{equation}\label{spec15}
\begin{aligned}
D(\widetilde{\mathcal L_1})&=\{v\in H^1_{per}: \exists w\in   L^2_{per}\; s.t.\; \forall z\in H^1_{per}, \;\mathcal Q_{\omega, Z}^1(v,z)=(w,z)\},\\
D(\widetilde{\mathcal L_2})&=\{v\in H^1_{per}: \exists w\in   L^2_{per}\; s.t.\; \forall z\in  H^1_{per},\; \mathcal Q_{\omega, Z}^2(v,z)=(w,z)\},
\end{aligned}
\end{equation}
and for $v\in D(\widetilde{\mathcal L_1})$ (resp. $v\in D(\widetilde{\mathcal L_2})$) we define $\widetilde{\mathcal L_1}v\equiv w$ (resp.  $\widetilde{\mathcal L_2}v\equiv w$), where $w$ is the (unique) function of $L^2_{per}([-L,L])$ which satisfies $\mathcal Q_{\omega, Z}^1(v,z)=(w,z)$ (resp. $\mathcal Q_{\omega, Z}^2(v,z)=(w,z)$) for all $z\in H^1_{per}$.

Next, we describe explicitly the self-adjoint operators $\widetilde{\mathcal L_1}$ and $\widetilde{\mathcal L_2}$.

\begin{lema}\label{repres} The domain for both  $\widetilde{\mathcal L_1}$ and $\widetilde{\mathcal L_2}$  in $L^2_{per}([-L,L])$ is
\begin{equation}
\begin{aligned}
D_Z=\{\zeta\in H^1_{\text{per}}([-L, L])&\cap H^2((-L,L)-\{0\})\cap H^2((2nL, 2(n+1)L)):\\
& \zeta'(0+)-\zeta'(0-)=-Z \zeta(0)\},
\end{aligned}
\end{equation}
and for $v\in D_Z$ one has that 
\begin{equation}\label{repres1}
\widetilde{\mathcal L_1}v =-\frac{d^2}{dx^2} v+\omega v-3\varphi^2_{\omega,Z} v,\quad
\widetilde{\mathcal L_2}v=-\frac{d^2}{dx^2}v+\omega v-\varphi^2_{\omega, Z}v.
\end{equation}
\end{lema}

\begin{proof} Since the proof of $\widetilde{\mathcal L_2}$ is similar to the one of $\widetilde{\mathcal L_1}$, we only deal with $\widetilde{\mathcal L_1}$. We decompose  the form $\mathcal Q_{\omega, Z}^1$ as $\mathcal Q_{\omega, Z}^1= \mathcal Q_{Z}^1 + \mathcal Q_{\omega}^1$ with $\mathcal Q_{Z}^1:  H^1_{\text{per}}([-L, L])\times  H^1_{\text{per}}([-L, L])\to \mathbb R$ and $\mathcal Q_{\omega}^1:  L^2_{\text{per}}([-L, L])\times  L^2_{\text{per}}([-L, L])\to \mathbb R$ defined by
\begin{equation}\label{repres2}
\begin{aligned}
\mathcal Q_{Z}^1(v,z)&=\int_{-L}^Lv_x z_xdx  -Zv(0)z(0),\\
\mathcal Q_{\omega}^1(v,z)&=\omega \int_{-L}^Lvzdx -\int_{-L}^L3\varphi_{\omega, Z}^2vzdx.
\end{aligned}
\end{equation}
We denote by $\mathcal T_1$ (resp. $\mathcal T_2$) the self-adjoint operator on 
$L^2_{per}([-L,L])$ (see Kato \cite{kato}, VI. Section 2.1) associated with $\mathcal Q_{Z}^1$ (resp. $\mathcal Q_{\omega}^1$). thus,  $D(\mathcal T_2)=L^2_{per}([-L,L])$ and $D(\mathcal T_1)=D(\widetilde{\mathcal L_1})$. We claim  that $\mathcal T_1$ is a self-adjoint extension of the operator $A^0$ defined in Lemma \ref{sp0}. Let $v\in H^2_{\text{per}}([-L, L])$ such that $v(0)=0$, and  define $w\equiv - v_{xx}\in L^2_{\text{per}}([-L, L])$. Then for every $z\in  H^1_{\text{per}}([-L, L])$ we have $ \mathcal Q_{Z}^1(v,z)=(w,z)$. 
Thus,  $v\in D(\mathcal T_1)$ and $\mathcal T_1 v=w=-\frac{d^2}{dx^2} v$. Hence, $
A^0\subset \mathcal T_1$. So, using Theorem \ref{self}  there exists $\beta\in \mathbb R$ such that $D(\mathcal T_1)=D(-\Delta_{\beta})$
which yields the claim. Next we shall show that $\beta=-Z$. Take $v\in D(\mathcal T_1)$ with $v(0)\neq 0$. Following the ideas in  \eqref{spec5} we obtain
\begin{equation*}
(\mathcal T_1 v,v)=[v'(0+)-v'(0-)]v(0) +\int_{-L}^L |v_{x}|^2dx=\int_{-L}^L |v_{x}|^2dx+\beta [v(0)]^2,
\end{equation*}
which should be equal to $\mathcal Q_{Z}^1(v,v)=\int_{-L}^L|v_x|^2dx-Z[v(0)]^2$.
Therefore $\beta=-Z$, and the lemma is proved.
\end{proof}

\item[\bf{(III)}] By Lemma \ref{repres} we can drop the tilde over $\widetilde{\mathcal L_1}$ and $\widetilde{\mathcal L_2}$ and  work with the operators  $\mathcal L_{1,Z}$ and $\mathcal L_{2,Z}$. The following Lemma verifies the analyticity of the families of operators $\mathcal L_{i, Z}$. 

\begin{lema}\label{analici} As a function of $Z$, $(\mathcal L_{1,Z})$ and $(\mathcal L_{2,Z})$ are two real-analytic families of self-adjoint operators of type (B) in the sense of Kato.
\end{lema}

\begin{proof} From Lemma \ref{repres}, Theorem VII-4.2 in \cite{kato}, it suffices to prove that the families of bilinear forms $(\mathcal Q_{\omega, Z}^1)$ and $(\mathcal Q_{\omega, Z}^2)$ defined in \eqref{spec14} are  real-analytic family of type (b). Indeed, since the form domains of these families are the same, namely $H^1_{per}$, for every $Z\in\mathbb R$, it is enough to prove that they are bounded  from below and closed, and that for any $v\in H^1_{per}$ the function $Z\to \mathcal Q_{\omega, Z}^i(v,v)$ is analytic. It is immediate that they are bounded  from below and closed. From the decomposition  of $\mathcal Q_{\omega, Z}^1$ into $\mathcal Q_{Z}^1$ and $\mathcal Q_{\omega}^1$, we see that $Z\to (\mathcal Q_{Z}^1v, v)$ is  real-analytic. From  Theorems \ref{s51b} and  \ref{negaZ3} we also have that $Z\to (\mathcal Q_{\omega}^1v, v)$ is real-analytic. The proof of the analyticity of the family $(\mathcal Q_{\omega, Z}^2)$ is similar to the one of $(\mathcal Q_{\omega, Z}^1)$.  
\end{proof}

{\bf{Remarks:}} \begin{enumerate}
\item The explicit resolvent  formula for $-\Delta_{-Z}$ in \eqref{resol2a} can be used to give another proof of the fact that  the families $(\mathcal L_{i,Z})$ are  real-analytic  in the sense of Kato.

\item We note  from   Theorems \ref{s51b} and \ref{negaZ3} that for $\omega>\pi^2/{2L^2}$  and $v, w\in H^1_{per}([0,L])$, 
\begin{equation}\label{quad0}
\begin{aligned}
\mathcal Q_{\omega, 0}^1(v,w)=\lim_{Z\to 0} \mathcal Q_{\omega, Z}^1(v,w)&=\int_{-L}^Lv_x w_x +\omega \int_{-L}^Lvw  -\int_{-L}^L3\phi_{\omega, 0}^2vw\\
\mathcal Q_{\omega, 0}^2(v,w)=\lim_{Z\to 0}\mathcal Q_{\omega, Z}^2(v,w)&=\int_{-L}^Lv_x w_x +\omega \int_{-L}^Lvw  -\int_{-L}^L\phi_{\omega, 0}^2vw.
\end{aligned}
\end{equation}
Here $\mathcal Q_{\omega, 0}^1$ is the bilinear form associated to the linear operator $\mathcal L_0$ defined in \eqref{L0}.
\end{enumerate}

\vskip0.1in

The following result of Angulo in \cite{angulo1} gives a precise description of the spectrum of the self-adjoint operator
\begin{equation}\label{L0}
\mathcal L_0\zeta \equiv -\frac{d^2}{dx^2}\zeta+\omega \zeta-3\phi^2_{\omega, 0}\zeta,
\end{equation}
on  $L^2_{per}([0,2L])$ and with domain $H^2_{per}([0,2L])$. Here $\omega>\frac{\pi^2}{2L^2}$ and $\phi_{\omega, 0}$ is the dnoidal traveling wave in \eqref{dnoidal2} which we want to perturb. 

\begin{lema}\label{specL0}  The operator $\mathcal L_0$ has exactly one negative simple isolated first eigenvalue $\tau_0$. The second eigenvalue is zero, and it is simple with associated eigenfunction $\frac{d}{dx}\phi_{\omega, 0}$. The rest of the spectrum is positive and discrete.
\end{lema}

{\bf{Remark:}} The Lemma \ref{specL0} can also be shown by using the method developed  by Angulo\&Natali in \cite{AnguloNatali1}.

\vskip0.2in

\begin{lema}\label{perteigen}  There  exist $Z_0>0$ and two analytic functions $\Pi : (-Z_0,Z_0)\to \mathbb R$ and $\Omega: (-Z_0,Z_0)\to L^2_{per}$ such that
\begin{enumerate}
\item[(i)] $\Pi(0)=0$ and $\Omega(0)=\frac{d}{dx}\phi_{\omega,0}$.

\item[(ii)] For all $Z\in (-Z_0,Z_0)$, $\Pi(Z)$ is the simple isolated second eigenvalue of $\mathcal L_{1,Z}$ and $\Omega(Z)$ is an associated eigenvector for $\Pi(Z)$.

\item[(iii)] $Z_0$ can be chosen small enough such that, except the two first eigenvalues, the spectrum of $\mathcal L_{1, Z}$ is positive.
\end{enumerate}
\end{lema}

\begin{proof} From Lemma \ref{specL0} we separate the spectrum $\sigma(\mathcal L_0)$ of the operator $\mathcal L_0$  in  \eqref{L0} into two parts $\sigma_0=\{\tau_0, 0\}$ and  $\sigma_1$ by a closed curve  $\Gamma$ (for example a circle) such that $\sigma_0$ belongs to the inner domain of $\Gamma$ and 
 $\sigma_1$ to the outer domain of $\Gamma$ (note that $\sigma_1\subset (a, +\infty)$ for $a>0$). From Lemma \ref{analici} follows that 
$\mathcal L_{1,Z}$ converges to $\mathcal L_0$ as $Z\to 0$ in the generalized sense, and so from Theorem IV-3.16 in \cite{kato} we have that $\Gamma\subset \rho(\mathcal L_{1,Z})$ for sufficiently small $|Z|$ and $\sigma (\mathcal L_{1,Z})$ is likewise separated by $\Gamma$ into two parts so that the part of $\sigma (\mathcal L_{1,Z})$ inside $\Gamma$ consists of a finite system of eigenvalues with total multiplicity (algebraic) two (we recall that zero is not eigenvalue of $\mathcal L_{1,Z}$). Next,  for $\epsilon$ small enough we consider the contours $\Gamma_1(\tau_0)=\{z\in \mathbb C: |z-\tau_0|<\epsilon\}$ and $\Gamma_2(0)=\{z\in \mathbb C: |z|<\epsilon\}$ such that $\Gamma_1(\tau_0)\cap \Gamma_2(0)=\emptyset $ and the only points of $\sigma(\mathcal L_0)$ in the  inner domain of $\Gamma_i$ are $\tau_0$ and $0$. Therefore from the nondegeneracy of  $\tau_0$ and $0$ we obtain from the Kato-Rellich Theorem (see Teorem XII.8 in \cite{RS}) the existence of two analytic functions $\Pi, \Omega$ defined in a neighborhood of zero such that we obtain the items (i), (ii) and (iii). This completes the proof of the Lemma.
\end{proof}

Next we shall study how the perturbed second eigenvalue $\Pi(Z)$ changes depending on the sign of $Z$.  For $Z$ small we have the following picture.

\begin{lema}\label{signeigen} There exists $0<Z_1<Z_0$ such that $\Pi(Z)<0$ for any $Z\in (-Z_1,0)$ and $\Pi(Z)>0$ for  any $Z\in (0, Z_1)$. Therefore, for $Z$ negative and small $\mathcal L_{1,Z}$ has exactly two negative eigenvalues and for  $Z$ positive and small $\mathcal L_{1,Z}$ has exactly one negative eigenvalue.
\end{lema}

\begin{proof} From Taylor's theorem we can write   the functions $\Pi$ and $\Omega$ of Lemma \ref{perteigen} around zero as
\begin{equation}\label{decomp1}
\begin{aligned}
\Pi(Z)&=\beta Z+ O(Z^2),\\
\Omega(Z)&=\phi'_{\omega,0}+ Z \psi_0 +  O(Z^2)
\end{aligned}
\end{equation}
where $\phi'_{\omega,0}=\frac{d}{dx} \phi_{\omega,0}$, $\beta\in \mathbb R$ ($\beta=\Pi'(0)$) and $\psi_0\in L^2_{per}$ ($\psi_0=\Omega'(0)$).  The desired result will follow  if we show that $\beta>0$.  From Theorems \ref{conver1}, \ref{s51b} and \ref{negaZ3} there exists $\chi_0\in H^1_{per}$  such that for $Z$ close to zero
\begin{equation}\label{decomp2}
\varphi_{\omega,Z}=\phi_{\omega,0}+Z\chi_0 +O(Z^2).
\end{equation}
Now, from    \eqref{p2} one has that for all $\psi\in H^1_{per}$ 
\begin{equation}\label{decomp3}
\langle -\varphi''_{\omega,Z}+\omega \varphi_{\omega,Z}-\varphi^3_{\omega,Z},\psi\rangle=Z\varphi_{\omega,Z}(0)\psi(0).
\end{equation}
So, inserting \eqref{decomp2} into \eqref{decomp3} and  differentiating with respect to $Z$, we obtain
\begin{equation}\label{decomp4}
\langle \mathcal L_0\chi_0,\psi\rangle=\phi_{\omega,0}(0)\psi(0)+ O(Z).
\end{equation}

We develop $\beta$ with respect to $Z$. We compute $\langle \mathcal L_{1,Z} \Omega(Z), \phi'_{\omega,0}\rangle$ in two different ways. 
\begin{enumerate}

\item Since $\mathcal L_{1,Z} \Omega(Z)=\Pi(Z) \Omega(Z)$ it follows from \eqref{decomp1} that
\begin{equation}\label{decomp5}
\langle \mathcal L_{1,Z} \Omega(Z),  \phi'_{\omega,0}\rangle=\beta Z\|\phi'_{\omega,0}\|^2 +O(Z^2).
\end{equation}

\item Since $\mathcal L_{1,Z} $ is  self-adjoint and $\phi'_{\omega,0}\in \mathcal D(\mathcal L_{1,Z})$ (in view of  $ \phi'_{\omega,0}\in H^n_{per}$ for all $n$ and $\phi_{\omega,0}$ is even), we obtain $\langle\mathcal L_{1,Z} \Omega(Z), \phi'_{\omega,0}\rangle=\langle\Omega(Z), \mathcal L_{1,Z}\phi'_{\omega,0}\rangle$.  Thus, from \eqref{decomp2},
\begin{equation}\label{decomp6}
\begin{aligned}
 \mathcal L_{1,Z} \phi'_{\omega,0}&= \mathcal L_0(\phi'_{\omega,0})+3(\phi^2_{\omega,0}-\varphi^2_{\omega,Z})\phi'_{\omega,0}=3(\phi^2_{\omega,0}-\varphi^2_{\omega,Z})\phi'_{\omega,0}\\
 &=-6Z \phi_{\omega,0} \phi'_{\omega,0} \chi_0+O(Z^2).
\end{aligned}
\end{equation}
Hence,  from \eqref{decomp1} and \eqref{decomp6} it follows that 
\begin{equation}\label{decomp7}
\langle \mathcal L_{1,Z} \Omega(Z), \phi'_{\omega,0}\rangle=- 6Z\langle \phi'_{\omega,0}, \chi_0\phi_{\omega,0}\phi'_{\omega,0}\rangle+O(Z^2).
\end{equation}
It is easy to see that
\begin{equation}\label{decomp8}
\mathcal L_{0}(\omega \phi_{\omega,0}-\phi^3_{\omega,0})=6\phi_{\omega,0}(\phi'_{\omega,0})^2,
\end{equation}
which combined with \eqref{decomp7} gives us the last equality
\begin{equation}\label{decomp9}
\begin{aligned}
\langle\mathcal L_{1,Z} \Omega(Z),  \phi'_{\omega,0}\rangle&=-Z\langle\mathcal L_{0} \chi_0,\omega \phi_{\omega,0}-\phi^3_{\omega,0}\rangle+O(Z^2)\\
&=-Z[\omega \phi^2_{\omega,0}(0)-\phi^4_{\omega,0}(0)]+O(Z^2).
\end{aligned}
\end{equation}
\end{enumerate}

Finally, a combination of \eqref{decomp5} and \eqref{decomp9} leads to 
\begin{equation}\label{decomp10}
\beta=-\frac{\omega \phi^2_{\omega,0}(0)-\phi^4_{\omega,0}(0)}{\|\phi'_{\omega,0}\|^2}+O(Z).
\end{equation}

Now, from Theorem \ref{conver0} we have $\phi_{\omega,0}(0)\in (0,\sqrt{\omega})$ and so $\beta>0$ for $Z$ small the same holds. Hence, the first equality in \eqref{decomp1} completes the proof.
\end{proof}

{\bf{Remark:}} The proof of Lemma \ref{signeigen} also  shows  the eigenvalue-mapping $Z\to \Pi(Z)$ is a strictly increasing function in a neighborhood of zero.

\vskip0.2in

\item[\bf{(IV)}]  Now we are in position for counting the number of negative eigenvalues of $\mathcal L_{i,Z}$ for all $Z$. using  a classical continuation argument based on the Riesz-projection. We denote the number of negatives eigenvalues of $\mathcal L_{i,Z}$ by $n(\mathcal L_{i,Z})$.

\begin{lema}\label{spec6a} Let $\omega$ such that $\omega>\frac{\pi^2}{2L^2}$ and $\omega>Z^2/4$. Then 
\begin{enumerate}
\item for $Z>0$,  $n(\mathcal L_{1,Z})=1$,

\item for $Z<0$, $n(\mathcal L_{1,Z})=2$.
\end{enumerate}
\end{lema}

\begin{proof} Let $Z<0$ and define $Z_\infty$ by
$$
Z_\infty=\inf \{z<0: \mathcal L_{1,Z}\;{\text{has exactly two negative eigenvalues for all}}\; Z\in (z,0)\}.
$$
From Lemma \ref{signeigen} one has that  $\mathcal L_{1,Z}$ has exactly two negative eigenvalues for all $Z\in (Z_1,0)$, so $Z_\infty$ is well defined and $Z_\infty\in [-\infty,0)$. We claim that $Z_\infty=-\infty$. Suppose that $Z_\infty> -\infty$. Let $N( \mathcal L_{1,Z_{\infty}})$ and $\Gamma$ a closed curve (for example a circle or a rectangle) such that $0\in \Gamma\subset \rho(\mathcal L_{1,Z_{\infty}})$ and such that all the negatives eigenvalues of  $\mathcal L_{1,Z_{\infty}}$ belong to the inner domain of $\Gamma$. From Lemma \ref{analici} it follows  that $\mathcal L_{1,Z}\to \mathcal L_{1,Z_{\infty}}$ as $Z\to Z_{\infty}$ in the generalized sense, and so there is a $\delta>0$ such that for $Z\in [Z_{\infty}-\delta, Z_{\infty}+\delta]$ we have $\Gamma\subset \rho(\mathcal L_{1, Z})$ and $\rho(\mathcal L_{1,Z})$ is likewise separated by $\Gamma$ into two parts so that the part of $\sigma(\mathcal L_{1, Z})$ inside $\Gamma$ consists of a system of eigenvalues with total multiplicity (algebraic) equal to $N$. This conclusion follows from the existence of an analytic family of Riesz-projections, $Z\to P(Z)$,  given by
$$
P(Z)=-\frac{1}{2\pi i}\int_{\Gamma} (\mathcal L_{1, Z}-\xi)^{-1}d\xi,
$$
which implies  that 
\begin{equation}\label{decomp11}
{\rm{dim}} ({\rm{Rank}}\; P(Z))={\rm{dim}} ({\rm{Rank}}\; P(Z_\infty))=N,\;\;\forall\;\,Z\in [Z_{\infty}-\delta, Z_{\infty}+\delta].
\end{equation}
{\it We observe that we can choose $\Gamma$ independently of the parameter $Z$} (see Remark below). Now by definition of $Z_{\infty}$, there exists $z_0$ such that $Z_\infty<z_0<Z_\infty+\delta$ and $\mathcal L_{1,Z}$ has exactly two negative eigenvalues for all $Z\in (z_0,0)$. Therefore $\mathcal L_{1, Z_\infty+\delta}$ has two negative eigenvalues and from \eqref{decomp11} it follows that $N=2$ and so $\mathcal L_{1,Z}$ has two negative eigenvalues for $Z\in (Z_{\infty}-\delta, 0)$ contradicting the definition of $Z_{\infty}$. Therefore, we have established the claim$Z_{\infty}=-\infty$. A similar analysis is applied to the case $Z>0$. This finishes the proof of the lemma.
\end{proof}

{\bf{Remark:}}  We can choose $\Gamma$ independently of the parameter $Z<0$ in the  beginning  of the proof of Lemma \ref{spec6a} in the following manner : since for all $Z$, $\varphi_{\omega,Z}\leqq \eta_{1, +}\leqq \sqrt{2\omega}$, for $\|f\|=1$ and  $f\in \mathcal D$
$$
\langle \mathcal L_{1,Z} f, f\rangle \geqq -3\int \varphi^2_{\omega,Z}f^2dx\geqq -6\omega.
$$
Therefore, $\inf \sigma(\mathcal L_{1,Z})\geqq -6\omega$ for all $Z<0$. So, $\Gamma$ can be chosen as the rectangle $\Gamma=\partial R$ for $R$ being
$$
R=\{ z\in\mathbb C: z=z_1+iz_2, (z_1,z_2)\in [-6\omega-1,0]\times [-a,a],\;\text{for some}\;a>0\}.
$$

\begin{lema}\label{spec6c} The function $\Omega(Z)$ defined in Lemma \ref{perteigen} and associated to the second negative eigenvalue of $ \mathcal L_{1, Z}$ can be extended to $(-\infty,\infty)$. Moreover, $\Omega(Z)\in H^1_{per}$ is an odd function for $Z\in (-\infty,\infty)$.
\end{lema}

\begin{proof}  From Lemma \ref{analici}  and Theorem XII.7 in \cite{RS} the set $\Gamma_0=\{(Z, \lambda)| Z\in \mathbb R, \lambda\in \rho( \mathcal L_{1, Z})\}$ is open and 
$$
(Z, \lambda)\in \Gamma_0\to  (\mathcal L_{1, Z}-\lambda)^{-1}
$$
is a holomorphic function in both variables. So, we can repeat the argument of
Lemma \ref{perteigen} at each point $Z$ and on each neighborhood of $Z$ to see that the functions $\Omega(Z)$ and $\Pi(Z)$ are holomorphic for every $Z\in \mathbb R$. Next we consider $Z<0$ (the case $Z>0$ is similar). We know from Lemma \ref{pone} and Lemma \ref{perteigen} that the eigenvectors  $\Omega(Z)$ are even or odd and $\Omega(0)=\frac{d}{dx}\phi_{\omega,0}$ is odd. Then, from the equality
$$
\lim_{Z\to 0}\langle \Omega(Z), \Omega(0)\rangle=\|\Omega(0)\|^2\neq 0,
$$
one has that $\langle \Omega(Z), \Omega(0)\rangle \neq 0$ for $Z$ close to $0$. Thus  $\Omega(Z)$ is odd. Let $z_\infty$ be
$$
z_\infty=\{z<0: \Omega(Z)\;\text{is odd for any}\;Z\in (z,0]\}.
$$
Suposse now that $z_\infty>-\infty$. If $\Omega(z_\infty)$ is odd, then by continuity  there exists $\delta>0$ such that $\Omega(z_\infty-\delta)$ is odd which is a contradiction. Thus Lemma \ref{pone} implies that $\Omega(z_\infty)$ is even. Now, since  $\Omega(z_\infty)$  is the limit of odd functions we obtain that $\Omega(z_\infty)$ is odd. Therefore $\Omega(z_\infty)\equiv 0$, which is a contradiction because  $\Omega(z_\infty)$ is an eigenvector. This concludes the proof of the Lemma.
\end{proof}

\subsection{Convexity condition}

Here, we shall prove the increasing property of the mapping $\omega\to \|\varphi_{\omega, Z}\|^2$, for all $Z$,  which suffices  for our stability/instability results for the orbit defined in \eqref{orb}. For technical reasons we can only show this property for $\omega$ large. But we believe that this property should be true  for every $\omega$ admissible.

\begin{teo}\label{con1}  Let $Z\in \mathbb R-\{0\}$, $\omega>Z^2/4$ and $\omega$ large. Then for the dnoidal-peak smooth  curve $\omega\to \varphi_{\omega, Z}$ given in \eqref{pp} we have
$$
\frac{d}{d\omega}\|\varphi_{\omega, Z}\|^2>0.
$$
\end{teo}

\begin{proof} For $Z>0$ we have $\varphi_{\omega, Z}=\phi_{\omega, Z}$. Then via a change of variable and from Theorem \ref{s51b} we have for $a=a(\omega)$, $\eta_1= \eta_{1, +}$, $k=k(\omega)$ and $K-a=\frac{\eta_1}{\sqrt{2}}L$ the equality
\begin{equation}\label{conc2}
\begin{aligned}
\|\phi_{\omega, Z}\|^2&=\eta_1^2\int_{-L}^L dn^2\Big(\frac{\eta_1}{\sqrt 2} |\xi| +a;k\Big)d\xi=
2\sqrt{2}\eta_1\int_{a}^{K(k)} dn^2(y;k)dy\\
&=2\sqrt{2}\eta_1[E(k)-E(a)]=2\sqrt{2}\eta_1[E(k)-E(\varphi_a,k)].
\end{aligned}
\end{equation}
Here $E(\varphi_a,k)$  is the normal elliptic integral of the second kind defined for  $\sin \varphi_a= sn(a)$ by
\begin{equation}\label{con3}
E(\varphi_a,k)=\int_{0}^{\varphi_a}\sqrt{1-k^2 \sin^2 \theta}\;d\theta=\int_{0}^a dn^2 (u;k)\;du=E(a),
\end{equation}
and  $E(k)=E(\pi/2,k)$. Next, we consider  the identity
\begin{equation}\label{con2}
\begin{aligned}
\frac{d}{d\omega}\|\phi_\omega\|^2&=2\sqrt{2}\frac{d\eta_1}{d\omega}[E(k)-E(\varphi_a,k)]\\
&+2\sqrt{2}\eta_1\Big[\Big(E'(k)-\frac{\partial E}{\partial k}\Big) \frac{dk}{d\omega}-\frac{\partial E}{\partial \varphi_a}\frac{d\varphi_a}{d\omega}\Big].
\end{aligned}
\end{equation}

We shall calculate the differentiation terms in \eqref{con2}.
\begin{enumerate}

\item From \eqref{con3} one has that $
\frac{\partial E}{\partial \varphi_a}(\varphi_a,k)=\sqrt{1-k^2 sn^2(a)}=dn(a)$.

\item From (\cite{byrd}) we obtain
$$
\frac{\partial E}{\partial k}(\varphi_a,k)= \frac{E(\varphi_a,k)-F(\varphi_a,k)}{k}=\frac{E(a)-a}{k},
$$
where $F(\varphi_a,k)$ is the normal elliptic integral of the first kind such that for 
$\sin \varphi_a= sn(a)$ it follows that $F(\varphi_a,k)=a$.

\item Next, since $sn(u+K)=\frac{cn (u)}{dn (u)}\equiv cd(u)$ one has that
$\varphi_a(\omega)=\sin ^{-1}[cd(\eta_1L/\sqrt{2})]$. So,
\begin{equation}\label{con5}
\frac{d}{d\omega}\varphi_a=\frac{dn}{k'sn} \frac{d}{d\omega}cd\Big(\frac{\eta_1}{\sqrt{2}}L;k).
\end{equation}
Now, from using \cite{byrd} again one finds that 
$$
\begin{aligned}
\frac{d}{d\omega}cd\Big(\frac{\eta_1}{\sqrt{2}}L;k)&=\frac{L}{\sqrt{2}}\frac{\partial}{\partial u} cd\Big(\frac{\eta_1}{\sqrt{2}}L;k\Big)\frac{d\eta_1}{d\omega}+\frac{\partial}{\partial k} cd \Big(\frac{\eta_1}{\sqrt{2}}L;k\Big)\frac{dk}{d\omega}\\
&=-\frac{k'^2 L}{\sqrt{2}}\frac{d\eta_1}{d\omega}\frac{sn}{dn^2}+\frac{sn}{kdn^2}\Big[E\Big(\frac{\eta_1}{\sqrt{2}}L\Big)-k'^2\frac{\eta_1}{\sqrt{2}}L\Big]\frac{dk}{d\omega}.
\end{aligned}
$$
So, from \eqref{con5} and from the equality $dn(u+K)=k'/{(dn u)}$
\begin{equation}\label{con6}
\frac{d}{d\omega}\varphi_a=dn a\Big[- \frac{ L}{\sqrt{2}}\frac{d\eta_1}{d\omega}+ 
\frac{1}{kk'^2}\Big[E\Big(\frac{\eta_1}{\sqrt{2}}L\Big)-k'^2\frac{\eta_1}{\sqrt{2}}L\Big]\frac{dk}{d\omega}\Big ]
\end{equation}

\item Combining the identities
$$
\frac{d}{dk}K(k)=\frac{E(k)-k'^2K(k)}{kk'^2},\;\;\;\;\frac{ L}{\sqrt{2}}\frac{d\eta_1}{d\omega}=\frac{d}{dk}K(k)\frac{dk}{d\omega}-a'(\omega),
$$
and
$$
E\Big(\frac{\eta_1}{\sqrt{2}}L\Big)-E(k)+k'^2 a=\int_{K-a}^K[k'^2-dn^2(u)]du=-k^2\int_{K-a}^Kcn^2(u)du
$$
it follows that
\begin{equation}\label{con7}
\begin{aligned}
\frac{d}{d\omega}\varphi_a=dn(a) \Big [a'(\omega)-\frac{k}{k'^2}\int_{K-a}^Kcn^2(u)du \frac{dk}{d\omega}\Big]\equiv dn(a) A(\omega).
\end{aligned}
\end{equation}
We observe  that $ A(\omega)<0$ and so $\frac{d}{d\omega}\varphi_a<0$.
\end{enumerate}

Then, gathering the information \eqref{con2} and from (1)-(4) above we obtain that
\begin{equation}\label{con8}
\begin{aligned}
&\frac{d}{d\omega}\|\phi_\omega\|^2=\frac{4}{L}\Big[K'(k)[E(k)-E(a)]+E'(k)[K(k)-a]\Big ]\frac{dk}{d\omega}\\
&\;\;-\frac{4}{L} a'(\omega)[E(k)-E(a)] +\frac{4}{L}[K(k)-a]\frac{a-E(a)}{k}\frac{dk}{d\omega} -2\sqrt{2}\eta_1 dn^2(a) A(\omega).
\end{aligned}
\end{equation}

Now, since 
$$
a-E(a)=\int_0^a [1-dn^2(u)]du=k^2\int_0^a sn^2(u)du>0,\;\;\;\;
E(k)-E(a) >0,
$$
$a'(\omega)<0$, $A(\omega)<0$ and $\frac{dk}{d\omega}>0$ we obtain that the expression on the second line in \eqref{con8} is positive. Therefore from \eqref{con8} one concludes that
\begin{equation}\label{con9}
\begin{aligned}
\frac{L}{4}\frac{d}{d\omega}\|\phi_{\omega, Z}\|^2&>
 \frac{d}{d\omega}[K(k)E(k)]-E(a)\frac{d}{d\omega} K(k)-a\frac{d}{d\omega}E(k)\\
&>\frac{d}{d\omega}[K(k)E(k)]-a\frac{d}{d\omega}[K(k)+E(k)]\\
&>\frac{d}{d\omega}[K(k)E(k)-\frac12(K(k)+E(k))]
\end{aligned}
\end{equation}
where $\omega$ is chosen large enough such that $a(\omega)\leqq \frac12$. We note that here we have used that the mapping $k\to K(k)+E(k)$ is increasing and so $\frac{d}{d\omega}[K(k)+E(k)]=\frac{d}{dk}[K(k)+E(k)]\frac{dk}{d\omega}>0$. Since 
$$
\frac{d}{dk}\Big[K(k)E(k)-\frac12(K(k)+E(k))\Big]>0,
$$
it follows from \eqref{con9} that $\frac{d}{d\omega}\|\phi_{\omega, Z}\|^2>0$ for $\omega$ large.

Next, we consider the case $Z<0$. For  $\varphi_{\omega, Z}=\zeta_{\omega, Z}$ and  $\beta=\sqrt{2}/{\eta_1}$ one has that 
\begin{equation}\label{zcon}
\begin{aligned}
\|\zeta_{\omega, Z}\|^2&=\eta_1^2\int_{-L}^L dn^2\Big(\frac{\eta_1}{\sqrt 2} |\xi| -a\Big)d\xi=\frac{4}{\beta}\int_{-a}^{\frac{L}{\beta}-a} dn^2(y)dy\\
&=\frac{4}{\beta}\int_{-a}^{K} dn^2(y)dy\equiv G(\beta),
\end{aligned}
\end{equation}
using that  $K+a=\frac{\eta_1}{\sqrt{2}}L$. So,
\begin{equation}\label{zcon1}
\frac{d}{d\omega}\|\zeta_{\omega, Z}\|^2=G'(\beta)\frac{d\beta}{d\omega}=-\frac{\sqrt{2}}{\eta^2_1}\frac{d\eta_1}{d\omega}G'(\beta),
\end{equation}
where
\begin{equation}\label{zcon2}
G'(\beta)=4\beta^{-2}\Big [ -\int_{-a}^{K} dn^2(y)dy +\beta \frac{d}{d\beta}\int_{-a}^{K} dn^2(y)dy\Big]\equiv 4\beta^{-2}H(\beta).
\end{equation}
The idea now is to show that $H(\beta)<0$. Indeed, from Section 5 we have   $\omega \to \eta_2(\omega)$ is a positive decreasing function, then  for $\omega\to +\infty$ follows $\eta_2^2/{2\omega}\to 0$. So, Theorem \ref{negaZ3} implies that $k^2\to 1$ and $\eta_1^2/{2\omega}\to 1$ for $\omega\to +\infty$. Thus, $\beta\to 0$ as  $\omega\to +\infty$. Hence, $a(\beta)=a(\eta_1^{-1}(\sqrt{2}/\beta))\to 0$ as $\beta\to 0$ (see Corollaries \ref{s62a} and \ref{s62b}). Since $dn(x;1)=sech(x)$ and $K(1)=+\infty$ we obtain 
\begin{equation}\label{con10}
H(0)=-\int_0^{\infty}sech^2(y)dy<0.
\end{equation}
Therefore $H(\beta)<0$  for $\beta$ close to zero. This completes the proof of th Theorem.
\end{proof}

\subsection{Stability results}

From the last subsections our stability results associated to the orbit in \eqref{orb1} generated by the dnoidal-peak solution profile $\varphi_{\omega,Z}$ in \eqref{pp} can be now established. As it was pointed the abstract theory of Grillakis, Shatah and Strauss \cite{grillakis2} shall be use, and so  we briefly discuss the criterion for obtaining stability or instability in our case. Consider the linear operator $H_{\omega, Z}$ defined in \eqref{spec0} and  denote by $n(H_{\omega, Z})$ the number of negative eigenvalues of $H_{\omega, Z}$. Define
\begin{equation}\label{st}
p_Z(\omega_0)=
\begin{cases}
\begin{aligned}
&1, \quad  {\text{if}} \;\; \partial_\omega\|\varphi_{\omega,Z}\|^2>0,\;\;at\;\; \omega=\omega_0,\\
&0, \quad  {\text{if}} \;\; \partial_\omega\|\varphi_{\omega,Z}\|^2<0,\;\;at\;\; \omega=\omega_0.
\end{aligned}
\end{cases}
\end{equation}
Then, having established the Assumption 1, Assumption 2 and Assumption 3 of \cite{grillakis2}, namely, the existence of global solutions (Proposition \ref{lwp1}), the existence of a smooth curve of standing-wave, $\omega \to \varphi_{\omega, Z}$ (Theorem \ref{s51b} - Theorem \ref{negaZ3}), and $Ker(\mathcal L_{1,Z})=\{0\}$, $Ker(\mathcal L_{2,Z})=[\varphi_{\omega, Z}]$, the next Theorem follows from the Instability Theorem and Stability Theorem in \cite{grillakis2}.
\begin{teo}\label{st0} 
Let $\omega_0>\frac{\pi^2}{2L^2}$ and  $\omega_0>\frac{Z^2}{4}$.
\begin{enumerate}
\item If $n(H_{\omega_0, Z})=p_Z(\omega_0)$, then the  dnoidal-peak standing wave $e^{i\omega_0 t}\varphi_{\omega_0,Z}$ is stable in $H^1_{per}([-L,L])$.

\item If $n(H_{\omega_0, Z})-p_Z(\omega_0)$ is odd, then the  dnoidal-peak standing wave $e^{i\omega_0 t}\varphi_{\omega_0,Z}$ is unstable in $H^1_{per}([-L,L])$.
\end{enumerate}
\end{teo}

Now we can prove our main result Theorem \ref{main}

\begin{proof} From Theorem \ref{con1} follows that $p_Z(\omega)=1$ for all $Z\in \mathbb R-\{0\}$ and $\omega$ large. Next, from Lemma \ref{spec6} we have that $\mathcal L_{2, Z}$ has zero as a simple eigenvalue and from Lemma \ref{spec6a1} we have $\mathcal L_{1, Z}$ has a trivial kernel. Thus, from Theorem \ref{st0}, Lemma \ref{spec6a} we obtain the item (1) and item (2).

Lemma \ref{spec6c} assures that the second eigenvalue of $\mathcal L_{1, Z}$ considered in the whole space $L^2_{per}([-L,L])$ is associated with an odd eigenfunction, and thus dissapears when the problem is restricted to subspace of even periodic functions. Moreover, since $\varphi_{\omega,Z}$ is an even function and trivially satisfies that $\langle  \mathcal L_{1, Z} \varphi_{\omega,Z},\varphi_{\omega,Z}\rangle<0$, for  $Z<0$, we obtain that the first negative eigenvalue of $\mathcal L_{1, Z} $ is still present when the problem is restricted to the subspace of even periodic function of $H^1_{per}([-L,L])$, namely, $H^1_{per, even}([-L,L])$. So we obtain in this case that $n(H_{\omega, Z}|_{H^1_{per, even}([-L,L])})=1$. Therefore item (3) of the Theorem follows from item (1) of Theorem \ref{st0} and Proposition \ref{lwp1}. This finishes the proof of the Theorem.
\end{proof}

\end{enumerate}
\section{Appendix}

We shall establish some properties of the function $\psi$ defined in \eqref{s10} and which has been used in subsection 5.2. Many of these properties are immediate and so we omit the proof. (a) If $\psi(0)<\frac{\pi}{2}$ then $\psi(r)<\frac{\pi}{2}$  for all $r$. (b) $\eta_1>\phi(\xi)$ implies that $\psi(\xi)\neq 0$ for all $\xi$. So, without loss of generality, we suppose $\frac{\pi}{2}>\psi(\xi)> 0$. (c) Since $\phi$ is $1$-periodic then $\psi$ is also $1$-periodic. (d) Zeros of $\psi'$:  from \eqref{s11} and  (a)-(b) follows that $\psi'(\xi)=0$ if and only if $\phi'(\xi)=0$.  (e) $\phi'(\xi)=0$ if and only if $\phi(\xi)=\eta_2$. (f) There is a unique $s\in (0,1)$ such that $\phi(s)=\eta_2$. Indeed, consider $s<s_0$ s.t. $\phi(s_0)=\eta_2$. Then there is a $r\in [s,s_0]$ where $\phi(r)$ is a maximum. So $\phi(r)\geqq \phi(x)\geqq \eta_2$ for every $x\in [s, s_0]$. Since $\phi'(r)=0$ then $\phi(r)=\eta_2$. Therefore, $\phi(x)\equiv \eta_2$ for every $x\in [s, s_0]$. Then from \eqref{p3}-(3) it follows $\omega=\eta_2^2$ and so it follows from \eqref{s2} the equality $\eta_1=\eta_2$, which is a contradiction. (g) $\psi'(\xi)=0$ if and only if $\xi=s$, where $s$ is the unique point  in $(0,1)$ s.t. $\phi(s)=\eta_2$. (h) Let $s$ be such that $\phi(s)=\eta_2$ (so $s$ is a minimum for $\phi$) then
\begin{equation}\label{s18}
\sin^2\psi(s)=\frac{\eta_1^2-\eta_2^2}{\eta_1^2}.
\end{equation}
Hence from \eqref{s17} and \eqref{s18} $s$ is a maximum of $\psi$. Indeed,  for every $\xi\in (0,1)-\{s\}$,
$\sin^2\psi(\xi)< \frac{1}{\eta^2}=\sin^2\psi(s)$. Then, since $0<\psi(\xi)<\frac{\pi}{2}$ we obtain that $\psi(\xi)< \psi(s)$. (i) If $\phi$ is even then $\psi$ is also even. (j) For $Z>0$ it follows from \eqref{s4} the inequality $\phi'(0+)<0$ (so by evenness we have a peak in zero for $\phi$ in the form ``$\wedge$''). Now,  \eqref{s17} implies  $|\psi'(0+)|=|\psi'(0-)|$ and so from \eqref{s13}  $\psi'(0+)=-\psi'(0-)$. Therefore, \eqref{s12} implies that $
0<Z \phi^2(0)=\eta_1^2\psi'(0+)\sin2\psi(0)$, and so $\psi'(0+)>0$. Then for $\xi\in (0,s)$, $\psi'(\xi)>0$. By evenness we have a peak in  zero for $\psi$ in the form ``$\vee$''

\bigskip\noindent{\bf Acknowledgments:} J. Angulo was partially supported by CNPq/Brazil grant and CAPES/Brazil grant, and G. Ponce was supported by a NSF grant. This work started while J. A. was visiting the Mathematics Department of the University of California at Santa Barbara whose hospitality he gratefully acknowledges.


\begin{thebibliography}{22}

\bibitem{AS} Ablowitz, M.J. and Segur, H., \textit{Solitons and Inverse
Scattering}, SIAM Publication, (1981).

\bibitem{Ag} Agrawal, G.,  \textit{Nonlinear fiber optics},  Academic Press, (2001).

\bibitem{agfhr} Albeverio, S., Gesztesy, F., Hoegh-Krohn, R., Holden, H., \textit{Solvable models in quantum mechanics}, Texts and Monographs in Physics. Springer-Verlag, New York, (1988). 

\bibitem{ak} Albeverio, S. and  P. Kurasov, \textit{Singular Perturbations of Differential Operators}, London Mathematical Society, Lecture Note Series, 271, Cambridge University Press, (2000).

\bibitem{angulo1} Angulo, J., \textit{Non-linear stability of periodic travelling-wave equation for the Schr\"odinger and modified Korteweg-de Vries equation}, J. of
Differential Equations, 235 (2007), 1--30.

\bibitem{angulo4} Angulo, J., \textit{Nonlinear Dispersive Equations: Existence and Stability of Solitary and Periodic Travelling Wave Solutions},  Mathematical Surveys and Monographs (SURV), AMS, (2009).

\bibitem{angulo3} Angulo, J., Bona, J.L. and Scialom, M., \textit{Stability of cnoidal waves},
Advances in Differential Equations, 11 (2006), 1321--1374.

\bibitem{AnguloNatali1} Angulo, J. and Natali, F., \textit{Positivity properties and stability of
periodic travelling waves solutions}, SIAM, J. Math. Anal., 40 (2008),  1123--1151.

\bibitem{AnguloNatali2}  Angulo, J. and Natali, F., \textit{Stability and instability of periodic travelling wave solutions for the critical Korteweg-de Vries and nonlinear Schr\"odinger equations},  Phys. D, 238 (2009),  603--621.

\bibitem{Bou} Bourgain, J., \textit{Global solutions of
nonlinear Schr\"odinger equations}, American Mathematical Society Colloquium Publications, AMS, Providence, RI., 46, (1999).

\bibitem{BK}  Brazhnyi, V. A. and Konotop, V. V. \textit{Theory of nonlinear matter waves in optical lattices},  N. Akhmediev (Ed.). Dissipative Solitons. vol. 18, (2005)  627.

\bibitem{BR} Bronski, J. and Rapti, Z., \textit{Modulation instability for nonlinear Schr\"odinger equations with a periodic potential}, Dynamics of PDE, 2 (2005), 335--355.

\bibitem{byrd} Byrd, P.F. and Friedman, M.D., \textit{Handbook of elliptic integrals
for engineers and scientists}, 2nd ed., Springer, NY, (1971).

\bibitem{CMM}  Cai, D., McLaughlin, D. W. and  McLaughlin, K. T. R., \textit{The nonlinear Schr\"odinger Equation as both a PDE and a dynamical system}, In handbook of dynamical systems, North-Holland, Amsterdam, vol 2 (2002), 599--675.

\bibitem{CM}  Cao,  X. D. and Malomed, B. A., \textit{ Soliton-defect collisions in the nonlinear Schr\"odinger Equation},  Phys. Lett. A 206 (1995), 177--182.

\bibitem{Ca} Cazenave, T., \textit{Semilinear Schr\"odinger Equation}, Courant Lecture Notes in Mathematics, vol. 10, AMS, Courant Institute of Mathematical Science, (2003).

\bibitem{CL}  Cazenave, T. and Lions, P.-L., \textit{Orbital stability of standing waves for some nonlinear Schr\"odinger equations}, Comm. Math. Phys. 85 (1982), 549--561.

\bibitem{DH} Datchev, K. and Holmer, J., \textit{Fast soliton scattering by attractive delta impurities}, pre-print.

\bibitem{DMADDKK} Davis, K. B., Mewes, M.O., Andrews, M. R., van Druten, N. J., Durfee, D.S., Kurn, D.M. and Ketterle, W., \textit{Bose-Einstein condensation in gas of sodium atoms}, Phys. Rev. Lett., 74(22) (1995), 3969--3973.

\bibitem{Ea} Eastham, M.S.P., \textit{The Spectral Theory of Periodic Differential Equations}, Scottish Academic Press, London, UK, (1973).

\bibitem{fuje} Fukuizumi, R. and Jeanjean, L., \textit{Stability of standing waves for a nonlinear Schr\"odinger equation with a repulsive Dirac delta potential}, Discrete Contin. Dyn. Syst., 21 (2008),  121--136.

\bibitem{fuoht} Fukuizumi, R., Ohta, M., and Ozawa, T. \textit{Nonlinear Schr\"odinger equation with a point defect},  Ann. Inst. H. Poincar\'e Anal. Non Lin\'eaire, 25 (2008), 837--845. 

\bibitem{GH1}  Gallay, T. and H$\breve{\text{a}}$r$\breve{\text{a}}$gu\c s, M., \textit{Stability of small periodic waves for the nonlinear Schr\"odinger equation}, J. Differential Equations 234 (2007), 544--581.

\bibitem{GH2}   Gallay, T. and H$\breve{\text{a}}$r$\breve{\text{a}}$gu\c s, M., \textit{Orbital stability of periodic waves for the nonlinear Schr\"odinger equation},  J. Dyn. Diff. Eqns. 19 (2007), 825-865.

\bibitem{ghw} Goodman, R. H.,  Holmes,  J. and Weinstein, M., \textit{Strong NLS soliton-defect interactions}, Phys. D 192 (2004),  215--248.
\bibitem{grillakis1} Grillakis, M., Shatah, J. and Strauss, W., \textit{Stability theory
of solitary waves in the presence of symmetry I}, J. Functional
Anal., 74 (1987),  160-197.
\bibitem{grillakis2} Grillakis, M., Shatah, J. and Strauss, W., \textit{Stability theory
of solitary waves in the presence of symmetry II}, J. Functional
Anal., 94 (1990),  308-348.

\bibitem{HMZ1}  Holmer, J., Marzuola, J. and Zworski, M., \textit{Fast soliton scattering by delta impurities}, Comm. Math. Phys., 274(91) (2007), 187--216.

\bibitem{HMZ2}  Holmer, J., Marzuola, J. and Zworski, M., \textit{Soliton alignedting by external delta potentials},  J. Nonlinear Sci., 17(4) (2007),  349--367.

\bibitem{HZ1}  Holmer, J. and Zworski, M., \textit{Slow soliton interaction with external delta potentials},  J. Modern Dynam., 1 (2007),  689--718.

\bibitem{HZ2}  Holmer, J. and Zworski, M., \textit{Soliton interaction with slowly varying potentials},  IMRN, 2008, Article ID rnn026, 36 pages (2008).


\bibitem{kato} Kato, T., \textit{Perturbation Theory for Linear Operators}, 2nd edition, Springer, Berlin, 1984.

\bibitem{lffgs}  Le Coz, S., Fukuizumi, R., Fibich, G., Ksherim, B. and  Sivan, Y., \textit{Instability of bound states of a nonlinear Schr�dinger equation with a Dirac potential}, Phys. D 237 (2008), 1103--1128.

\bibitem{LP}  Linares, F. and Ponce, G., \textit{Introduction to Nonlinear Dispersive Equations}, Springer New York (2009) 

\bibitem{MA}  Ma, Y. C. and Ablowitz, M. J.,  \textit{The periodic cubic Schr�dinger equation}, Stud. Appl.  Math. 65(2) (1981), 113--158.




\bibitem{Mn}  Menyuk, C. R.,  \textit{Soliton robustness in optical fibers}, J. Opt. Soc. Am. B, 10(9)  (1993), 1585--1591.

\bibitem{MN} Moloney, J. and  Newell, A., \textit{Nonlinear optics}, Westview Press. Advanced Book Program, Boulder, 

\bibitem{RS}  Reed, S. and Simon, B., \textit{Methods of Modern Mathematical Physics: Analysis of Operator},  Academic Press, Vol. IV, 1978.

\bibitem{Row}  Rowlands, G., \textit{On the stability of solutions of nonlinear Schr\"odinger equation},  IMA J. Appl. Math.  13, (1974),  367--377.

\bibitem{ST}  Sakaguchi, H. and Tamura, M., \textit{Scattering and trapping of nonlinear Schr\"odinger solitons in external potentials},  J. Phys. Soc. Japan,  73, (2004),  2003.

\bibitem{SCH}  Seaman, B. T.,  Car, L. D. and Holland, M. J., \textit{Effect of a potential step or impurity on the Bose-Einstein condensate mean field},  Phys. Rev. A,  71, (2005).


\bibitem{SuSu}  Sulem, C. and Sulem, P-L., \textit{Nonlinear Schr\"odinger Equations: Self-Focusing and Wave Collapse}, Applied Mathematical Sciences, vol. 139, Springer, New York, (1999)


\bibitem{Tao}  Tao, T., \textit{Local And Global Analysis of Nonlinear Dispersive And Wave Equations}, CBMS Regional Conference Series in Mathematics, AMS, vol. 106,
Providence, RI., (2006)




\bibitem{w3} Weinstein, M.I., \textit{Nonlinear Schr\"odinger equation and sharp
interpolation estimates}. Comm. Math. Phys., 87, (1983), p.
567-576.


\bibitem{ZaSh} Zakharov, V. E., and Shabat, A. B., , \textit{Exact theory of two dimensional and one dimensional self modulation of waves in nonlinear media}. Sov. Phys. J.E.T.P. 34, (1972), p.
62-69.

\end{thebibliography}
\end{document}